\documentclass{m2an}
\usepackage{bm}
\usepackage{amssymb}
\usepackage{graphicx}
\usepackage[colorlinks = true,
            linkcolor = blue,
            urlcolor  = blue,
            citecolor = blue,
            anchorcolor = blue]{hyperref}
\usepackage{caption}
\usepackage[margin=1.2in]{geometry}
\usepackage{algorithm,algorithmic}
\newcommand{\dd}{\mathrm{d}}
\newcommand*\mathinhead[2]{\texorpdfstring{$\boldsymbol{#1}$}{#2}}
\allowdisplaybreaks
\begin{document}
\title{Discontinuous Galerkin discretization in time of systems of second-order nonlinear hyperbolic equations}
\author{Aili Shao}\address{Department of Mathematics, University of Oxford, Radcliffe Observatory Quarter, Woodstock Road, Oxford,UK OX2 6GG; \email: {aili.shao@maths.ox.ac.uk}}
\date{17 December 2021}

\begin{abstract} In this paper we study the finite element approximation of systems of second-order nonlinear hyperbolic equations. The proposed numerical method combines a $hp$-version discontinuous Galerkin finite element approximation in the time direction with an $H^1(\Omega)$-conforming finite element approximation in the spatial variables. Error bounds at the temporal nodal points are derived under a weak restriction on the temporal step size in terms of the spatial mesh size. Numerical experiments are presented to verify the theoretical results. \end{abstract}
%
%
\maketitle
\section{Introduction}
This paper aims to show how a discontinuous Galerkin time-stepping method can be used to approximate solutions of second-order quasilinear hyperbolic systems, which arise in a range of relevant applications, namely elastodynamics and general relativity. There has been a substantial body of research devoted to both the theoretical and numerical analysis of solutions of second-order hyperbolic equations. In particular, Kato \cite{Ka76} established the existence of solutions to the initial-boundary-value problem for quasilinear hyperbolic equations using semigroup theory. Building on Kato's work \cite{Ka76}, Hughes, Kato and Marsden \cite{HKM} analyzed the existence, uniqueness and well-posedness for a more general class of quasilinear second-order hyperbolic systems on a short time interval. They also applied these results to elastodynamics and Einstein's equations for the Lorentz metric $g_{\alpha, \beta}$ on $\mathbb{R}^4$, $0\leq\alpha,\beta\leq 3$. In contrast to the semigroup approach, Dafermos and Hrusa \cite{DH} used energy methods to establish local in time existence of smooth solutions to initial-boundary-value problems for such hyperbolic systems on a bounded domain $\Omega\subseteq\mathbb{R}^d$ where $d=1,2,3.$ In the case of $d=3$, Chen and Von Wahl proved an existence theorem for similar initial-boundary-value problems in \cite{CW}. Concerning numerical approximations of second-order hyperbolic equations, fully discrete schemes based on Galerkin finite element approximations in space for the linear case can be found in \cite{BDS, Ba84, BB}. These time-discrete schemes are generated from the rational approximations to either the cosine or the exponential and were later generalized by Bales and Dougalis \cite{BD,Ba} to approximate nonlinear hyperbolic problems. In \cite{Ba}, Bales considered a scalar nonlinear wave equation and introduced a class of single-step fully discrete schemes, which have temporal accuracy up to fourth-order. Dupont \cite{Du} and Dendy \cite{De} also showed optimal-order $L^2$ and $H^1$ error estimates for scalar nonlinear wave equations in semi-discrete Galerkin schemes. Makridakis \cite{Ma} proved optimal $L^2$ error estimates for both a semi-discrete and a class of fully discrete schemes for systems of second-order nonlinear hyperbolic equations. In \cite{OS}, Ortner and S\"{u}li developed the convergence analysis of semidiscrete discontinuous Galerkin finite element approximations of second-order quasilinear hyperbolic systems. Hockbruck and Maier\cite{HM} proved error estimates for space discretizations of a general class of first- and second-order quasilinear wave-type problems. In this paper, we will focus on the equations of nonlinear elastodynamics, though the ideas and techniques can be easily generalized to other second-order nonlinear hyperbolic equations, for instance, the Einstein's equations, provided the assumptions on the nonlinearity assumed herein are satisfied.

We begin by formulating the time-dependent problem resulting from nonlinear elasticity. Let $\Omega$ be a bounded domain in $\mathbb{R}^d$ for $d=1,2,3$, with sufficiently smooth boundary $\partial\Omega$, and let $0<T<\infty$. We consider the following initial-boundary-value problem:
\begin{equation}\label{setup}
\ddot u_i(x,t)-\sum_{\alpha=1}^d \partial_{\alpha} S_{i\alpha}(\nabla \mathbf{u}(x,t))=f_{i}(x,t) \mbox{ in } \Omega\times (0, T],
\end{equation}
for each $i=1,\ldots, d,$ where $\mathbf{u}=[u_1,\ldots, u_d]^{\mathrm{T}}$ represents the displacement field and $\mathbf{f}=[f_1,\ldots,f_d]^{\mathrm{T}}$ is the given body force which is sufficiently smooth, and 
\begin{equation}\label{bc}
\mathbf{u}(x,t)=\mathbf{0} \mbox{ on } \partial{\Omega}\times (0,T],
\end{equation}
\begin{equation}\label{ic}
\mathbf{u}(x,0)=\mathbf{u}_0(x)\in [H_0^1(\Omega)]^d\cap [H^m(\Omega)]^d, \quad \mathbf{\dot u}(x,0)=\mathbf{u}_1(x)\in[H^{m-1}(\Omega)]^d,
\end{equation}
are prescribed boundary and initial conditions, and $m$ is an integer which will be specified later. Here the dots over $\mathbf{u}$ denote differentiation with respect to time $t$, and $\partial_{\alpha}$ is the partial derivative with respect to $x_{\alpha}$. $S$ is a given smooth $d\times d$ matrix-valued function defined on $\mathbb{R}^{d\times d}$, which characterizes the Piola--Kirchhoff stress tensor. For a complete discussion of the relevant mechanical background, we refer the reader to \cite{An, Gu}.

For hyperelastic materials, $S$ is the gradient of a scalar-valued `stored energy function'. Hence, if $$A_{i\alpha j\beta}(\eta):=\frac{\partial}{\partial \eta_{j\beta}} S_{i\alpha}(\eta), \: \eta\in \mathbb{R}^{d\times d},$$ the \textit{elasticities} $A_{i\alpha j\beta}$ satisfy
\begin{equation}\tag{S1a}
A_{i\alpha j\beta}=A_{j\beta i\alpha}, 1\leq i, \alpha, j, \beta\leq d.
\end{equation}
We assume that $A_{i\alpha j\beta}$ satisfy the \emph{strict Legendre--Hadamard condition} 
\begin{equation}\tag{S1b}
\sum_{i, \alpha,\beta,j=1}^d A_{i\alpha j\beta}(\eta)\zeta_{\alpha}\zeta_{\beta}\xi_{i}\xi_{j}\geq M_0 |\zeta|^2|\xi|^2\quad \mbox{ for all } \eta \in \mathcal{O} \mbox{ and } \zeta, \xi\in\mathbb{R}^d,
\end{equation}
for some real number $M_0>0$, where $\mathcal{O}$ is the domain of definition of $A_{i\alpha j\beta}$ and here $|\cdot|$ denotes the Euclidean norm on $\mathbb{R}^d$. This condition (S1b) is indeed satisfied by the constitutive relations of the standard material models on sizeable portions of the displacement gradient space \cite{DH}.

The initial-boundary-value problem (\ref{setup})--(\ref{ic}) does not have a global smooth solution as a result of breaking waves and shocks no matter how smooth $\mathbf{u}_0$, $\mathbf{u}_1$ and $\mathbf{f}$ are. It was proved by Dafermos and Hrusa \cite{DH} that there exists a unique local solution to the problem (\ref{setup})--(\ref{ic}) provided that (S1a,b) are satisfied. We summarize this existence result in the following theorem.
\begin{thrm}\label{existence theorem}
Let $\Omega$ be a bounded domain in $\mathbb{R}^d$ with smooth boundary $\partial\Omega$. Assume that (S1a,b) hold, that $A_{i \alpha j\beta}$ and $\mathbf{f}$ are sufficiently smooth, and that $\mathbf{u}_0\in [H^m(\Omega)]^d$ and $\mathbf{u}_1\in [H^{m-1}(\Omega)]^d$ for some integer $m\geq[\frac{d}{2}]+3.$ Assume further that the initial values of the time derivatives of $\mathbf{u}$ up to order $m-1$ vanish on $\partial\Omega$ and that $\nabla \mathbf{u}_0[\overline{\Omega}]\subset\mathcal{O}.$ Then, there exists a finite time $T>0$ for which (\ref{setup})--(\ref{ic}) has a unique solution $\mathbf{u}$ such that 
\begin{equation}\label{regularity}
\mathbf{u}\in \bigcap_{s=0}^m C^{m-s}([0,T]; [H^s(\Omega)]^d).
\end{equation}
\end{thrm}  
By the Sobolev embedding theorem, (\ref{regularity}) implies that  $$\mathbf{u}\in [C^{\beta}( [0,T]\times \overline{\Omega})]^d:=\bigcap_{s=0}^{\beta} C^{\beta-s}([0,T]; C^s(\overline{\Omega})^d),$$
where $\beta=m-[\frac{d}{2}]-1$. Note that the assumption on $m$ implies that $\beta\geq 2.$ We shall assume throughout the paper that the above assumptions are satisfied for $m$ sufficiently large so that a unique solution of (\ref{setup}) exists. In fact, by Theorem \ref{existence theorem}, this unique solution satisfies  
$$\mathbf{u}\in C^1([0,T],[H^{m-1}(\Omega)]^d)\cap C([0,T],[H^m(\Omega)]^d).$$

There have been a few contributions to the literature discussing the numerical approximations of (\ref{setup})--(\ref{ic}). In particular, Makridakis \cite{Ma} proved an optimal $L^2$ error estimate for the classical conforming method under the assumption that $p>\frac{d}{2}$, where $p$ is the polynomial degree in space, and for appropriate initial approximations. In \cite{OS}, S\"{u}li and Ortner showed an optimal error estimate based on the broken $H^1$ norm for the semi-discrete discontinuous Galerkin finite element approximations of a similar problem, but with mixed Dirichlet-Neumann boundary conditions and weaker Lipschitz assumptions on the nonlinear term. Following \cite{Ba,BD}, Makridakis also proposed and analyzed two different fully discrete schemes to approximate (\ref{setup})--(\ref{ic}). The first scheme is based on second-order accurate approximations of the cosine while the second group of fully discrete methods that have temporal order of accuracy up to fourth-order are based on rational approximations of the exponential function. Discontinuous Galerkin in time methods, which are the focus of this paper are, however, arbitrarily high-order accurate. In contrast with the above-mentioned finite difference time integration schemes, for which the solution at the current time step depends on the previous steps, this discontinuous-in-time scheme on the time interval $(t_n, t_{n+1}]$ only depends on the solution at $t_n^{-}$. Since the local polynomial degree is free to vary between time steps, this method is also naturally suited for an adaptive choice of the time discretization parameters. To the best of our knowledge, the analysis of 
a fully discrete scheme based on such  discontinuous-in-time discretization of  second-order quasilinear hyperbolic systems has not been previously considered in the literature. 

Discontinuous Galerkin methods \cite{RH, LR} have been widely and successfully used for the numerical approximations of PDEs. They have been first introduced by Reed and Hill \cite{RH} to solve the hyperbolic neutron transport equations. Simultaneously, but independently, they were proposed as non-standard numerical schemes for solving elliptic and parabolic problems by Babu\v{s}ka $\&$ Zl\'{a}mal \cite{BZ}, Baker \cite{Bak}, Wheeler \cite{Wh}, Arnold \cite{Ar} and Rivi\`{e}re \cite{Ri} etc. In recent years, there has been considerable interest in applying discontinuous Galerkin finite element methods to nonlinear hyperbolic PDEs. In particular, Antonietti \emph{et al.} \cite{AMMNW} developed a high-order discontinuous Galerkin scheme for the spatial discretization of nonlinear acoustic waves. Muhr \emph{et al.}\cite{MWN} also proposed and analyzed a hybrid discontinuous Galerkin coupling approach for the semi-discrete nonlinear elasto-acoustic problem. However, there has been little work on the construction and mathematical analysis of fully discrete discontinuous Galerkin schemes for second-order nonlinear hyperbolic PDEs. In this article, a high-order discontinuous Galerkin finite element method for the time integration will be proposed and analyzed. To construct such a scheme, we first discretize with respect to the spatial variables by the means of a Galerkin finite element method, which results in a system of ordinary differential equations (ODEs) in time; then we discretize the resulting ODE system using discontinuous Galerkin method in time  ( e.g. see Antonietti \emph{et al.} \cite{AMSQ}). The resulting weak formulation in time is based on weakly imposing the continuity of the approximate displacements and velocities between time steps by penalizing jumps in these quantities in the definition of the numerical method. 

The paper is structured as follows. The next section sets up the assumptions required for the numerical approximation. Section \ref{sect_scheme} dicusses the construction of a fully discrete scheme for the approximations of (\ref{setup})--(\ref{ic}) using a time-discontinuous Galerkin method. In Section \ref{elasto_convergence}, we perform the convergence analysis of the discontinuous-in-time scheme under the hypotheses (S2a,b). Building on the work of Makridakis \cite{Ma}, this convergence proof is based on Banach's fixed point theorem and a nonlinear elliptic projection operator whose approximation properties will be analyzed in Section \ref{W}. Finally, numerical experiments are presented in Section \ref{numerical} to verify the theoretical results.
\section{Definition and assumptions}
In order to find a numerical approximation to the solution of the hyperbolic system (\ref{setup})--(\ref{ic}), we discretize it in space using a continuous Galerkin method, and then apply a discontinuous Galerkin method in time. For the sake of showing the well-posedness of the resulting numerical method, we consider the substitution $\mathbf{u}=\mathrm{e}^{\gamma t}\mathbf{v}$, with $\gamma>0$, resulting in the equivalent equation:
\begin{equation}\label{setup2}
\ddot v_i(x,t)+2\gamma \dot v_i(x,t)+\gamma^2 v_i(x,t)-\mathrm{e}^{-\gamma t}\sum_{\alpha=1}^d \partial_{\alpha} S_{i\alpha}(\mathrm{e}^{\gamma t}\nabla \mathbf{v}(x,t))=\tilde{f}_{i}(x,t) \mbox{ in } \Omega\times (0, T],
\end{equation}
for each $i=1,\ldots, d,$ where $\mathbf{\tilde{f}}=\mathrm{e}^{-\gamma t} \mathbf{f}$, $\gamma>0$ is a fixed constant, 
\begin{equation}\label{bc2}
\mathbf{v}(x,t)= \mathbf{0} \mbox{ on } \partial{\Omega}\times (0,T],
\end{equation}
\begin{equation}\label{ic2}
\mathbf{v}(x,0)=\mathbf{v}_0(x)\in [H^{m}(\Omega)]^d \cap [H_0^1(\Omega)]^d, \quad \mathbf{\dot v}(x,0)=\mathbf{v}_1(x)\in [H^{m-1}(\Omega)]^d,
\end{equation}
where $\mathbf{v}_0(x)=\mathbf{u}_0(x)$ and $\mathbf{v}_1(x)=\mathbf{u}_1(x)-\gamma\mathbf{u}_0(x).$
From now on, we focus on the system of equations (\ref{setup2})--(\ref{ic2}) only.  By Theorem \ref{existence theorem}, we have
$$\mathbf{v}=\mathbf{u}\mathrm{e}^{-\gamma t}\in C^1([0,T],[H^{m-1}(\Omega)]^d)\cap C([0,T],[H^m(\Omega)]^d).$$
This shows that the initial conditions stated as (\ref{ic2}) in the above definition are meaningful.

Before describing its discretization, we first fix the notation. We use the symbol $:=$ to indicate an equality in which the left hand side is defined by the right hand side. We denote by $(\cdot, \cdot)_{L^2}$ the inner product in $L^2(\Omega)$ and $[L^2(\Omega)]^d$. Following standard notational conventions, we shall write $W^{s,p}:=[W^{s,p}(\Omega)]^d$ for $s\in\mathbb{Z}, p\in\mathbb{R}^{+}\cup\{\infty\}$, and put $H^s:=W^{s,2}$. Similarly, $H_0^1:=[H_0^1(\Omega)]^d$ and $L^2:=[L^2(\Omega)]^d.$

We define the following time-dependent semilinear form
$$a(\mathbf{v}(t),\bm{\varphi}):=\sum_{i,\alpha=1}^d \mathrm{e}^{-\gamma t}(S_{i\alpha}( e^{\gamma t}\nabla \mathbf{v}(t)),\partial_{\alpha}\bm{\varphi}_i)_{L^2}, \quad \text{ for }\bm{\varphi}\in H_0^1.$$
Since we also need to approximate the gradient of the solution $\nabla \mathbf{u}=\mathrm{e}^{\gamma t}\nabla\mathbf{v}$, we assume that there exists an open convex set $\mathcal{M}$ with $\overline{\mathcal{M}}\subset \mathcal{O}$ such that $\nabla \mathbf{u}( [\overline{\Omega}\times [0,T]])\subset \mathcal{M}$. If the distance of $\overline{\mathcal{M}}$ from $\partial\mathcal{O}$ is $\delta$, we consider the set 
\begin{equation}\label{M}
\mathcal{M}_{\delta}:=\{\eta\in\mathbb{R}^{d\times d}\colon \inf_{\sigma\in \mathcal{M}} |\eta-\sigma|\leq \delta \},
\end{equation}
where $|\cdot|$ denotes the Frobenius norm on $\mathbb{R}^{d\times d}$ defined, for $\eta\in\mathbb{R}^{d\times d}$, by $|\eta|=(\eta:\eta)^{\frac{1}{2}}.$ Notice that the set $\mathcal{M}_{\delta}$ is convex (cf. Lemma 1 in \cite{OS2}). Since we only require $S$ to be locally Lipschitz continuous in $\mathcal{M}_{\delta}$, we define the local Lipschitz constant of $S$ in $\mathcal{M}_{\delta}$ by 
\begin{equation}\label{lipschitz_const}
K_{\delta}:=\sup_{\eta\in \mathcal{M}_{\delta}}\left(\sum_{i,\alpha,j,\beta=1}^d |A_{i\alpha j\beta}(\eta)|^2 \right)^{\frac{1}{2}},
\end{equation}
and the local Lipschitz constant of the fourth-order elasticity tensor $A=\nabla S$ by
\begin{equation}
L_{\delta}:=\sup_{\eta_1,\eta_2\in\mathcal{M}_{\delta}, \:\eta_1\neq \eta_2} |\eta_1-\eta_2|^{-1}\left( \sum_{i,\alpha,j,\beta=1}^d |A_{i\alpha j\beta}(\eta_1)-A_{i\alpha j\beta}(\eta_2)|^2\right)^{\frac{1}{2}}.
\end{equation}
Since the set $\mathcal{M}_{\delta}$ is a compact subset of $\mathbb{R}^{d\times d}$ for every $\delta>0$ and $A_{i\alpha j\beta}$ is sufficiently smooth (and in particular continuously differentiable on $\mathcal{M}_{\delta})$, it follows that $K_{\delta}$ and $L_{\delta}$ are finite. 
We also define
\begin{equation}\label{Z}
\mathcal{Z}_{\delta}:=\{\Phi\in L^{\infty}(\overline{\Omega})^{d\times d}\colon \Phi(x)\in \mathcal{M}_{\delta}, x\in\overline{\Omega}\}.
\end{equation} 
This set $\mathcal{Z}_{\delta}$ is expected to contain the gradients of approximations of $\mathbf{u}$. We define 
\begin{equation}
\tilde{a}(\bm{\varphi};\bm{\phi},\bm{\psi}):=\sum_{i,\alpha,j,\beta=1}^d \left( A_{i\alpha j\beta} (\nabla \bm{\varphi})\partial_{\beta}\bm{\phi}_j,\partial_{\alpha}\bm{\psi}_i\right)_{L^2}, \quad \mbox{ for } \bm{\varphi},\bm{\phi},\bm{\psi}\in H_0^1.
\end{equation}
By the definition of $\mathcal{Z}_{\delta}$ and (S1a), we have
\begin{equation}\tag{S2a}
\tilde{a}(\bm{\varphi};\bm{\phi},\bm{\psi})=\tilde{a}(\bm{\varphi};\bm{\psi},\bm{\phi}), \quad\mbox{ for } \bm{\varphi},\bm{\phi},\bm{\psi} \in H_0^1, \nabla\bm{\varphi}\in \mathcal{Z}_{\delta}.
\end{equation}
We also assume that there exists a real number $M_1>0$ such that 
\begin{equation}\tag{S2b}
\tilde{a}(\bm{\varphi};\bm{\phi},\bm{\phi})\geq M_1 \|\nabla \bm{\phi}\|_{L^2}^2, \quad\mbox{ for } \bm{\varphi},\bm{\phi} \in H_0^1, \nabla\bm{\varphi}\in \mathcal{Z}_{\delta}.
\end{equation}
Note that (S2b) is a stronger assumption than (S1b). In general, (S1b) does not imply (S2b) for $d>1$. We refer the reader to \cite{Zh, Sv, DM} for counterexamples. In fact, (S1b) only implies the following G\r{a}rding's inequality:
\begin{equation}\label{garding}
\tilde{a}(\bm{\varphi};\bm{\phi},\bm{\phi})\geq \frac{1}{2}M_0 \|\nabla \bm{\phi}\|_{L^2}^2-\mu\| \bm{\phi}\|_{L^2}^2 \quad \mbox{ for } \mu \geq 0, \bm{\varphi}, \bm{\phi}\in H_0^1, \nabla\bm{\varphi}\in\mathcal{Z}_{\delta},
\end{equation}
cf. Theorem 6.5.1 in \cite{Mo} and  Lemma 5 in \cite{OS2}. We note that the techniques of this paper can be extended so that our results are still valid under this weak condition.
\section{Numerical scheme}\label{sect_scheme}
\subsection{Semi-discrete approximation}
We shall discretize the problem (\ref{setup2})--(\ref{ic2}) in space using a continuous Galerkin method. For the spatial discretization parameter $h\in(0,1)$, we define $\mathcal{V}_h$ to be a given family of finite-dimensional subspaces of $H_0^1\cap H^m$ with polynomial degree $p\geq 1$. We shall assume that the triangulation $\{\mathcal{T}_h\}_{h>0}$ of $\Omega$ into $d$-dimensional simplices, which are possibly curved along the boundary $\partial\Omega$, is shape-regular and quasi-uniform.  It follows from Bernardi's work \cite{Be} that 
\begin{equation}\tag{i}
\inf_{\mathbf{v}_h\in \mathcal{V}_h} \{\| \mathbf{v}-\mathbf{v}_h \|_{L^2}+h\|\mathbf{v}-\mathbf{v}_h\|_{H^1}\}\leq C h^{r+1}\|\mathbf{v}\|_{H^{r+1}}, \quad 1\leq r\leq \min(p, m-1), \quad \mathbf{v}\in H^m\cap H_0^1.
\end{equation}
Further, the following inverse inequalities follow directly from the quasi-uniformity of the triangulation.

There exists a positive constant $C_0$ such that, for every $\mathbf{v}_h\in \mathcal{V}_h$,
\begin{equation}\tag{ii,a}
\|\nabla \mathbf{v}_h\|_{L^2}\leq C_0 h^{-1}\|\mathbf{v}_h\|_{L^2} \quad\mbox{ and }\quad \|\nabla \mathbf{v}_h\|_{L^{\infty}}\leq C_0h^{-1}\|\mathbf{v}_h\|_{L^{\infty}}.
\end{equation}
There exists a positive constant $C_1$ such that, for every $\mathbf{v}_h\in \mathcal{V}_h$,  
\begin{equation}\tag{ii,b}
\|\nabla \mathbf{v}_h\|_{L^{\infty}}\leq C_1 h^{-\frac{d}{2}}\|\nabla \mathbf{v}_h\|_{L^2} \quad\mbox{ and }\quad \|\mathbf{v}_h\|_{L^{\infty}}\leq C_1 h^{-\frac{d}{2}}\|\mathbf{v}_h\|_{L^2}.
\end{equation}
With the above assumptions, we are ready to construct the continuous-in-time finite element approximation $\mathbf{v}_h$ of $\mathbf{v}$. The semi-discrete approximation $\mathbf{v}_h\colon [0,T]\to \mathcal{V}_h$ of the solution of (\ref{setup2})--(\ref{ic2}) satisfies the following initial-value problem in $\mathcal{V}_h$:
\begin{equation}\label{semi-discrete}
\left(\mathbf{\ddot{v}}_h(t), \bm{\varphi}\right)_{L^2}+a(\mathbf{v}_h(t),\bm{\varphi})+ 2\gamma \left(\mathbf{\dot{v}}_h(t), \bm{\varphi}\right)_{L^2} +\gamma^2 \left(\mathbf{v}_h(t), \bm{\varphi} \right)_{L^2}= (\tilde{\mathbf{f}}(t),\bm{\varphi})_{L^2}
\end{equation}
for all $\bm{\varphi}\in \mathcal{V}_h, 0\leq t\leq T,$
\begin{equation}\label{semi-ic}
\mathbf{v}_h(0)=\mathbf{v}_{0,h}\in \mathcal{V}_h, \quad \mathbf{\dot v}_h(0)=\mathbf{v}_{1,h}\in \mathcal{V}_h,
\end{equation}
where $\mathbf{v}_{0,h}$ and $\mathbf{v}_{1,h}$ are specially chosen initial values.
It was proved by Makridakis \cite{Ma} that the semi-discrete problem (\ref{semi-discrete}), (\ref{semi-ic}) with $\gamma=0$ (the semi-discrete form based on a continuous finite element approximation of the original problem (\ref{setup})--(\ref{ic})) has a locally unique solution and that the optimal-order $L^2$ error estimate 
\begin{equation}\label{L^2 error estimate}
\max_{0\leq t\leq T}\|\mathbf{v}(t)-\mathbf{v}_h(t)\|_{L^2}\leq C(\mathbf{v}) h^{p+1}
\end{equation}
holds for sufficiently smooth initial data. Here $p$ is the polynomial degree of the elements of the finite-dimensional space $\mathcal{V}_h$, which satisfies $p>\frac{d}{2}$. The proofs of these assertions for $\gamma>0$ are completely analogous and are therefore omitted.
\subsection{Discontinuous-in-time fully discrete scheme}\label{dgintime}
In this section we shall construct a fully discrete approximation of the solution of (\ref{setup2})--(\ref{ic2}) by applying a discontinuous Galerkin method in time. For this purpose, we partition the time interval $I=(0,T]$ into $N$ sub-intervals $I_n=(t_{n-1}, t_n]$ having length $k_n=t_n-t_{n-1}$ for $n=1,2,\ldots, N$, with $t_0=0$ and $t_N=T.$ To deal with the discontinuity at each $t_n$ in the numerical approximation to $\mathbf{v}$, we introduce the jump operator
$$[\mathbf{v}_h]_n:=\mathbf{v}_h(t_n^{+})-\mathbf{v}_h(t_n^{-})\quad \mbox{ for } n=0,1,\ldots, N-1,$$
where 
$$\mathbf{v}_h(t_n^{\pm}):=\lim_{\varepsilon\to 0^{\pm}}\mathbf{v}_h(t_n+\varepsilon) \quad\mbox{ for } n=0,1,\ldots, N-1.$$
By convention, we assume that $\mathbf{v}_h(0^{-})=\mathbf{v}_{0,h}$ and $\mathbf{\dot v}_h(0^{-})=\mathbf{v}_{1,h}$. 
Moreover, we define $\mathbf{v}_{h,n}^{+}:=\mathbf{v}_h(t_n^{+})$ and $\mathbf{v}_{h,n}^{-}:=\mathbf{v}_h(t_n^{-}).$ 
To deal with the nonlinear term, we apply Taylor's theorem with an integral remainder to have
\begin{align*}
S_{i\alpha}(\nabla \mathbf{v}_h(t)\mathrm{e}^{\gamma t})&=S_{i\alpha}(\mathbf{0})+\sum_{j,\beta=1}^d \partial_{\beta} \mathbf{v}_{h,j}(t)\mathrm{e}^{\gamma t}\int_0^1 A_{i\alpha j\beta}(\tau \nabla \mathbf{v}_h(t)\mathrm{e}^{\gamma t})\,\dd \tau. 
\end{align*}
By assuming that $S(\mathbf{0})=\mathbf{0}$, we can write the semilinear term as 
\begin{flalign*}
&\int_{t_{n-1}}^{t_n} \sum_{i,\alpha=1}^d \mathrm{e}^{-\gamma t}\left(S_{i\alpha}(\nabla \mathbf{v}_h(t)\mathrm{e}^{\gamma t}), \partial_{\alpha} \mathbf{\dot v}_{h,i}(t)\right)_{L^2}\, \dd t &&\\
&=\int_{t_{n-1}}^{t_n}\sum_{i,\alpha, j,\beta=1}^d \left( \partial_{\beta}\mathbf{v}_{h,j}(t)\int_0^1 A_{i\alpha j\beta}(\tau\nabla \mathbf{v}_h(t) \mathrm{e}^{\gamma t})\,\dd \tau, \partial_{\alpha}\mathbf{\dot v}_{h,i}(t)\right)_{L^2} \,\dd t &&\\
&=\int_{t_{n-1}}^{t_n}\int_0^1 \sum_{i,\alpha,j,\beta=1}^d\left( A_{i\alpha j\beta}(\tau \nabla \mathbf{v}_h(t)\mathrm{e}^{\gamma t})\partial_{\beta}\mathbf{v}_{h,j}(t), \partial_{\alpha}\mathbf{\dot v}_{h,i}(t)\right)_{L^2} \,\dd \tau\,\dd t &&\\
&=\int_{t_{n-1}}^{t_n} \bigg\{ \int_0^1 \tilde{a}(\tau \mathbf{v}_h(t) \mathrm{e}^{\gamma t}; \mathbf{v}_h(t), \mathbf{\dot v}_h(t)) \,\dd \tau \bigg\}\,\dd t,
\end{flalign*}
where 
$$\tilde{a}(\tau \mathbf{v}_h(t) \mathrm{e}^{\gamma t}; \mathbf{v}_h(t), \mathbf{\dot v}_h(t)):=\sum_{i,\alpha,j,\beta=1}^d\left( A_{i\alpha j\beta}(\tau \nabla \mathbf{v}_h(t)\mathrm{e}^{\gamma t})\partial_{\beta}\mathbf{v}_{h,j}(t), \partial_{\alpha}\mathbf{\dot v}_{h,i}(t)\right)_{L^2}.$$
We focus on the generic time interval $I_n$ and assume that the solution on $I_{n-1}$ is known. Following the discontinuous-in-time numerical scheme introduced in \cite{AMSQ}, we first test the equation (\ref{setup2}) against $\bm{\dot \varphi}$ for $\bm{\varphi}\in H^1(I_n; H_0^1)$ and integrate on $I_n$ to obtain the following weak formulation:
\begin{flalign}\label{weak_form}
&\int_{t_{n-1}}^{t_n} \left( \mathbf{\ddot v}(t), \bm{\dot\varphi}(t) \right)_{L^2}\: \dd t+\int_{t_{n-1}}^{t_n}\sum_{i,\alpha=1}^d \mathrm{e}^{-\gamma t}\left( S_{i\alpha}(\nabla \mathbf{v}(t)\mathrm{e}^{\gamma t}),\partial_{\alpha}\bm{\dot{\varphi}}_i(t) \right)_{L^2}\: \dd t+2\gamma \int_{t_{n-1}}^{t_n} \left(\mathbf{\dot v}(t), \bm{\dot\varphi}(t)\right)_{L^2}\,\dd t&& \nonumber\\ 
&\quad+\gamma^2 \int_{t_{n-1}}^{t_n} \left(\mathbf{v}(t), \bm{\dot \varphi}(t)\right)_{L^2} \dd t=\int_{t_{n-1}}^{t_n} \left(\tilde{\mathbf{f}}(t), \bm{\dot{\varphi}}(t)\right)_{L^2} \dd t.&&
\end{flalign}
Now we rewrite (\ref{weak_form}) by adding suitable (strongly consistent) terms:
\begin{flalign}\label{weak_form2}
&\int_{t_{n-1}}^{t_n} \left( \mathbf{\ddot v}(t), \bm{\dot \varphi}(t) \right)_{L^2}\dd t+\left( [\mathbf{\dot v}(t)]_{n-1}, \bm{\dot \varphi}(t_{n-1}^{+})\right)_{L^2}+\int_{t_{n-1}}^{t_n}\sum_{i,\alpha=1}^d e^{-\gamma t}\left( S_{i\alpha}(\nabla \mathbf{v}(t)e^{\gamma t}),\partial_{\alpha}\bm{\dot{\varphi}}_i(t) \right)_{L^2} \dd t &&\nonumber\\
&\quad+\frac{1}{2}\bigg\{\int_0^1 \tilde{a}(\tau \mathbf{v}_{h,n-1}^{+}\mathrm{e}^{\gamma t_{n-1}}; \mathbf{v}_{h,n-1}^{+}, \mathbf{v}_{h,{n-1}}^{+})\:\dd \tau -\int_0^1 \tilde{a}(\tau \mathbf{v}_{h,n-1}^{-}\mathrm{e}^{\gamma t_{n-1}}; \mathbf{v}_{h,n-1}^{-}, \mathbf{v}_{h,n-1}^{-})\:\dd \tau\bigg\}&&\nonumber\\
&\quad +2\gamma \int_{t_{n-1}}^{t_n} \left(\mathbf{\dot v}(t), \bm{\dot \varphi}(t)\right)_{L^2}\dd t+\gamma^2 \int_{t_{n-1}}^{t_n} \left(\mathbf{v}(t), \bm{\dot \varphi}(t)\right)_{L^2} \dd t+\gamma^2\left( [\mathbf{v}(t)]_{n-1}, \bm{\varphi}_{n-1}^{+}\right)_{L^2}&&\nonumber\\
{}&=\int_{t_{n-1}}^{t_n} \left(\tilde{\mathbf{f}}(t), \bm{\dot{\varphi}}(t)\right)_{L^2} \,\dd t.
\end{flalign}
Summing over all time intervals in (\ref{weak_form2}) leads us to define  the following semilinear form $\mathcal{A}\colon \mathcal{H}\times \mathcal{H}\to\mathbb{R}$ with $\mathcal{H}:= H^2(0,T; H^m\cap H_0^1)$
by
\begin{flalign*}
\mathcal{A}(\mathbf{v},\bm{\varphi})&:=\sum_{n=1}^N \int_{t_{n-1}}^{t_n} \left( \mathbf{\ddot v}(t), \bm{\dot \varphi}(t) \right)_{L^2} \dd t+\sum_{n=1}^{N-1} \left([\mathbf{\dot v}(t) ]_n, \bm{\dot \varphi}(t_n^{+}) \right)_{L^2}+\left(\mathbf{\dot v}(t_0^{+}), \bm{\dot \varphi}(t_0^{+}) \right)_{L^2} &&\\
&\quad+\sum_{n=1}^N \int_{t_{n-1}}^{t_n}\sum_{i,\alpha=1}^d \mathrm{e}^{-\gamma t}\left( S_{i\alpha}(\nabla \mathbf{v}(t)\mathrm{e}^{\gamma t}),\partial_{\alpha}\bm{\dot \varphi}_i(t) \right)_{L^2} \: \dd t+\frac{1}{2}\int_0^1 \tilde{a}(\tau \mathbf{v}_0^{+}; \mathbf{v}_0^{+}, \bm{\varphi}_0^{+})\:\dd \tau  &&\nonumber\\
&\quad+\sum_{n=1}^{N-1} \frac{1}{2}\bigg\{\int_0^1 \tilde{a}(\tau \mathbf{v}_{n}^{+}\mathrm{e}^{\gamma t_n}; \mathbf{v}_{n}^{+}, \bm{\varphi}(t_n^{+}))\:\dd \tau -\int_0^1 \tilde{a}(\tau \mathbf{v}_{n}^{-}\mathrm{e}^{\gamma t_n}; \mathbf{v}_{n}^{-}, \bm{\varphi}(t_n^{-}))\:\dd \tau\bigg\} &&\nonumber\\
&\quad+2\gamma\sum_{n=1}^N \int_{t_{n-1}}^{t_n} \left( \mathbf{\dot v}(t), \bm{\dot \varphi}(t) \right)_{L^2} \,\dd t+\gamma^2 \sum_{n=1}^N \int_{t_{n-1}}^{t_n} \left( \mathbf{ v}(t), \bm{\dot \varphi}(t)\right)_{L^2}\: \dd t +\gamma^2\sum_{n=1}^{N-1} \left([\mathbf{v}(t)]_n,  \bm{\varphi}_n^{+} \right)_{L^2} &&\nonumber\\
&\quad +\gamma^2\left(\mathbf{v}_0^{+},  \bm{\varphi}_0^{+} \right)_{L^2} 
\end{flalign*}
for $\bm{\varphi}\in\mathcal{H}.$ Let, further, $F$ be the linear functional defined by
\begin{flalign*}
F(\bm{\varphi}):=&\sum_{n=1}^N \int_{t_{n-1}}^{t_n} \left(\tilde{\mathbf{f}}(t), \bm{\dot \varphi}(t)\right)_{L^2} \dd t+\left(\mathbf{v}_1, \bm{\dot \varphi}_0^{+}\right)_{L^2}+\gamma^2\left(\mathbf{v_0},\bm{\varphi}_0^{+}\right)_{L^2} +\frac{1}{2}\int_0^1 \tilde{a}(\tau \mathbf{v_0}; \mathbf{v_0}, \bm{\varphi}_0^{-})\:\dd \tau.&&
\end{flalign*}
Now we introduce the finite-dimensional space
$$\mathcal{V}_{kh}^{q_n}:=\{\mathbf{v}\colon [0,T]\to \mathcal{V}_h;\mathbf{v}\vert_{I_n}=\sum_{j=0}^{q_n} \mathbf{v}_j t^j, \mathbf{v}_j\in \mathcal{V}_h\},$$
with $q_n\geq 2$ for each $1\leq n\leq N.$ For $\mathbf{q}:=[q_1,q_2,\ldots,q_N]^T\in\mathbb{N}^N$, we then define the space  
$$\mathcal{V}_{kh}^{\mathbf{q}}:=\{\mathbf{v}\colon [0,T]\to \mathcal{V}_h; \mathbf{v}\vert_{I_n}\in \mathcal{V}_{kh}^{q_n} \mbox{ for } n=1,2,\ldots, N \}.$$  Then, the discontinuous-in-time fully discrete approximation of the problem reads as follows: find $\mathbf{v}_{\mathrm{DG}}\in\mathcal{V}_{kh}^{\mathbf{q}}$ such that 
\begin{equation}\label{discrete_vform}
\mathcal{A}(\mathbf{v}_{\mathrm{DG}}, \bm{\varphi})=\tilde{F}(\bm{\varphi})\quad \mbox{ for all } \bm{\varphi}\in \mathcal{V}_{kh}^{\mathbf{q}},
\end{equation}
where $\tilde{F}$ is a modified version of $F$ defined as 
\begin{align*}
\tilde{F}(\bm{\varphi}):=&\sum_{n=1}^N \int_{t_{n-1}}^{t_n} \left(\tilde{\mathbf{f}}(t), \bm{\dot \varphi}(t)\right)_{L^2} \,\dd t+\left(\mathbf{v}_{1,h}, \bm{\dot \varphi}_0^{+}\right)_{L^2}+\gamma^2\left(\mathbf{v_0}_h,\bm{\varphi}_0^{+}\right)_{L^2} +\frac{1}{2}\int_0^1 \tilde{a}(\tau \mathbf{v_0}_{h}; \mathbf{v_0}_{h}, \bm{\varphi}_0^{-})\:\dd \tau.
\end{align*}
\section{Convergence analysis}\label{elasto_convergence}
By using the ideas introduced in \cite{Ma} based on Banach's fixed point theorem, we will show the existence and uniqueness of $\mathbf{v}_{\mathrm{DG}}$. We shall also prove \emph{a priori} error estimates as summarized in the following theorem.
\begin{thrm}\label{convergence_thm}
Let $\mathbf{v}\in W^{s,\infty}([0,T];H^{m}\cap H_0^1)$, with $\frac{d}{2}+1< r\leq\min(p,m-1)$, $s\geq q_i+1$ for each $i=1,2,\ldots N$, be the solution of (\ref{setup2})--(\ref{ic2}). Assume that $k_i^{q_i-\frac{1}{2}}=o(h^{\frac{d}{2}+1})$ and there exist positive constants $\mu_i$, $\nu_i$ such that  $\mu_i k_i\leq h^2\leq \nu_i k_i$ for each $i=1,2,\ldots N$. Suppose that we choose the initial data $\mathbf{v}_{0,h}$, $\mathbf{v}_{1,h}\in\mathcal{V}_{h}$ to be
\begin{equation}
\mathbf{v}_{1,h}=\mathbf{W}(0), \quad \mathbf{v}_{1,h}=\mathbf{\dot W}(0),
\end{equation}
where $\mathbf{W}(t)\in\mathcal{V}_h$ is the nonlinear elliptic projection of $\mathbf{v}(t)$ such that 
\begin{equation}\label{nonlinear projection}
a(\mathbf{W}(t), \bm{\varphi})=a(\mathbf{v}(t),\bm{\varphi})\quad \mbox{ for all } \bm{\varphi}\in \mathcal{V}_h.
\end{equation}
Then we have for the solutions of (\ref{discrete_vform}) that
\begin{equation}\label{error}
\|\mathbf{v}_{\mathrm{DG}}(t_j^{-})-\mathbf{v}(t_j^{-})\|_{L^2}+\|\mathbf{\dot v}_{\mathrm{DG}}(t_j^{-})-\mathbf{\dot v}(t_j^{-})\|_{L^2}\leq C(\mathbf{v})\left( h^{2r+2}+ \sum_{n=1}^j\frac{k_n^{2q_n+1}}{q_n^{2(s-1)}} \right)^{\frac{1}{2}} 
\end{equation}
for each  $j=1,\ldots,N,$ where $C(\mathbf{v})$ is a positive constant depending on the solution $\mathbf{v}$. 
\end{thrm}
\begin{rmrk}\label{uniform_grid}
If we use uniform time intervals $k_n=k=h^2$, and uniform polynomial degrees $q_n=q\geq 2$, for $n=1,\ldots, N$, then the error bound at the end nodal point  becomes 
$$\|\mathbf{v}_{\mathrm{DG}}(t_N^{-})-\mathbf{v}(t_N^{-})\|_{L^2}+\|\mathbf{\dot{v}}_{\mathrm{DG}}(t_N^{-})-\mathbf{\dot{v}}(t_N^{-})\|_{L^2}\leq C(\mathbf{v}) \left( k^{\frac{r+1}{2}}+ k^{q}\right).$$
\end{rmrk}
\begin{rmrk}
The assumptions that $k_i^{q_i-\frac{1}{2}}=o(h^{1+\frac{d}{2}})$ and $\mu_i k_i\leq h^2\leq \nu_i k_i$ for each $i=1,\ldots N$ require that $q_i>1+\frac{d}{4}$ for each $i=1,\ldots,N.$ That is, we need the polynomial degree in time satisfies $q_i\geq 2$ for $d=1,2,3$ on each time interval $I_i$, with $i=1,\ldots, N.$ 
\end{rmrk}
\begin{rmrk}\label{reg_v}
By the Sobolev embedding theorem, $\mathbf{v}\in W^{s,\infty}([0,T]; H^{m})$ for $m>\frac{d}{2}+2$ implies that $\mathbf{v}\in W^{s,\infty}([0,T]; C^{2,\alpha}(\overline{\Omega})^d)$ for some $\alpha\in (0,1).$ Note that the assumption $m>\frac{d}{2}+2$ is consistent with the assumption $m\geq [\frac{d}{2}]+3$ in Theorem \ref{existence theorem}. That is, we need $m\geq 3$ for $d=1$ and $m\geq 4$ for $d=2,3$.
\end{rmrk}
It will be assumed throughout the convergence analysis that 
$$\mathbf{v}\in W^{s,\infty}([0,T]; H^{m}\cap H_0^1).$$ 
\subsection{Definition of the fixed point map}\label{fixpointsetup}
It is known, see \cite{DR} and \cite{Ra}, that (\ref{nonlinear projection}) has, for $h$ sufficiently small, a locally unique solution $\mathbf{W}(t)\in\mathcal{V}_h$ for $0\leq t\leq T$. Furthermore, $\mathbf{W}$ satisfies the following properties, which are established in Section \ref{W}.

There exists a constant $C_r(\mathbf{v})$ depending on $\mathbf{v}$ such that, for $\frac{d}{2}+1<r \leq \min(p,m-1)$,
\begin{equation}\tag{iii,a}
\|\nabla \mathbf{v}^{(j)}(t)-\nabla\mathbf{W}^{(j)}(t)\|_{L^2}\leq C_r(\mathbf{v})h^{r}, \quad 0\leq t\leq T,
\end{equation}
for $j=0, 1$, where $\mathbf{v}^{(j)}:=\frac{\dd^{j}\mathbf{v}}{\dd t^{j}}.$
In addition, we shall prove for the time-derivatives of $\mathbf{W}$ there holds
\begin{equation}\tag{iii,b}
\|\mathbf{v}^{(j)}(t)-\mathbf{W}^{(j)}(t)\|_{L^2}\leq \tilde{C}_{r}(\mathbf{v})h^{r+1}, \quad 0\leq t\leq T, 
\end{equation}
for $j=0, 1,2$. We can also show that there exist constants $c_{0}$ and $c_1$, independent of $h$, such that
\begin{equation}\tag{iii,c}
\| \nabla \mathbf{W}(t)\|_{L^{\infty}}\leq c_0 \mbox{ and } \| \nabla \mathbf{\dot W}(t)\|_{L^{\infty}}\leq c_1, \quad 0\leq t\leq T.
\end{equation}
Let $\Pi_k=\Pi_{I_n}^{q_n}$ denote the modified $L^2$-projector in time direction. That is, for each $n=1,2,\ldots,N$,
\begin{align}
&(\Pi_k\mathbf{W}-\mathbf{W})(x,t_{n-1}^{+})=0;\label{prop1}\\
&(\Pi_k\mathbf{W}-\mathbf{W})(x,t_{n}^{-})=0;&&\label{prop12}\\
&\partial_t(\Pi_k\mathbf{W}-\mathbf{W})(x,t_{n}^{-})=0;\label{prop21}\\
&\int_{t_{n-1}}^{t_n} \left(\partial_t(\Pi_k\mathbf{W}-\mathbf{W}), \bm{\chi}\right)_{L^2}\, \dd t= 0 \quad\text{ for } \bm{\chi}\in \mathcal{V}_{kh}^{q_n-2}.\label{prop2}
\end{align}
It was first proved in \cite{SS} and further studied in \cite{Sh} that for each $\mathbf{W}\in H^s(I_n;L^2)$, there exists a positive constant $C$ such that, 
\begin{equation}\label{prop5}
\int_{t_{n-1}}^{t_n} \| \partial_{tt}(\mathbf{W}( \cdot, t)-\Pi_k\mathbf{W}(\cdot,t))\|_{L^2}^2\, \dd t \leq  C\frac{k_n^{2(\mu-2)}}{q_n^{2(s-3)}}\|\mathbf{W}\|_{H^s(I_n,  L^2)}^2, 
\end{equation} 
\begin{equation}\label{prop6}
\int_{t_{n-1}}^{t_n} \| \partial_{t}(\mathbf{W}( \cdot, t)-\Pi_k\mathbf{W}(\cdot,t))\|_{L^2}^2 \,\dd t \leq  C\frac{k_n^{2(\mu-1)}}{q_n^{2(s-1)}}\|\mathbf{W}\|_{H^s(I_n,  L^2)}^2, 
\end{equation} 
\begin{equation}\label{prop7}
\int_{t_{n-1}}^{t_n} \| \mathbf{W}( \cdot, t)-\Pi_k\mathbf{W}(\cdot,t)\|_{L^2}^2 \,\dd t \leq  C\frac{k_n^{2\mu}}{q_n^{2(s-1)}}\|\mathbf{W}\|_{H^s(I_n,  L^2)}^2, 
\end{equation} 
where $\mu=\min(q_n+1,s)$ and $q_n$ is the polynomial degree with respect to the variable $t$. If we change the spatial function space from $L^2$ to $H_0^1$ in (\ref{prop5})--(\ref{prop7}), analogous estimates follow. Note that we can also get inverse inequalities with respect to the time derivatives in an analogous manner as (ii,a). That is, there exists a positive constant $C_2$ such that, for each fixed $x\in\Omega$, for every $\bm{\varphi}\in \mathcal{V}_{kh}^{\mathbf{q}}$,
\begin{align}
\|\partial_t \bm{\varphi}(x,\cdot)\|_{L^2(I_n)} &\leq C_2 k_n^{-1}\| \bm{\varphi}(x,\cdot)\|_{L^2(I_n)} \label{inveqtime}\\
\| \partial_t\bm{\varphi}(x,\cdot)\|_{L^{\infty}(I_n)}&\leq C_2 k_n^{-1}\| \bm{\varphi}(x,\cdot)\|_{L^{\infty}(I_n)}\label{inveqtime2},
\end{align}
for each $n=1,2,\ldots,N.$ We now decompose the error as 
\begin{align*}
\mathbf{v}_{\mathrm{DG}}(t)-\mathbf{v}(t) =&\:(\mathbf{v}_{\mathrm{DG}}(t)-\Pi_k\mathbf{W}(t))+(\Pi_k\mathbf{W}(t)- \mathbf{W}(t))+( \mathbf{W}(t)-\mathbf{v}(t))&&\\
:=&\:\bm{\theta}(t)+\bm{\rho}_1(t) +\bm{\rho}_2(t) 
\end{align*}
for $t\in I_n, n=1,2,\ldots,N.$ First note that 
\begin{align}
\int_{t_{n-1}}^{t_n} \|\bm{\rho}_1^{(j)}(t)\|_{L^2}^2 \,\dd t =  \int_{t_{n-1}}^{t_n} \| \Pi_k \mathbf{W}^{(j)}(t)-\mathbf{W} ^{(j)}(t) \|_{L^2}^2 \,\dd t\leq &\:C\frac{k_n^{2(\mu-j)}}{q_n^{2(s-1)}}\int_{t_{n-1}}^{t_n} \sum_{\alpha=0}^{s}\|\mathbf{W}^{(\alpha)}(t)\|_{L^2}^2 \,\dd t\nonumber\\
{}\leq &\: C_1(\mathbf{v}) \frac{k_n^{2(\mu-j)+1}}{q_n^{2(s-1)}}\quad\text{ for } j=0,1,
\end{align}
where we have applied (\ref{prop6}) and (\ref{prop7}). Here $\mu=\min(s,q_n+1)$ and $\Pi_k \mathbf{W}^{ (j)}=\frac{\dd^{j}\Pi_k \mathbf{W}}{\dd t^j}$, $\mathbf{W}^{(j)}=\frac{\dd^j \mathbf{W}}{\dd t^j},$ with $j=0,1.$ If we assume that the solution $\mathbf{v}\in W^{s,\infty}([0,T]; H^{m}\cap H_0^1)$ of (\ref{setup2})--(\ref{ic2}) is sufficiently smooth (i.e. $s>q_n+1$), then we can write
\begin{align}\label{rho1_approx}
\int_{t_{n-1}}^{t_n} \|\bm{\rho}_1^{(j)}(t)\|_{L^2}^2 \,\dd t &\leq \:C  \frac{k_n^{2(q_n+1-j)}}{q_n^{2(s-1)}} \int_{t_{n-1}}^{t_n}\sum_{\alpha=0}^{s} \|\mathbf{W}^{(\alpha)}(t)\|_{L^2}^2 \,\dd t\leq \: C_1(\mathbf{v})\frac{k_n^{2(q_n+1-j)+1}}{q_n^{2(s-1)}} 
\end{align}
for $j=0,1.$ By the property (iii,b) of the elliptic projection, we know that
\begin{align}\label{rho2_approx}
\int_{t_{n-1}}^{t_n} \|\bm{\rho}_2^{(j)}(t)\|_{L^2}^2 \,\dd t= & \int_{t_{n-1}}^{t_n}\|\mathbf{W}^{(j)}(t)-\mathbf{v}^{(j)}(t)\|_{L^2}^2 \,\dd t \leq  \:C_2(\mathbf{v}) k_n h^{2r+2}
\end{align}
for $j=0,1,2.$ Here $C_i(\mathbf{v})$ for $i=1,2$ are constants depending on the exact solution $\mathbf{v}$. Recall that the fully discrete scheme is
\begin{flalign}\label{fully_discrete}
&\sum_{n=1}^N \int_{t_{n-1}}^{t_n} \left( \mathbf{\ddot v}_{\mathrm{DG}}(t), \bm{\dot \varphi}(t) \right)_{L^2}\, \dd t+\sum_{n=1}^{N-1} \left([\mathbf{\dot v}_{\mathrm{DG}}(t) ]_n, \bm{\dot \varphi}(t_n^{+}) \right)_{L^2}+\left(\mathbf{\dot v}_{\mathrm{DG}}(t_0^{+}), \bm{\dot \varphi}(t_0^{+} )\right)_{L^2}&&\nonumber\\
&\quad+\gamma^2 \sum_{n=1}^N \int_{t_{n-1}}^{t_n} \left( \mathbf{v}_{\mathrm{DG}}(t), \bm{\dot \varphi}(t)\right)_{L^2} \,\dd t +\gamma^2\sum_{n=1}^{N-1} \left([\mathbf{v}_{\mathrm{DG}}(t)]_n,  \bm{\varphi}(t_n^{+}) \right)_{L^2}+\gamma^2\left(\mathbf{v}_{\mathrm{DG}}(t_0^{+}),  \bm{\varphi}(t_0^{+}) \right)_{L^2} &&\nonumber\\
&\quad+\sum_{n=1}^N \int_{t_{n-1}}^{t_n}\sum_{i,\alpha=1}^d \mathrm{e}^{-\gamma t}\left( S_{i\alpha}(\nabla \mathbf{v}_{\mathrm{DG}}(t) \mathrm{e}^{\gamma t}),\partial_{\alpha}\bm{\dot \varphi}_i(t) \right)_{L^2}\,\dd t +2\gamma\sum_{n=1}^N \int_{t_{n-1}}^{t_n} \left( \mathbf{\dot v}_{\mathrm{DG}}(t), \bm{\dot \varphi}(t) \right)_{L^2} \,\dd t &&\nonumber\\
&\quad+\sum_{n=0}^{N-1} \frac{1}{2}\bigg\{\int_0^1 \tilde{a}(\tau \mathbf{v}_{\mathrm{DG}}(t_n^{+})\mathrm{e}^{\gamma t_n}; \mathbf{v}_{\mathrm{DG}}(t_n^{+}), \bm{\varphi}(t_n^{+}))\,\dd \tau -\int_0^1 \tilde{a}(\tau \mathbf{v}_{\mathrm{DG}}(t_n^{-})\mathrm{e}^{\gamma t_n}; \mathbf{v}_{\mathrm{DG}}(t_n^{-}), \bm{\varphi}(t_n^{-}))\,\dd \tau\bigg\}&&\nonumber\\
&=\left( \mathbf{v}_{1,h}, \bm{\dot \varphi}(t_0^{+})\right)_{L^2}+\gamma^2\left(\mathbf{v}_{0,h},\bm{\varphi}(t_0^{+})\right)_{L^2}+\sum_{n=1}^N\int_{t_{n-1}}^{t_n} \left(\mathbf{\tilde{f}}(t), \bm{\dot \varphi}(t) \right)_{L^2}\, \dd t,\quad \mbox{ for } \bm{\varphi}\in\mathcal{V}_{kh}^{\mathbf{q}}. &&
\end{flalign}
The variational form of the original problem is written as
\begin{flalign}\label{original_vform}
&\sum_{n=1}^N \int_{t_{n-1}}^{t_n} \left( \mathbf{\ddot v}(t), \bm{\dot \varphi}(t) \right)_{L^2} \,\dd t+\sum_{n=1}^{N-1} \left([\mathbf{\dot v}(t) ]_n, \bm{\dot \varphi}(t_n^{+}) \right)_{L^2}+\left(\mathbf{\dot v}(t_0^{+}), \bm{\dot \varphi}(t_0^{+}) \right)_{L^2}&&\nonumber\\
&\quad+\gamma^2 \sum_{n=1}^N \int_{t_{n-1}}^{t_n} \left( \mathbf{ v}(t), \bm{\dot \varphi}(t)\right)_{L^2}\, \dd t +\gamma^2\sum_{n=1}^{N-1} \left([\mathbf{v}(t) ]_n,  \bm{\varphi}(t_n^{+}) \right)_{L^2}+\gamma^2\left(\mathbf{v}(t_0^{+}),  \bm{\varphi}(t_0^{+}) \right)_{L^2}&&\nonumber\\
&\quad+\sum_{n=1}^N \int_{t_{n-1}}^{t_n}\sum_{i,\alpha=1}^d \mathrm{e}^{-\gamma t}\left( S_{i\alpha}(\nabla \mathbf{v}(t)\mathrm{e}^{\gamma t}),\partial_{\alpha}\bm{\dot \varphi}_i (t)\right)_{L^2}\, \dd t+2\gamma\sum_{n=1}^N \int_{t_{n-1}}^{t_n} \left( \mathbf{\dot v}(t), \bm{\dot \varphi}(t) \right)_{L^2}\, \dd t &&\nonumber\\
&\quad+\sum_{n=0}^{N-1}\frac{1}{2}\bigg\{ \int_0^1 \tilde{a}(\tau \mathbf{v}(t_n^{+})\mathrm{e}^{\gamma t_n}; \mathbf{v}(t_n^{+}), \bm{\varphi}(t_n^{+}))\,\dd \tau- \int_0^1 \tilde{a}(\tau \mathbf{v}(t_n^{-})\mathrm{e}^{\gamma t_n}; \mathbf{v}(t_n^{-}), \bm{\varphi}(t_n^{-}))\,\dd \tau\bigg\} &&\nonumber\\
&=\left(\mathbf{v}_1, \bm{\dot \varphi}(t_0^{+})\right)_{L^2}+\gamma^2\left(\mathbf{v}_0,\bm{\varphi}(t_0^{+})\right)_{L^2}+\sum_{n=1}^N \int_{t_{n-1}}^{t_n} \left(\tilde{\mathbf{f}}(t), \bm{\dot \varphi}(t)\right)_{L^2} \,\dd t, \quad\mbox{ for } \bm{\varphi}\in\mathcal{V}_{kh}^{\mathbf{q}}. &&
\end{flalign}
By considering the nonlinear elliptic projection of $\mathbf{v}(t)$  (cf. equality (\ref{nonlinear projection})), we can replace 
\begin{align*}
\sum_{i,\alpha=1}^d\mathrm{e}^{-\gamma t}\left(S_{i\alpha}(\nabla\mathbf{v} (t)\mathrm{e}^{\gamma t}),\partial_{\alpha}\bm{\dot \varphi}_i(t)\right)_{L^2}\,\dd t\quad 
\mbox{ by }\quad
\sum_{i,\alpha=1}^d\mathrm{e}^{-\gamma t}\left(S_{i\alpha}(\nabla\mathbf{W}(t) \mathrm{e}^{\gamma t}),\partial_{\alpha}\bm{\dot \varphi}_i(t)\right)_{L^2}\,\dd t
\end{align*}
in (\ref{original_vform}). Using the continuity of $\mathbf{v}(t)$ and $\mathbf{W}(t)$ in time, we can also replace
\begin{align*}
\sum_{n=0}^{N-1} \bigg\{ \int_0^1 \tilde{a}(\tau \mathbf{v}(t_n^{+})\mathrm{e}^{\gamma t_n}; \mathbf{v}(t_n^{+}), \bm{\varphi}(t_n^{+}))\,\dd \tau- \int_0^1 \tilde{a}(\tau \mathbf{v}(t_n^{-})\mathrm{e}^{\gamma t_n}; \mathbf{v}(t_n^{-}), \bm{\varphi}(t_n^{-}))\,\dd \tau\bigg\}
\end{align*}
by 
\begin{align*}
\sum_{n=0}^{N-1} \bigg\{ \int_0^1 \tilde{a}(\tau \mathbf{W}(t_n^{+})\mathrm{e}^{\gamma t_n}; \mathbf{W}(t_n^{+}), \bm{\varphi}(t_n^{+}))\,\dd \tau- \int_0^1 \tilde{a}(\tau \mathbf{W}(t_n^{-})\mathrm{e}^{\gamma t_n}; \mathbf{W}(t_n^{-}), \bm{\varphi}(t_n^{-}))\,\dd \tau\bigg\}
\end{align*}
in (\ref{original_vform}). Subtracting the resulting equality from (\ref{fully_discrete}), we have
\begin{flalign}\label{error_vform}
&\sum_{n=1}^N\int_{t_{n-1}}^{t_n} \left(  \bm{\ddot\theta}(t)+ \bm{\ddot \rho_1}(t)+\bm{\ddot \rho_2}(t), \bm{\dot \varphi}(t)\right)_{L^2} \,\dd t +\sum_{n=1}^{N-1}\left( [\bm{\dot \theta}(t)+\bm{\dot \rho_1}(t)+\bm{\dot \rho_2}(t)]_{n}, \bm{\dot \varphi}(t_n^{+})\right)_{L^2}&&\nonumber\\
&\quad+\left( \mathbf{\dot v}_{\mathrm{DG}}(t_0^{+})-\mathbf{\dot v}(t_0^{+}), \bm{\dot \varphi}(t_0^{+})\right)_{L^2} +2\gamma \sum_{n=1}^N \int_{t_{n-1}}^{t_n} \left( \bm{\dot \theta}(t)+\bm{\dot\rho_1}(t)+\bm{\dot\rho_2}(t), \bm{\dot \varphi}(t)\right)_{L^2} \,\dd t &&\nonumber\\
&\quad+\gamma^2\sum_{n=1}^N \int_{t_{n-1}}^{t_n} \left( \bm{\theta}(t)+\bm{\rho_1}(t)+\bm{\rho_2}(t),\bm{\dot \varphi}(t)\right)_{L^2} \,\dd t +\gamma^2 \sum_{n=1}^{N-1}\left([\bm{\theta}(t)+\bm{\rho_1}(t)+\bm{\rho_2}(t)]_n, \bm{\varphi}(t_n^{+})\right)_{L^2} &&\nonumber\\
&\quad +\gamma^2\left(\mathbf{v}_{\mathrm{DG}}(t_0^{+})-\mathbf{v}(t_0^{+}), \bm{\varphi}(t_0^{+})\right)_{L^2}+\sum_{n=1}^N\int_{t_{n-1}}^{t_n}\sum_{i,\alpha=1}^d \mathrm{e}^{-\gamma t}\left(S_{i \alpha}(\nabla \mathbf{v}_{\mathrm{DG}}(t)\mathrm{e}^{\gamma t})-S_{i \alpha}(\nabla \mathbf{W}(t)\mathrm{e}^{\gamma t}),\partial_{\alpha}\bm{\dot \varphi}_i\right)_{L^2} \,\dd t&&\nonumber\\
&\quad+\sum_{n=0}^{N-1}\frac{1}{2} \bigg\{ \int_0^1 \tilde{a}(\tau \mathbf{v}_{\mathrm{DG}}(t_n^{+})\mathrm{e}^{\gamma t_n}; \mathbf{v}_{\mathrm{DG}}(t_n^{+}), \bm{\varphi}(t_n^{+}))\:\dd \tau -\int_0^1 \tilde{a}(\tau \mathbf{W}(t_n^{+})\mathrm{e}^{\gamma t_n}; \mathbf{W}(t_n^{+}), \bm{\varphi}(t_n^{+}))\:\:\dd \tau\bigg\}&&\nonumber\\
&\quad+\sum_{n=0}^{N-1} \frac{1}{2}\bigg\{ \int_0^1 \tilde{a}(\tau \mathbf{W}(t_n^{-})\mathrm{e}^{\gamma t_n}; \mathbf{W}(t_n^{-}), \bm{\varphi}(t_n^{-}))\:\dd \tau - \int_0^1 \tilde{a}(\tau \mathbf{v}_{\mathrm{DG}}(t_n^{-})\mathrm{e}^{\gamma t_n}; \mathbf{v}_{\mathrm{DG}}(t_n^{-}), \bm{\varphi}(t_n^{-}))\:\dd \tau \bigg\}&&\nonumber\\
&=\:\left(\mathbf{v}_{1,h}-\mathbf{v}_1,\bm{\dot \varphi}(t_0^{+})\right)_{L^2}+\gamma^2\left(\mathbf{v}_{0,h}-\mathbf{v},\bm{\varphi}(t_0^{+})\right)_{L^2}, \quad \text { for } \bm{\varphi}\in\mathcal{V}_{kh}^{\mathbf{q}}.&&
\end{flalign}
Now we consider the integral on $I_n=(t_{n-1}, t_n]$ only, 
\begin{flalign}\label{local_err}
&\int_{t_{n-1}}^{t_n} \left(  \bm{\ddot\theta}(t), \bm{\dot \varphi}(t)\right)_{L^2} \dd t +\left( [\bm{\dot \theta}(t)]_{n-1}, \bm{\dot \varphi}(t_{n-1}^{+})\right)_{L^2}+2\gamma \int_{t_{n-1}}^{t_n} \left( \bm{\dot \theta}(t), \bm{\dot \varphi}(t)\right)_{L^2}+\gamma^2\int_{t_{n-1}}^{t_n} \left( \bm{\theta}(t),\bm{\dot \varphi}(t)\right)_{L^2} \dd t && \nonumber\\
&\quad+\gamma^2 \left([\bm{\theta}(t)]_{n-1}, \bm{\varphi}(t_{n-1}^{+})\right)_{L^2}+\int_{t_{n-1}}^{t_n}\sum_{i,\alpha=1}^d \mathrm{e}^{-\gamma t}\left(S_{i \alpha}(\nabla \mathbf{v}_{\mathrm{DG}}(t)\mathrm{e}^{\gamma t})-S_{i \alpha}(\nabla \mathbf{W}(t)\mathrm{e}^{\gamma t}),\partial_{\alpha}\bm{\dot \varphi}_i(t)\right)_{L^2}\dd t &&\nonumber\\
&\quad+ \frac{1}{2}\bigg\{ \int_0^1 \tilde{a}(\tau \mathbf{v}_{\mathrm{DG}}(t_{n-1}^{+})\mathrm{e}^{\gamma t_{n-1}}; \mathbf{v}_{\mathrm{DG}}(t_{n-1}^{+}), \bm{\varphi}(t_{n-1}^{+}))\dd \tau -\int_0^1 \tilde{a}(\tau \mathbf{W}(t_{n-1}^{+})\mathrm{e}^{\gamma t_{n-1}}; \mathbf{W}(t_{n-1}^{+}), \bm{\varphi}(t_{n-1}^{+}))\dd \tau\bigg\}&&\nonumber\\
&\quad+\frac{1}{2}\bigg\{ \int_0^1 \tilde{a}(\tau \mathbf{W}(t_{n-1}^{-})\mathrm{e}^{\gamma t_{n-1}}; \mathbf{W}(t_{n-1}^{-}), \bm{\varphi}(t_{n-1}^{-}))\dd \tau - \int_0^1 \tilde{a}(\tau \mathbf{v}_{\mathrm{DG}}(t_{n-1}^{-})\mathrm{e}^{\gamma t_{n-1}}; \mathbf{v}_{\mathrm{DG}}(t_{n-1}^{-}), \bm{\varphi}(t_{n-1}^{-}))\dd \tau \bigg\}&&\nonumber\\
&=-\int_{t_{n-1}}^{t_n} \left( \bm{\ddot \rho_1}(t),\bm{\dot \varphi}(t)\right)_{L^2} \dd t-\gamma^2\int_{t_{n-1}}^{t_n}\left(\bm{\rho_1}(t),\bm{\dot \varphi}(t)\right)_{L^2} \dd t-2\gamma\int_{t_{n-1}}^{t_n}\left( \bm{\dot \rho_1}(t),\bm{\dot \varphi}(t) \right)_{L^2} \dd t &&\nonumber\\
&\quad- \left([ \bm{\dot \rho_1}]_{n-1},\bm{\dot \varphi}(t_{n-1}^{+})\right)_{L^2}-\gamma^2\left( [\bm{\rho_1}]_{n-1},\bm{\varphi}(t_{n-1}^{+})\right)_{L^2}  &&\nonumber\\
&\quad -\int_{t_{n-1}}^{t_n} \left( \bm{\ddot \rho_2}(t),\bm{\dot \varphi}(t)\right)_{L^2} \dd t-\gamma^2\int_{t_{n-1}}^{t_n}\left(\bm{\rho_2}(t),\bm{\dot \varphi}(t)\right)_{L^2} \dd t-2\gamma\int_{t_{n-1}}^{t_n}\left( \bm{\dot \rho_2}(t),\bm{\dot \varphi}(t) \right)_{L^2} \dd t&&\nonumber\\
&=\int_{t_{n-1}}^{t_n} \left( \bm{\dot \rho_1}(t),\bm{\ddot \varphi}(t)\right)_{L^2} \dd t-\gamma^2\int_{t_{n-1}}^{t_n}\left(\bm{\rho_1}(t),\bm{\dot \varphi}(t)\right)_{L^2} \dd t-2\gamma\int_{t_{n-1}}^{t_n}\left( \bm{\dot \rho_1}(t),\bm{\dot \varphi}(t) \right)_{L^2} \dd t &&\nonumber\\
&\quad- \left(\bm{\dot \rho_1}(t_{n}^{-}),\bm{\dot \varphi}(t_{n}^{-})\right)_{L^2}+\left(\bm{\dot \rho_1}(t_{n-1}^{-}),\bm{\dot \varphi}(t_{n-1}^{+})\right)_{L^2}-\gamma^2\left( [\bm{\rho_1}]_{n-1},\bm{\varphi}(t_{n-1}^{+})\right)_{L^2} &&\nonumber\\
&\quad -\int_{t_{n-1}}^{t_n} \left( \bm{\ddot \rho_2}(t),\bm{\dot \varphi}(t)\right)_{L^2} \dd t-\gamma^2\int_{t_{n-1}}^{t_n}\left(\bm{\rho_2}(t),\bm{\dot \varphi}(t)\right)_{L^2} \dd t-2\gamma\int_{t_{n-1}}^{t_n}\left( \bm{\dot \rho_2}(t),\bm{\dot \varphi}(t) \right)_{L^2} \dd t&&\nonumber\\
&=-\gamma^2\int_{t_{n-1}}^{t_n}\left(\bm{\rho_1}(t),\bm{\dot\varphi}(t)\right)_{L^2} \dd t-2\gamma\int_{t_{n-1}}^{t_n}\left( \bm{\dot \rho_1}(t),\bm{\dot \varphi}(t) \right)_{L^2} \dd t  -\int_{t_{n-1}}^{t_n} \left( \bm{\ddot \rho_2}(t),\bm{\dot \varphi}(t)\right)_{L^2} \dd t&&\nonumber\\
&\quad-\gamma^2\int_{t_{n-1}}^{t_n}\left(\bm{\rho_2}(t),\bm{\dot \varphi}(t)\right)_{L^2} \dd t-2\gamma\int_{t_{n-1}}^{t_n}\left( \bm{\dot \rho_2}(t),\bm{\dot \varphi}(t) \right)_{L^2} \dd t,&&
\end{flalign}
where we have used the fact that $\bm{\rho_2}(t)$ and $\bm{\dot \rho_2}(t)$ are continuous in time and properties (\ref{prop1})--(\ref{prop2}). By Taylor's theorem with an integral remainder, we have
\begin{flalign*}
&S_{i \alpha}(\nabla\mathbf{v}_{\mathrm{DG}}(t)\mathrm{e}^{\gamma t})= S_{i \alpha}(\nabla\mathbf{W}(t)\mathrm{e}^{\gamma t})\\
&\quad+\sum_{j,\beta=1}^d \mathrm{e}^{\gamma t} \partial_{\beta}(\mathbf{v}_{\mathrm{DG}}(t)-\mathbf{W}(t))_j \int_0^1 \frac{\partial}{\partial \eta_{j\beta}}S_{i \alpha}(\nabla \mathbf{W}(t)\mathrm{e}^{\gamma t}+\tau(\nabla \mathbf{v}_{\mathrm{DG}}(t)-\nabla \mathbf{W}(t))\mathrm{e}^{\gamma t})\dd \tau. 
\end{flalign*}
If $\nabla\mathbf{v}_{\mathrm{DG}}(t)\mathrm{e}^{\gamma t}\in\mathcal{Z}_{\delta}, \nabla\mathbf{W}(t)\mathrm{e}^{\gamma t}\in \mathcal{Z}_{\delta}$ for each $t\in [0,T]$, we have $\nabla \mathbf{W}(t)\mathrm{e}^{\gamma t}+\tau(\nabla \mathbf{v}_{\mathrm{DG}}(t)-\nabla \mathbf{W}(t))\mathrm{e}^{\gamma t}\in\mathcal{Z}_{\delta}$ for $0\leq \tau \leq 1$ by the convexity of $\mathcal{Z}_{\delta}.$ This implies that the term in the integral remainder is well-defined. Thus, we can write
\begin{flalign*}
&\sum_{i,\alpha=1}^d \mathrm{e}^{-\gamma t}\left(S_{i \alpha}(\nabla \mathbf{v}_{\mathrm{DG}}(t)\mathrm{e}^{\gamma t})-S_{i \alpha}(\nabla \mathbf{W}(t)\mathrm{e}^{\gamma t}), \partial_{\alpha}\bm{\dot \varphi}_i(t)\right)_{L^2} \dd t &&\\
&= \int_0^1 \sum_{i,\alpha,j,\beta=1}^d \left( A_{i\alpha j\beta}(\nabla \mathbf{W}(t)\mathrm{e}^{\gamma t}+\tau(\nabla\bm{\phi}(t)-\nabla\mathbf{W}(t)) \mathrm{e}^{\gamma t})\partial_{\beta}(\mathbf{v}_{\mathrm{DG}}(t)-\mathbf{W}(t))_j, \partial_{\alpha}\bm{\dot \varphi}_i(t) \right)_{L^2} \dd \tau &&\\
&=\int_0^1 \tilde{a}(\mathbf{W}(t)\mathrm{e}^{\gamma t}+\tau(\mathbf{v}_{\mathrm{DG}}(t)-\mathbf{W}(t))\mathrm{e}^{\gamma t}; \mathbf{v}_{\mathrm{DG}}(t)-\mathbf{W}(t),\bm{\dot \varphi}(t))\dd \tau.
\end{flalign*}
For simplicity of notation, we write $$\mathsf{A}(\mathbf{v}_{\mathrm{DG}}(t)\mathrm{e}^{\gamma t}; \mathbf{v}_{\mathrm{DG}}(t)-\mathbf{W}(t), \bm{\dot \varphi}(t)):=\int_0^1 \tilde{a}(\mathbf{W}(t)\mathrm{e}^{\gamma t}+\tau(\mathbf{v}_{\mathrm{DG}}(t)-\mathbf{W}(t))\mathrm{e}^{\gamma t}; \mathbf{v}_{\mathrm{DG}}(t)-\mathbf{W}(t),\bm{\dot \varphi}(t))\dd \tau.$$
Analogously,
\begin{flalign*}
& \int_0^1 \tilde{a}(\tau \mathbf{v}_{\mathrm{DG}}(t_{n-1}^{\pm})\mathrm{e}^{\gamma t_{n-1}}; \mathbf{v}_{\mathrm{DG}}(t_{n-1}^{\pm}), \bm{\varphi}(t_{n-1}^{\pm}))\,\dd \tau -\int_0^1 \tilde{a}(\tau \mathbf{W}(t_{n-1}^{\pm})\mathrm{e}^{\gamma t_{n-1}}; \mathbf{W}(t_{n-1}^{\pm}), \bm{\varphi}(t_{n-1}^{\pm})\,\dd \tau &&\\
&= \sum_{i,\alpha=1}^d \mathrm{e}^{-\gamma t_{n-1}}\bigg\{\left( S_{i \alpha}(\nabla \mathbf{v}_{\mathrm{DG}}(t_{n-1}^{\pm})\mathrm{e}^{\gamma t_{n-1}}), \partial_{\alpha}\bm{ \varphi}_i(t_{n-1}^{\pm})\right)_{L^2} -\left( S_{i \alpha}(\nabla \mathbf{W}(t_{n-1}^{\pm})\mathrm{e}^{\gamma t_{n-1}}), \partial_{\alpha}\bm{ \varphi}_i(t_{n-1}^{\pm})\right)_{L^2} \bigg\}&&\\
&= \sum_{i,\alpha=1}^d \mathrm{e}^{-\gamma t_{n-1}}\left( S_{i \alpha}(\nabla \mathbf{v}_{\mathrm{DG}}(t_{n-1}^{\pm})\mathrm{e}^{\gamma t_{n-1}})-S_{i \alpha}(\nabla \mathbf{W}(t_{n-1}^{\pm})\mathrm{e}^{\gamma t_{n-1}}), \partial_{\alpha}\bm{ \varphi}_i(t_{n-1}^{\pm})\right)_{L^2} &&\\
&= \:\mathsf{A}(\mathbf{v}_{\mathrm{DG}}(t_{n-1}^{\pm})\mathrm{e}^{\gamma t_{n-1}}; \mathbf{v}_{\mathrm{DG}}(t_{n-1}^{\pm})-\mathbf{W}(t_{n-1}^{\pm}), \bm{\varphi}(t_{n-1}^{\pm})).&&
\end{flalign*}
Since $\mathbf{W}(t_n^{\pm})-\Pi_k \mathbf{W}(t_n^{\pm})=0$ for each $n=1,2,\ldots, N$, our equation (\ref{local_err}) becomes
\begin{flalign}\label{simple_eq}    
&\int_{t_{n-1}}^{t_n} \left( \bm{\ddot \theta}(t), \bm{\dot \varphi}(t)\right)_{L^2}\, \dd t+\left( [\bm{\dot \theta}(t)]_{n-1}, \bm{\dot \varphi} (t_{n-1}^{+})\right)_{L^2}+2\gamma  \int_{t_{n-1}}^{t_n} \left( \bm{\dot \theta}(t), \bm{\dot \varphi}(t)\right)_{L^2} \,\dd t  &&\nonumber\\
&\quad+\gamma^2\int_{t_{n-1}}^{t_n} \left( \bm{\theta}(t),\bm{\dot \varphi}(t)\right)_{L^2} \dd t+\gamma^2\left([\bm{\theta}(t)]_{n-1}, \bm{\varphi}(t_{n-1}^{+})\right)_{L^2}&&\nonumber\\
&\quad+\int_{t_{n-1}}^{t_n} \mathsf{A}(\mathbf{v}_{\mathrm{DG}}(t)\mathrm{e}^{\gamma t}; \mathbf{v}_{\mathrm{DG}}(t)-\Pi_k\mathbf{W}(t), \bm{\dot \varphi}(t))\,\dd t +\frac{1}{2}\mathsf{A}(\mathbf{v}_{\mathrm{DG}}(t_{n-1}^{+})\mathrm{e}^{\gamma t_{n-1}};\bm{\theta}(t_{n-1}^{+}), \bm{\varphi}(t_{n-1}^{+})) &&\nonumber\\
&=-\gamma^2\int_{t_{n-1}}^{t_n}\left(\bm{\rho_1}(t),\bm{\dot \varphi}(t)\right)_{L^2} \,\dd t-2\gamma\int_{t_{n-1}}^{t_n}\left( \bm{\dot \rho_1}(t),\bm{\dot \varphi}(t) \right)_{L^2} \dd t -\int_{t_{n-1}}^{t_n} \left( \bm{\ddot \rho_2}(t),\bm{\dot \varphi}(t)\right)_{L^2} \,\dd t &&\nonumber\\
&\quad-\gamma^2\int_{t_{n-1}}^{t_n}\left(\bm{\rho_2}(t),\bm{\dot \varphi}(t)\right)_{L^2} \dd t-2\gamma\int_{t_{n-1}}^{t_n}\left( \bm{\dot \rho_2}(t),\bm{\dot \varphi}(t) \right)_{L^2} \,\dd t &&\nonumber\\
&\quad+\int_{t_{n-1}}^{t_n}\mathsf{A}(\mathbf{v}_{\mathrm{DG}}(t)\mathrm{e}^{\gamma t}; \mathbf{W}(t)-\Pi_k\mathbf{W}(t),\bm{\dot \varphi}(t))\, \dd t + \frac{1}{2}\mathsf{A}(\mathbf{v}_{\mathrm{DG}}(t_{n-1}^{-})\mathrm{e}^{\gamma t_{n-1}};\bm{\theta}(t_{n-1}^{-}), \bm{\varphi}(t_{n-1}^{-}))&&\nonumber
\end{flalign}
for $\bm{\varphi}\in\mathcal{V}_{kh}^{\mathbf{q}}$. Consider the following subset of $\mathcal{V}_{kh}^{\mathbf{q}}$ defined by
\begin{flalign*}
\mathcal{F}:=&\bigg\{\bm{\psi}\in\mathcal{V}_{kh}^{\mathbf{q}} \mid \mbox{ for each } j=1,2,\ldots,N, \,\|\bm{\psi}(t_j^{-})-\Pi_k\mathbf{W}(t_j^{-})\|_{H^1}^2 +\|\partial_t(\bm{ \psi}(t_j^{-})-\Pi_k\mathbf{W}(t_j^{-}))\|_{L^2}^2 &&\\
&+\sum_{n=1}^j\int_{t_{n-1}}^{t_n}\|\partial_t (\bm{\psi}(t)-\Pi_k\mathbf{W}(t))\|_{L^2}^2\dd t\leq  C_{\ast}(\mathbf{v})\left(\sum_{n=1}^j k_n h^{2r+2}+ \frac{k_n^{2q_n+1}}{q_n^{2(s-1)}}\right), \mathrm{e}^{\gamma t}\nabla \bm{\psi}(t)\in \mathcal{Z}_{\delta}\bigg\},&&
\end{flalign*}
where $C_{\ast}(\mathbf{v})$ is a positive constant depending on the solution $\mathbf{v}$, which will be specified later. First note that $\mathcal{F}$ is non-empty since $\Pi_k\mathbf{W}\in\mathcal{F}$. In addition, $\mathcal{F}$ is a closed and convex subset of $\mathcal{V}_{kh}^{\mathbf{q}}$ in the topology induced by the norm $\|\cdot\|_{\mathcal{F}}$, which is defined by 
$$\|\bm{\varphi}\|_{\mathcal{F}}=\max_{t\in I_n, 1\leq n\leq N}(\|\bm{\varphi}(t)\|_{H^1}+\|\dot{\bm{\varphi}}(t)\|_{L^2})$$ for $\bm{\varphi}\in\mathcal{V}_{kh}^{\mathbf{q}}.$  With this notation, we are ready to define a fixed point mapping $\mathcal{N}$ on $\mathcal{F}$ as follows: if $\bm{\phi}\in\mathcal{F}$, the image $\mathbf{v}_{\phi}:=\mathcal{N}(\bm{\phi})$ is given by the relation
\begin{equation}
\mathbf{v}_{\phi}(0)=\mathbf{v}_{0,h}, \:\: \mathbf{\dot v}_{\phi}(0)=\mathbf{v}_{1,h},
\end{equation}
\begin{flalign}\label{fixed_point}
&\int_{t_{n-1}}^{t_n} \left( \bm{\ddot \theta}_{\phi}(t), \bm{\dot \varphi}(t)\right)_{L^2} \dd t+\left( [\bm{\dot \theta}_{\phi}(t)]_{n-1}, \bm{\dot \varphi}(t_{n-1}^{+})\right)_{L^2}+2\gamma  \int_{t_{n-1}}^{t_n} \left( \bm{\dot \theta}_{\phi}(t), \bm{\dot \varphi}(t)\right)_{L^2} \,\dd t&&\nonumber\\
& \quad+\gamma^2\int_{t_{n-1}}^{t_n} \left( \bm{\theta}_{\phi}(t),\bm{\dot \varphi}(t)\right)_{L^2} \dd t +\gamma^2\left([\bm{\theta}_{\phi}(t)]_{n-1}, \bm{\varphi}(t_{n-1}^{+})\right)_{L^2}\nonumber&&\\
&\quad+\int_{t_{n-1}}^{t_n} \mathsf{A}(\bm{\phi}(t)\mathrm{e}^{\gamma t}; \bm{\theta}_{\phi}(t), \bm{\dot \varphi}(t))\,\dd t+\frac{1}{2}\mathsf{A}(\bm{\phi}(t_{n-1}^{+})\mathrm{e}^{\gamma t_{n-1}};\bm{\theta}_{\phi}(t_{n-1}^{+}), \bm{\varphi}(t_{n-1}^{+}))&&\nonumber\\
&=-\gamma^2\int_{t_{n-1}}^{t_n}\left(\bm{\rho_1}(t),\bm{\dot \varphi}(t)\right)_{L^2} \dd t-2\gamma\int_{t_{n-1}}^{t_n}\left( \bm{\dot\rho_1}(t),\bm{\dot \varphi}(t) \right)_{L^2}\, \dd t  -\int_{t_{n-1}}^{t_n} \left( \bm{\ddot \rho_2}(t),\bm{\dot \varphi}(t)\right)_{L^2} \,\dd t&&\\
&\quad-\gamma^2\int_{t_{n-1}}^{t_n}\left(\bm{\rho_2}(t),\bm{\dot \varphi}(t)\right)_{L^2} \dd t-2\gamma\int_{t_{n-1}}^{t_n}\left( \bm{\dot \rho_2}(t),\bm{\dot \varphi}(t) \right)_{L^2} \,\dd t&&\nonumber\\
&\quad+\int_{t_{n-1}}^{t_n}\mathsf{A}(\bm{\phi}(t)\mathrm{e}^{\gamma t}; \mathbf{W}(t)-\Pi_k\mathbf{W}(t),\bm{\dot \varphi}(t) ) \,\dd t+\frac{1}{2}\mathsf{A}(\bm{\phi}(t_{n-1}^{-})\mathrm{e}^{\gamma t_{n-1}};\bm{\theta}_{\phi}(t_{n-1}^{-}), \bm{\varphi}(t_{n-1}^{-})),&&\nonumber
\end{flalign}
where $\bm{\theta}_{\phi}=\mathbf{v}_{\phi}-\Pi_k\mathbf{W}.$

In order to complete the proof of the theorem, it suffices to show that, for each $n=1,\ldots, N$, the map $\mathcal{N}$ defined by (\ref{fixed_point}) has a unique fixed point in $\mathcal{F}$. If $\mathbf{v}_{\mathrm{DG}}\in\mathcal{F}$ is this fixed point, then $\mathbf{v}_{\mathrm{DG}}$ is a solution to (\ref{discrete_vform}). 
\subsection{Auxiliary results}
If we take $\bm{\varphi}=\bm{\theta_{\phi}}$ in (\ref{fixed_point}), then the nonlinear term inside the integral becomes
$$\mathsf{A}(\bm{\phi}(t)\mathrm{e}^{\gamma t}; \bm{\theta}_{\phi}(t), \bm{\dot \varphi}(t))=\mathsf{A}(\bm{\phi}(t)\mathrm{e}^{\gamma t}; \bm{\theta}_{\phi}(t), \bm{\dot \theta}_{\phi}(t)).$$
Following the proof in \cite{Ma}, it is crucial to replace the expression $\mathsf{A}(\bm{\phi}(t)\mathrm{e}^{\gamma t}; \bm{\theta}_{\phi}(t), \bm{\dot \theta}_{\phi}(t))$ by $$\frac{1}{2} \frac{\dd }{\dd t}\mathsf{A}(\bm{\phi}(t)\mathrm{e}^{\gamma t}; \bm{\theta}_{\phi}(t), \bm{\theta}_{\phi}(t))-\frac{1}{2}\sum_{i,\alpha,j,\beta=1}^d \int_0^1 \left( \partial_t  A_{i\alpha j \beta}^{\tau}\partial_{\beta}\bm{\theta}_{\phi,j}(t),\partial_{\alpha} \bm{\theta}_{\phi, i}(t)\right)_{L^2}\,\dd \tau,$$
where $A_{i\alpha j \beta}^{\tau}:=A_{i\alpha j \beta}(\nabla \mathbf{W}(t)\mathrm{e}^{\gamma t}+\tau(\nabla\bm{\phi}(t)-\nabla\mathbf{W}(t)) \mathrm{e}^{\gamma t})$ and $t\in I_n, n=1,2\ldots N.$ We shall need an estimate on the expression
$$\mathsf{A}_t(\bm{\phi}(t)\mathrm{e}^{\gamma t}; \bm{\varphi}(t),\bm{\psi}(t)):=\frac{1}{2}\sum_{i,\alpha,j,\beta=1}^d\int_0^1 \left( \partial_t  A_{i\alpha j \beta}^{\tau}\partial_{\beta}\bm{\varphi}_{j}(t),\partial_{\alpha} \bm{\psi}_{i}(t)\right)_{L^2}\,\dd \tau$$ 
for $\bm{\varphi}, \bm{\psi}\in\mathcal{V}_{hk}^{\mathbf{q}}$, $t\in I_n, n=1,2,\ldots N.$
\begin{lmm}\label{bdlemma}
Under the assumptions stated in Theorem \ref{convergence_thm}, there exists a constant $C_{\tau}>0$ such that, for $t\in I_n$, $n=1,2\ldots N$,
\begin{equation}\label{lemmaie}
\vert \mathsf{A}_{t}(\bm{\phi}(t)\mathrm{e}^{\gamma t}; \bm{\varphi}(t),\bm{\psi}(t))\vert\leq C_{\tau}\|\nabla \bm{\varphi}(t)\|_{L^2}\|\nabla\bm{\psi}(t)\|_{L^2}.
\end{equation}
\end{lmm}
\begin{proof}
See Section \ref{lemm_pf}.
\end{proof}

\subsection{Convergence proof}\label{convergenceproof}
We will establish the existence of a unique fixed point in $\mathcal{F}$ by showing that the pair $\mathcal{F}$ and $\mathcal{N}$ satisfies the assumptions of Banach's fixed point theorem, namely that
\begin{enumerate}
\item[(a)] $\mathcal{N}(\mathcal{F})\subset\mathcal{F}.$
\item[(b)] 
$\mathcal{N}$ is a contraction with respect to $d(\cdot,\cdot)$ where for $\bm{\phi},\bm{\varphi}\in\mathcal{F}$, 
$$d(\bm{\phi},\bm{\varphi}):=\max_{t\in I_n, 1\leq n\leq N}\left(\|\bm{\phi}(t)-\bm{\varphi}(t)\|_{H^1}+\|\bm{\dot{\phi}}(t)-\bm{\dot{\varphi}}(t)\|_{L^2}\right).$$
\end{enumerate}
\subsubsection{Existence of a fixed point of $\mathcal{N}$ in $\mathcal{F}$}
For (a), we first observe that $\mathcal{N}$ is well-defined. Indeed, if $\bm{\phi}\in\mathcal{F}$, since $\nabla \mathbf{W}\mathrm{e}^{\gamma t}\in\mathcal{Z}_{\delta}$, $\nabla\mathbf{W}\mathrm{e}^{\gamma t}+\tau(\nabla \bm{\phi}-\nabla\mathbf{W})\mathrm{e}^{\gamma t}\in\mathcal{Z}_{\delta}$ for $0\leq \tau\leq 1$, and the bilinear form $\mathsf{A}(\bm{\phi}(t)\mathbf{e}^{\gamma t}; \cdot, \cdot)$ is symmetric and positive definite.
Taking $\bm{\varphi}=\bm{\theta}_{\phi}$ in (\ref{fixed_point}) and replacing $\mathsf{A}(\bm{\phi}(t)\mathrm{e}^{\gamma t}; \bm{\theta}_{\phi}(t), \bm{\dot \theta}_{\phi}(t))$ by 
$$\frac{1}{2} \frac{\dd }{\dd t}\mathsf{A}(\bm{\phi}(t)\mathrm{e}^{\gamma t}; \bm{\theta}_{\phi}(t), \bm{\theta}_{\phi}(t))-\frac{1}{2}\sum_{i,\alpha,j,\beta=1}^d \int_0^1 \left( \partial_t  A_{i\alpha j \beta}^{\tau}\partial_{\beta}\bm{\theta}_{\phi,j}(t),\partial_{\alpha} \bm{\theta}_{\phi, i}(t)\right)_{L^2}\dd \tau,$$
we obtain
\begin{flalign}\label{local_energy}
&\|\bm{\dot \theta}_{\phi} (t_n^{-})\|_{L^2}^2+\|\bm{\dot \theta}_{\phi}(t_{n-1}^{+})\|_{L^2}^2+\gamma^2\|\bm{\theta}_{\phi} (t_n^{-})\|_{L^2}^2+\gamma^2\|\bm{\theta}_{\phi}(t_{n-1}^{+})\|_{L^2}^2+4\gamma\int_{t_{n-1}}^{t_n} \|\bm{\dot \theta}_{\phi}(t)\|_{L^2}^2\, \dd t &&\nonumber\\
&\quad+\mathsf{A}(\bm{\phi}(t_n^{-})\mathrm{e}^{\gamma t_n}; \bm{\theta}_{\phi}(t_{n}^{-}), \bm{\theta}_{\phi}(t_{n}^{-}))-\mathsf{A}(\bm{\phi}(t_{n-1}^{-}) e^{\gamma t_{n-1}}; \bm{\theta}_{\phi}(t_{n-1}^{-}), \bm{\theta}_{\phi}(t_{n-1}^{-}))&&\nonumber\\
&= -2\gamma^2\int_{t_{n-1}}^{t_n}\left(\bm{ \rho_1}(t),\bm{\dot\theta_{\phi}}(t)\right)_{L^2} \dd t-4\gamma\int_{t_{n-1}}^{t_n}\left( \bm{\dot \rho_1}(t),\bm{\dot \theta_{\phi}}(t) \right)_{L^2}\, \dd t &&\\
&\quad-2\int_{t_{n-1}}^{t_n} \left( \bm{\ddot \rho_2}(t),\bm{\dot \theta_{\phi}}(t)\right)_{L^2} \,\dd t-2\gamma^2\int_{t_{n-1}}^{t_n}\left(\bm{\rho_2}(t),\bm{\dot \theta_{\phi}}(t)\right)_{L^2} \dd t -4\gamma\int_{t_{n-1}}^{t_n}\left( \bm{\dot \rho_2}(t),\bm{\dot \theta_{\phi}}(t) \right)_{L^2} \,\dd t &&\nonumber\\
&\quad+2\left(\bm{\dot \theta_{\phi}}(t_{n-1}^{-}), \bm{\dot \theta_{\phi}}(t_{n-1}^{+})\right)_{L^2}+2\gamma^2\left(\bm{\theta_{\phi}}(t_{n-1}^{-}), \bm{\theta_{\phi}}(t_{n-1}^{+})\right)_{L^2}&&\nonumber\\
&\quad+2\int_{t_{n-1}}^{t_n}\mathsf{A}(\bm{\phi}(t)\mathrm{e}^{\gamma t}; \mathbf{W}(t)-\Pi_k\mathbf{W}(t),\bm{\dot \theta_{\phi}}(t)) \,\dd t+\sum_{i,\alpha,j,\beta=1}^d\int_{t_{n-1}}^{t_n}\int_0^1 \left(\partial_t A_{i\alpha j \beta}^{\tau}\partial_{\beta}\bm{\theta}_{\phi,j}(t), \partial_{\alpha}\bm{\theta}_{\phi,i}(t)\right)_{L^2} \,\dd \tau\, \dd t.&&\nonumber
\end{flalign}
Now we need to bound the terms on the right-hand side of equation (\ref{local_energy}).
By using (\ref{rho1_approx}), (\ref{rho2_approx}) and Young's inequality, we have
\begin{flalign}\label{app1}
&\bigg\vert \:-2\gamma^2\int_{t_{n-1}}^{t_n}\left(\bm{\rho_1}(t),\bm{\dot\theta_{\phi}}(t)\right)_{L^2} \dd t-4\gamma\int_{t_{n-1}}^{t_n}\left( \bm{\dot \rho_1}(t),\bm{\dot \theta_{\phi}}(t) \right)_{L^2}\, \dd t &&\nonumber\\
&\quad -2\int_{t_{n-1}}^{t_n} \left( \bm{\ddot \rho_2}(t),\bm{\dot \theta_{\phi}}(t)\right)_{L^2}\, \dd t-2\gamma^2\int_{t_{n-1}}^{t_n}\left(\bm{\rho_2}(t),\bm{\dot \theta_{\phi}}(t)\right)_{L^2} \dd t-4\gamma\int_{t_{n-1}}^{t_n}\left( \bm{\dot \rho_2}(t),\bm{\dot \theta_{\phi}}(t) \right)_{L^2} \,\dd t\bigg\vert \nonumber\\
&\leq  3\gamma \int_{t_{n-1}}^{t_n} \|\bm{\dot \theta_{\phi}}(t)\|_{L^2}^2\, \dd t+C_1(\gamma)\int_{t_{n-1}}^{t_n} \bigg(\|\bm{\dot \rho_1}(t)\|_{L^2}^2+\|\bm{\rho_1}(t)\|_{L^2}^2\bigg)\,\dd t \nonumber\\
&\quad+C_2(\gamma)\int_{t_{n-1}}^{t_n} \bigg(\|\bm{\ddot \rho_2}(t)\|_{L^2}^2+\|\bm{\dot \rho_2}(t)\|_{L^2}^2+\|\bm{\rho_2}(t)\|_{L^2}^2\bigg)\,\dd t\nonumber\\
&\leq 3\gamma \int_{t_{n-1}}^{t_n} \|\bm{\dot \theta_{\phi}}(t)\|_{L^2}^2\, \dd t +c_1(\gamma,\mathbf{v})\frac{k_n^{2q_n+1}}{q_n^{2(s-1)}}+ c_2(\gamma,\mathbf{v}) k_n h^{2r+2} \nonumber\\
&\leq   3\gamma \int_{t_{n-1}}^{t_n} \|\bm{\dot \theta_{\phi}}(t)\|_{L^2}^2\, \dd t +C(\gamma,\mathbf{v})\left(k_n h^{2r+2}+\frac{k_n^{2q_n+1}}{q_n^{2(s-1)}}\right),
\end{flalign}
where $C_i(\gamma)$ for $i=1,2$ are constants depending on $\gamma$ only, while $C(\gamma,\mathbf{v})$ and $c_i(\gamma,\mathbf{v})$ for $i=1,2$ are constants depending on both $\gamma$ and the exact solution $\mathbf{v}$. By Cauchy--Schwarz inequality, we obtain
\begin{align}\label{cs1}
2\left( \bm{\dot \theta}_{\phi}(t_{n-1}^{-}),\bm{\dot \theta}_{\phi}(t_{n-1}^{+})\right)_{L^2}\leq \|\bm{\dot \theta}_{\phi}(t_{n-1}^{+})\|_{L^2}^2+\|\bm{\dot \theta}_{\phi}(t_{n-1}^{-})\|_{L^2}^2,
\end{align}
and 
\begin{align}\label{cs2}
2\gamma^2\left( \bm{\theta}_{\phi}(t_{n-1}^{-}),\bm{\theta}_{\phi}(t_{n-1}^{+})\right)_{L^2}\leq \gamma^2\|\bm{\theta}_{\phi}(t_{n-1}^{+})\|_{L^2}^2+\gamma^2\|\bm{\theta}_{\phi}(t_{n-1}^{-})\|_{L^2}^2.
\end{align}
Note that 
\begin{flalign*}
&\mathsf{A}(\phi(t)\mathrm{e}^{\gamma t}; \mathbf{W}(t)-\Pi_k\mathbf{W}(t), \bm{\dot \theta}_{\phi}(t))=\frac{\dd }{\dd t}\mathsf{A}(\phi(t)\mathrm{e}^{\gamma t}; \mathbf{W}(t)-\Pi_k\mathbf{W}(t), \bm{\theta}_{\phi}(t)) &&\\
&-\mathsf{A}(\phi(t)\mathrm{e}^{\gamma t};\partial_t (\mathbf{W}(t)-\Pi_k\mathbf{W}(t)), \bm{\theta}_{\phi}(t))-\sum_{i,\alpha,j,\beta=1}^d \int_0^1 \bigg(\partial_t A_{i\alpha j \beta}^{\tau}\partial_{\beta}(\mathbf{W}-\Pi_k\mathbf{W})_{j}, \partial_{\alpha} \bm{\theta}_{\phi,i}(t)\bigg)_{L^2} \,\dd \tau.&&
\end{flalign*}
Using the fact that $\left(\mathbf{W}-\Pi_k\mathbf{W}\right)(t_{n}^{-})=\left(\mathbf{W}-\Pi_k\mathbf{W}\right)(t_{n-1}^{+})=0$ for $n=1,2,\ldots,N,$ we have
\begin{flalign}\label{non-symmetric}
&\int_{t_{n-1}}^{t_n}\mathsf{A}(\bm{\phi}(t)\mathrm{e}^{\gamma t}; \mathbf{W}(t)-\Pi_k\mathbf{W}(t), \bm{\dot \theta_{\phi}}(t)) \,\dd t= \:-\int_{t_{n-1}}^{t_n}\mathsf{A}(\bm{\phi}(t)\mathrm{e}^{\gamma t}; \partial_t\left(\mathbf{W}(t)-\Pi_k\mathbf{W}(t)\right), \bm{\theta_{\phi}}(t))\,\dd t &&\nonumber\\
&-\int_{t_{n-1}}^{t_n} \sum_{i,\alpha,j,\beta=1}^d \int_0^1 \bigg(\partial_t A_{i\alpha j \beta}^{\tau}\partial_{\beta}(\mathbf{W}(t)-\Pi_k\mathbf{W}(t))_{j}, \partial_{\alpha} \bm{\theta}_{\phi,i}(t) \bigg)_{L^2} \,\dd \tau\, \dd t.
\end{flalign}
Then
\begin{flalign}
&\bigg \vert \int_{t_{n-1}}^{t_n}\mathsf{A}(\bm{\phi}(t)\mathrm{e}^{\gamma t}; \mathbf{W}(t)-\Pi_k\mathbf{W}(t), \bm{\dot \theta_{\phi}}(t)) \,\dd t\bigg\vert \leq \bigg \vert -\int_{t_{n-1}}^{t_n}\mathsf{A}(\bm{\phi}(t)\mathrm{e}^{\gamma t}; \partial_t\left(\mathbf{W}(t)-\Pi_k\mathbf{W}(t)\right), \bm{\theta_{\phi}}(t)) \, \dd t \bigg \vert &&\nonumber\\  
&\quad+\bigg \vert\int_{t_{n-1}}^{t_n} \sum_{i,\alpha,j,\beta=1}^d \int_0^1 \bigg(\partial_t A_{i\alpha j \beta}\partial_{\beta}(\mathbf{W}(t)-\Pi_k\mathbf{W}(t))_{j}, \partial_{\alpha} \bm{\theta}_{\phi, i}(t) \bigg)_{L^2} \,\dd \tau \,\dd t \bigg\vert&&\nonumber\\
&\leq  \int_{t_{n-1}}^{t_n} K_{\delta}\|\nabla \partial_t (\mathbf{W}(t)-\Pi_k\mathbf{W}(t))\|_{L^2}\|\nabla \bm{\theta_{\phi}}(t)\|_{L^2} \,\dd t+ C_{\tau}  \int_{t_{n-1}}^{t_n} \|\nabla(\mathbf{W}(t)-\Pi_k\mathbf{W}(t))\|_{L^2}\|\nabla \bm{\theta_{\phi}}(t)\|_{L^2} \,\dd t&&\nonumber\\
&\leq \left( \frac{K_{\delta}}{2}+\frac{C_{\tau}}{2}\right)\int_{t_{n-1}}^{t_n} \|\nabla \bm{\theta}_{\phi}(t)\|_{L^2}^2 \,\dd t+\frac{K_{\delta}}{2}\int_{t_{n-1}}^{t_n}\|\nabla(\partial_t (\mathbf{W}(t)-\Pi_k \mathbf{W}(t)))\|_{L^2}^2 \,\dd t &&\nonumber\\
&\quad+ \frac{C_{\tau}}{2} \int_{t_{n-1}}^{t_n}\|\nabla(\mathbf{W}(t)-\Pi_k\mathbf{W}(t))\|_{L^2}^2\, \dd t &&\nonumber\\
&\leq \left( \frac{K_{\delta}}{2}+\frac{C_{\tau}}{2}\right)\int_{t_{n-1}}^{t_n} \|\nabla \bm{\theta_{\phi}}(t)\|_{L^2}^2 \,\dd t+C(\mathbf{v})\frac{k_n^{2q_n+1}}{q_n^{2(s-1)}} \quad\text{ for } \mathbf{W}\in W^{s, \infty}([0,T]; H_0^1).
\end{flalign}
To bound the terms involving $\partial_tA_{i\alpha j \beta}^{\tau}$, we apply Lemma \ref{bdlemma} to get
\begin{flalign}\label{At1}
\bigg\vert \int_{t_{n-1}}^{t_n}\sum_{i,\alpha,j,\beta=1}^d\int_0^1 \left(\partial_t A_{i\alpha j \beta}^{\tau}\partial_{\beta}\bm{\theta}_{\phi,j}(t), \partial_{\alpha}\bm{\theta}_{\phi,i}(t)\right)_{L^2} \dd \tau \dd t\bigg\vert \leq \:C_{\tau} \int_{t_{n-1}}^{t_n} \|\nabla\bm{\theta}_{\phi}(t)\|_{L^2}^2 \,\dd t.&&
\end{flalign}
Combining (\ref{local_energy})--(\ref{At1}), we obtain  
\begin{flalign}\label{induction}
&\|\bm{\dot \theta}_{\phi}(t_n^{-})\|_{L^2}^2+\gamma^2 \|\bm{ \theta}_{\phi}(t_n^{-})\|_{L^2}^2+\gamma\int_{t_{n-1}}^{t_n} \|\bm{\dot \theta}_{\phi}(t)\|_{L^2}^2\dd t +\mathsf{A}(\bm{\phi}(t_n^{-})\mathrm{e}^{\gamma t_n}; \bm{\theta}_{\phi}(t_{n}^{-}), \bm{\theta}_{\phi}(t_{n}^{-}))&&\nonumber\\
&\leq  \mathsf{A}(\bm{\phi}(t_{n-1}^{-}); \bm{\theta}_{\phi}(t_{n-1}^{-}), \bm{\theta}_{\phi}(t_{n-1}^{-}))+ \tilde{C}\int_{t_{n-1}}^{t_n} \|\nabla\bm{\theta}_{\phi}(t)\|_{L^2}^2 \dd t+ \|\bm{\dot\theta}_{\phi}(t_{n-1}^{-})\|_{L^2}^2+\gamma^2 \|\bm{ \theta}_{\phi}(t_{n-1}^{-})\|_{L^2}^2&&\nonumber\\
&\quad+C(\gamma,\mathbf{v})\left(k_n h^{2r+2}+\frac{k_n^{2q_n+1}}{q_n^{2(s-1)}}\right)+ C(\mathbf{v})\frac{k_n^{2q_n+1}}{q_n^{2(s-1)}}.
\end{flalign}
Summing up over $n=1,\ldots, j$, we obtain 
\begin{flalign}\label{sum}
&\|\bm{\dot \theta}_{\phi}(t_j^{-})\|_{L^2}^2+\gamma^2 \|\bm{ \theta}_{\phi}(t_j^{-})\|_{L^2}^2+\gamma \sum_{n=1}^j \int_{t_{n-1}}^{t_n}\|\bm{\dot \theta}_{\phi}(t)\|_{L^2}^2\dd t+\mathsf{A}(\bm{\phi}(t_j^{-})\mathrm{e}^{\gamma t_j}; \bm{\theta}_{\phi}(t_j^{-}), \bm{\theta}_{\phi}(t_j^{-})) &&\nonumber\\
&\leq \tilde{C}\sum_{n=1}^j\int_{t_{n-1}}^{t_n} \|\nabla\bm{\theta}_{\phi}(t)\|_{L^2}^2 \dd t+C(\mathbf{v})\sum_{n=1}^j\left(k_n h^{2r+2}+\frac{k_n^{2q_n+1}}{q_n^{2(s-1)}}\right).
\end{flalign}
Using the the coercivity of $\mathsf{A}(\bm{\phi}(t_j^{-})\mathrm{e}^{\gamma t_n}; \bm{\theta}_{\phi}(t_j^{-}), \bm{\theta}_{\phi}(t_j^{-}))$ (i.e. assumption (S2b)), we have 
\begin{flalign}\label{sumforgronwall}
&\|\bm{\dot \theta}_{\phi}(t_j^{-})\|_{L^2}^2+\gamma^2 \|\bm{ \theta}_{\phi}(t_j^{-})\|_{L^2}^2+\gamma \sum_{n=1}^j \int_{t_{n-1}}^{t_n}\|\bm{\dot \theta}_{\phi}(t)\|_{L^2}^2\dd t+ M_1\|\nabla\bm{\theta}_{\phi}(t_j^{-})\|_{L^2}^2 &&\nonumber\\
&\leq \tilde{C}\sum_{n=1}^j\int_{t_{n-1}}^{t_n} \|\nabla\bm{\theta}_{\phi}(t)\|_{L^2}^2 \dd t+C(\mathbf{v})\sum_{n=1}^j\left(k_n h^{2r+2}+\frac{k_n^{2q_n+1}}{q_n^{2(s-1)}}\right).
\end{flalign}
By the fundamental theorem of calculus and the triangle inequality, we have for each $t\in I_n$, with $n=1,\ldots N$,
\begin{align}\label{FCT}
\|\nabla \bm{\theta}_{\phi}(t)\|_{L^2}^2&\leq\left(\|\nabla \bm{\theta}_{\phi}(t_n^{-})\|_{L^2}+\int_{t_{n-1}}^{t_n}\|\nabla\bm{\dot \theta}_{\phi}(t)\|_{L^2} \dd t\right)^2\leq  2\|\nabla\bm{\theta}_{\phi}(t_n^{-})\|_{L^2}^2+2C_0^2h^{-2}k_n\int_{t_{n-1}}^{t_n}\|\bm{\dot \theta}_{\phi}(t)\|_{L^2}^2\dd t.
\end{align}
Substituting (\ref{FTC}) into (\ref{sumforgronwall}), we have
\begin{flalign}\label{sumforgronwall2}
&\|\bm{\dot \theta}_{\phi}(t_j^{-})\|_{L^2}^2+\gamma^2 \|\bm{ \theta}_{\phi}(t_j^{-})\|_{L^2}^2+\gamma \sum_{n=1}^j \int_{t_{n-1}}^{t_n}\|\bm{\dot \theta}_{\phi}(t)\|_{L^2}^2\dd t+ M_1\|\nabla\bm{\theta}_{\phi}(t_j^{-})\|_{L^2}^2 &&\nonumber\\
&\leq  2\tilde{C}\sum_{n=1}^{j-1} k_n\|\nabla\bm{\theta}_{\phi}(t_n^{-})\|_{L^2}^2  +2\tilde{C}k_j\|\nabla\bm{\theta}_{\phi}(t_j^{-})\|_{L^2}^2 + 2\tilde{C} C_0^2\sum_{n=1}^{j}  h^{-2}k_n^2\int_{t_{n-1}}^{t_n}\|\bm{\dot \theta}_{\phi}(t)\|_{L^2}^2\dd t&&\nonumber\\
&\quad+C(\mathbf{v})\sum_{n=1}^j\left(k_n h^{2r+2}+\frac{k_n^{2q_n+1}}{q_n^{2(s-1)}}\right) &&\nonumber\\
&\leq  2\tilde{C}\sum_{n=1}^{j-1} k_n\|\nabla\bm{\theta}_{\phi}(t_n^{-})\|_{L^2}^2  +2\tilde{C}k_j\|\nabla\bm{\theta}_{\phi}(t_j^{-})\|_{L^2}^2 + \hat{C}\sum_{n=1}^{j}k_n\int_{t_{n-1}}^{t_n}\|\bm{\dot \theta}_{\phi}(t)\|_{L^2}^2\dd t&&\nonumber\\
&\quad+C(\mathbf{v})\sum_{n=1}^j\left(k_n h^{2r+2}+\frac{k_n^{2q_n+1}}{q_n^{2(s-1)}}\right),
\end{flalign}
where the last inequality follows from the assumption that $\mu_i k_i\leq  h^2$ for each $i=1,\ldots, N$, with $\hat{C}=2\tilde{C}C_0^2\max_{1\leq i\leq j}\frac{1}{\mu_i}.$ 
The term $2\tilde{C}k_j\|\nabla\bm{\theta}_{\phi}(t_j^{-})\|_{L^2}^2$ and the sum of integrals on the right-hand side of (\ref{sumforgronwall2}) can be absorbed into the third and fourth terms of the left-hand side of (\ref{sumforgronwall2}) if we choose each time step $k_n$ is sufficiently small. That is,
\begin{flalign}\label{sumforgronwall3}
&\|\bm{\dot \theta}_{\phi}(t_j^{-})\|_{L^2}^2+\gamma^2 \|\bm{ \theta}_{\phi}(t_j^{-})\|_{L^2}^2+\sum_{n=1}^j(\gamma-\hat{C}k_n)\int_{t_{n-1}}^{t_n}\|\bm{\dot \theta}_{\phi}(t)\|_{L^2}^2\dd t+ (M_1-2\tilde{C}k_n)\|\nabla\bm{\theta}_{\phi}(t_n^{-})\|_{L^2}^2 &&\nonumber\\
&\leq  2\tilde{C}\sum_{n=1}^{j-1} k_n\|\nabla\bm{\theta}_{\phi}(t_n^{-})\|_{L^2}^2 +C(\mathbf{v})\sum_{n=1}^j\left(k_n h^{2r+2}+\frac{k_n^{2q_n+1}}{q_n^{2(s-1)}}\right).
\end{flalign}
By choosing $k_n\leq \min\{\frac{\gamma}{2\hat{C}}, \frac{M_1}{4\tilde{C}} \}$ for each $n=1,\ldots, N$ and applying the discrete Gr\"{o}nwall lemma, we have 
\begin{flalign}
&\|\bm{\dot \theta}_{\phi}(t_j^{-})\|_{L^2}^2+ \|\bm{ \theta}_{\phi}(t_j^{-})\|_{L^2}^2+\sum_{n=1}^j\int_{t_{n-1}}^{t_n} \|\bm{\dot \theta}_{\phi}(t)\|_{L^2}^2\dd t +\|\nabla\bm{\theta}_{\phi}(t_j^{-}))\|_{L^2}^2\nonumber &&\\
&\leq C_j(\mathbf{v})\exp\left(C \sum_{n=1}^j k_n\right)\sum_{n=1}^j\left(k_n h^{2r+2}+\frac{k_n^{2q_n+1}}{q_n^{2(s-1)}}\right)\leq \: C_{\mathrm{max}}(\mathbf{v})\sum_{n=1}^j\left(k_n h^{2r+2}+\frac{k_n^{2q_n+1}}{q_n^{2(s-1)}}\right),&&
\end{flalign}
where $C_{\mathrm{max}}(\mathbf{v})=\max_{1\leq n\leq N} C_n(\mathbf{v})\exp(C T)$. Now tracing back constants through the previous estimates, we notice that $C_{\mathrm{max}}(\mathbf{v})$ does not depend on $C_{\ast}(\mathbf{v})$, so we can define $C_{\ast}(\mathbf{v}):=C_{\mathrm{max}}(\mathbf{v}).$ Note that
$$\|\nabla\mathbf{v}_{\phi}(t)-\nabla\mathbf{v}(t)\|_{L^{\infty}}\leq \|\nabla\mathbf{v}_{\phi}(t)-\nabla\mathbf{W}(t)\|_{L^{\infty}}+\|\nabla\mathbf{v}(t)-\nabla\mathbf{W}(t)\|_{L^{\infty}}.$$  By the inverse estimate (ii,b), the error bound (iii,a), and the approximation properties of $\mathcal{P}_h$ in the $W^{1,\infty}$ and $H^1$ semi-norms, we can find an $h_1>0$ such that, for $h<h_1$, 
\begin{flalign*}
\|\nabla \mathbf{W}(t)-\nabla \mathbf{v}(t)\|_{L^{\infty}} &\leq  \|\nabla \mathbf{W}(t)-\nabla \mathcal{P}_h\mathbf{v}(t)\|_{L^{\infty}}+\|\nabla\mathcal{P}_h\mathbf{v}(t)-\mathbf{v}(t)\|_{L^{\infty}}&&\\
&\leq C_1 h^{-\frac{d}{2}}\| \nabla(\mathbf{W}(t)-\mathcal{P}_h\mathbf{v}(t))\|_{L^2}+ C(\mathbf{v})h^{r-\frac{d}{2}} \quad\mbox{ (by (ii,b))} && \\
&\leq  C_1 h^{-\frac{d}{2}}\| \nabla(\mathbf{W}(t)-\mathbf{v}(t))\|_{L^2}+C_1 h^{-\frac{d}{2}}\| \nabla(\mathcal{P}_h\mathbf{v}(t)-\mathbf{v}(t))\|_{L^2}+ C(\mathbf{v})h^{r-\frac{d}{2}}&&\\
&\leq \tilde{C}(\mathbf{v}) h^{r-\frac{d}{2}}\leq\frac{\delta}{2}\mathrm{e}^{-\gamma T}. &
\end{flalign*} 
Since $r>\frac{d}{2}+1$, $k_i^{q_i-\frac{1}{2}}=o(h^{1+\frac{d}{2}})$ and $\mu_i k_i\leq h^2$ for each $i=1,\ldots N$, we can also choose $h_2>0$ such that, for  $t\in I_n, n=1,2\ldots,N,$ for $h<h_2$,
\begin{flalign*}
\|\nabla \mathbf{v}_{\phi}(t)-\nabla \mathbf{W}(t)\|_{L^{\infty}}&\leq C_1  h^{-\frac{d}{2}}\bigg( \|\nabla  \bm{\theta}_{\phi}(t)\|_{L^2}+\|\nabla\Pi_k\mathbf{W}(t)-\nabla\mathbf{W}(t)\|_{L^2}\bigg)&&\\
&\leq  C_1 h^{-\frac{d}{2}}\bigg( \|\nabla \bm{\theta}_{\phi}(t_n^{-})\|_{L^2}+ \int_{t_{n-1}}^{t_n} \| \partial_t  \nabla\bm{\theta}_{\phi}(t)\|_{L^2} \,\dd t&&\\
&\quad+\| \nabla \Pi_k\mathbf{W}(t_n^{-})-\nabla\mathbf{W}(t_n^{-})\|_{L^2}+\int_{t_{n-1}}^{t_n} \| \partial_t( \nabla \Pi_k\mathbf{W}-\nabla\mathbf{W})(t)\|_{L^2}\, \dd t \bigg)&&\\
&\leq C_1 h^{-\frac{d}{2}}\bigg( \|\nabla \bm{\theta}_{\phi}(t_n^{-})\|_{L^2}+ \int_{t_{n-1}}^{t_n} \| \partial_t  \nabla\bm{\theta}_{\phi}(t)\|_{L^2} \,\dd t\bigg)&&\\
&\quad+ C_1 h^{-\frac{d}{2}}\sqrt{k_n}\left(\int_{t_{n-1}}^{t_n} \| \partial_t( \nabla \Pi_k\mathbf{W}-\nabla\mathbf{W})(t)\|_{L^2}^2 \,\dd t \right)^{\frac{1}{2}}&&\\
&\leq  C_1 h^{-\frac{d}{2}}\bigg( \|\nabla \bm{\theta}_{\phi}(t_n^{-})\|_{L^2}+ C_0h^{-1}\sqrt{k_n}\int_{t_{n-1}}^{t_n} \| \partial_t\bm{\theta}_{\phi}(t)\|_{L^2}^2 \,\dd t\bigg)^{\frac{1}{2}}+C(\mathbf{W})h^{-\frac{d}{2}} \frac{k_n^{q_n+1}}{q_n^{s-1}}&&\\
&\leq C(\mathbf{v})h^{-\frac{d}{2}}\left(\sum_{i=1}^n k_i h^{2r+2}+\frac{k_i^{2q_i+1}}{q_i^{2(s-1)}} \right)^{\frac{1}{2}}< \frac{\delta}{2}\mathrm{e}^{-\gamma T}.
\end{flalign*}
By choosing $h<h_{\ast}=\min\{h_0, h_1,h_2\}$, we obtain $\nabla \mathbf{v}_{\phi}\mathrm{e}^{\gamma t}\in\mathcal{Z}_{\delta}.$ 
\subsubsection{Verification that $\mathcal{N}$ is a contraction mapping}
To show the contraction property (b), we consider $\mathbf{R}=\bm{\phi}-\bm{\phi'}$ and $\bm{\Theta}=\mathbf{v}_{\phi}-\mathbf{v}_{\phi'}$ where $\bm{\phi}, \bm{\phi'}\in\mathcal{F}.$  Replacing $\bm{\phi}$ in (\ref{fixed_point}) by $\bm{\phi'}$ and subtracting the new equation from (\ref{fixed_point}), we have
\begin{flalign}\label{theta_eqn}
&\int_{t_{n-1}}^{t_n}\left( \bm{\ddot \Theta}(t),\bm{\dot \varphi}(t)\right)_{L^2} \dd t+\left([\bm{\dot \Theta }(t)]_{n-1}, \bm{\dot \varphi}(t_{n-1}^{+})\right)_{L^2}+2\gamma\int_{t_{n-1}}^{t_n} \left(\bm{\dot \Theta}(t),\bm{\dot \varphi}(t)\right)_{L^2} \dd t+\gamma^2 \int_{t_{n-1}}^{t_n} \left( \bm{\Theta}(t), \bm{\dot \varphi}(t) \right)_{L^2} \,\dd t &&\nonumber\\
&\quad+\gamma^2 \left([\bm{\Theta}(t)]_{n-1}, \bm{\varphi}(t_{n-1}^{+})\right)_{L^2}+\int_{t_{n-1}}^{t_n} \bigg(\mathsf{A}(\bm{\phi}(t)\mathrm{e}^{\gamma t}; \bm{\theta}_{\phi}(t), \bm{\dot \varphi}(t))-\mathsf{A}(\bm{\phi'}(t)\mathrm{e}^{\gamma t}; \bm{\theta}_{\phi'}(t), \bm{\dot \varphi}(t))\bigg)\,\dd t &&\nonumber\\
&\quad+\frac{1}{2}\bigg\{\mathsf{A}(\bm{\phi}(t_{n-1}^{+})\mathrm{e}^{\gamma t_{n-1}};\bm{\theta}_{\phi}(t_{n-1}^{+}), \bm{\varphi}(t_{n-1}^{+}))-\mathsf{A}(\bm{\phi'}(t_{n-1}^{+})\mathrm{e}^{\gamma t_{n-1}};\bm{\theta}_{\phi'}(t_{n-1}^{+}), \bm{\varphi}(t_{n-1}^{+}))\bigg\}&&\nonumber\\
&=\int_{t_{n-1}}^{t_n} \bigg(\mathsf{A}(\bm{\phi}(t)\mathrm{e}^{\gamma t}; \mathbf{W}(t)-\Pi_k\mathbf{W}(t),\bm{\dot \varphi}(t))-\mathsf{A}(\bm{\phi'}(t)\mathrm{e}^{\gamma t}; \mathbf{W}(t)-\Pi_k\mathbf{W}(t),\bm{\dot \varphi}(t)) \bigg )\,\dd t&&\\
&\quad+\frac{1}{2}\bigg\{\mathsf{A}(\bm{\phi}(t_{n-1}^{-})\mathrm{e}^{\gamma t_{n-1}};\bm{\theta}_{\phi}(t_{n-1}^{-}), \bm{\varphi}(t_{n-1}^{-}))-\mathsf{A}(\bm{\phi'}(t_{n-1}^{-})\mathrm{e}^{\gamma t_{n-1}};\bm{\theta}_{\phi'}(t_{n-1}^{-}), \bm{\varphi}(t_{n-1}^{-}))\bigg\}.&&\nonumber
\end{flalign}
Taking $\bm{\varphi}=\bm{\Theta}$ in equation (\ref{theta_eqn}),
replacing 
$$\mathsf{A}(\bm{\phi}(t)\mathrm{e}^{\gamma t}; \bm{\theta}_{\phi}(t), \bm{\dot \Theta}(t))- \mathsf{A}(\bm{\phi'}(t)\mathrm{e}^{\gamma t}; \bm{\theta}_{\phi'}(t), \bm{\dot \Theta}(t))$$
by $$\mathsf{A}(\bm{\phi}(t)\mathrm{e}^{\gamma t}; \bm{\Theta}(t), \bm{\dot \Theta}(t))+ \mathsf{A}(\bm{\phi}(t)\mathrm{e}^{\gamma t}; \bm{\theta}_{\phi'}(t), \bm{\dot \Theta}(t))- \mathsf{A}(\bm{\phi'}(t)\mathrm{e}^{\gamma t}; \bm{\theta}_{\phi'}(t), \bm{\dot \Theta}(t)),$$ and writing
\begin{flalign*}
\int_{t_{n-1}}^{t_n} \mathsf{A}(\bm{\phi}(t)\mathrm{e}^{\gamma t}; \bm{\Theta}(t), \bm{\dot \Theta}(t)) \dd t &=   \frac{1}{2}\int_{t_{n-1}}^{t_n} \frac{\dd }{\dd t} \mathsf{A}(\bm{\phi}(t)\mathrm{e}^{\gamma t}; \bm{\Theta}(t), \bm{ \Theta}(t)) \,\dd t &&\\
&\quad-\frac{1}{2 } \int_{t_{n-1}}^{t_n} \int_0^1 \sum_{i,\alpha,j,\beta=1}^d\bigg( \partial_t A_{i\alpha j \beta}^{\tau}\partial_{\beta}\bm{\Theta}_j(t), \partial_{\alpha}\bm{\Theta}_i(t) \bigg)_{L^2}\,\dd \tau \,\dd t, 
\end{flalign*}
we have
\begin{flalign*}
&\|\bm{\dot \Theta}(t_n^{-})\|_{L^2}^2+\|\bm{\dot \Theta}(t_{n-1}^{+})\|_{L^2}^2+\gamma^2\|\bm{\Theta}(t_{n}^{-})\|_{L^2}^2+\gamma^2\|\bm{\Theta}(t_{n-1}^{+})\|_{L^2}^2+4\gamma\int_{t_{n-1}}^{t_n} \|\bm{\dot \Theta}(t)\|_{L^2}^2\, \dd t&&\\
&\quad+\mathsf{A}(\bm{\phi}(t_n^{-})\mathrm{e}^{\gamma t_n}; \bm{\Theta}(t_{n}^{-}), \bm{\Theta}(t_n^{-}))-\mathsf{A}(\bm{\phi}(t_{n-1}^{-})\mathrm{e}^{\gamma t_{n-1}};\bm{\Theta}(t_{n-1}^{-}),\bm{\Theta}(t_{n-1}^{-}))&&\\
&=2\left( \bm{\dot \Theta}(t_{n-1}^{-}),\bm{\dot \Theta} (t_{n-1}^{+})\right)_{L^2}+2\gamma^2 \left( \bm{\Theta}(t_{n-1}^{-}),\bm{\Theta} (t_{n-1}^{+})\right)_{L^2}&&\\
&\quad+\int_{t_{n-1}}^{t_n} \int_0^1 \sum_{i,\alpha,j,\beta=1}^d \bigg( \partial_t A_{i\alpha j \beta}^{\tau}\partial_{\beta}\bm{\Theta}_j(t), \partial_{\alpha}\bm{\Theta}_i(t) \bigg)_{L^2}\dd \tau \,\dd t &&\\
&\quad+ 2\int_{t_{n-1}}^{t_n} \bigg(\mathsf{A}(\bm{\phi}(t)\mathrm{e}^{\gamma t}; \mathbf{W}(t)-\Pi_k\mathbf{W}(t),\bm{\dot \Theta}(t))-\mathsf{A}(\bm{\phi'}(t)\mathrm{e}^{\gamma t}; \mathbf{W}(t)-\Pi_k\mathbf{W}(t),\bm{\dot \Theta}(t)) \bigg )\,\dd t&&\\
&\quad+2\int_{t_{n-1}}^{t_n}  \bigg(\mathsf{A}(\bm{\phi'}(t)\mathrm{e}^{\gamma t}; \bm{\theta}_{\phi'}(t), \bm{\dot \Theta}(t))- \mathsf{A}(\bm{\phi}(t)\mathrm{e}^{\gamma t}; \bm{\theta}_{\phi'}(t), \bm{\dot \Theta}(t))\bigg) \,\dd t &&\\
&\quad+\bigg(\mathsf{A}(\bm{\phi'}(t_{n-1}^{+})\mathrm{e}^{\gamma t_{n-1}};\bm{\theta}_{\phi'}(t_{n-1}^{+}), \bm{\Theta}(t_{n-1}^{+}))-\mathsf{A}(\bm{\phi}(t_{n-1}^{+})\mathrm{e}^{\gamma t_{n-1}};\bm{\theta}_{\phi'}(t_{n-1}^{+}),\bm{\Theta}(t_{n-1}^{+}))\bigg)&&\\
&\quad+\bigg(\mathsf{A}(\bm{\phi}(t_{n-1}^{-})\mathrm{e}^{\gamma t_{n-1}};\bm{\theta}_{\phi'}(t_{n-1}^{-}), \bm{\Theta}(t_{n-1}^{-}))-\mathsf{A}(\bm{\phi'}(t_{n-1}^{-})\mathrm{e}^{\gamma t_{n-1}};\bm{\theta}_{\phi'}(t_{n-1}^{-}),\bm{\Theta}(t_{n-1}^{-}))\bigg).&&
\end{flalign*}
Note that
\begin{flalign*}
&\int_{t_{n-1}}^{t_n}\mathsf{A}(\bm{\phi'}(t)\mathrm{e}^{\gamma t}; \bm{\theta}_{\phi'}(t),\bm{\dot{\Theta}}(t)) \, \dd t+\mathsf{A}(\bm{\phi'}(t_{n-1}^{+})\mathrm{e}^{\gamma t_{n-1}};\bm{\theta}_{\phi'}(t_{n-1}^{+}), \bm{\Theta}(t_{n-1}^{+}))&&\\
&\quad-\mathsf{A}(\bm{\phi'}(t_{n-1}^{-})\mathrm{e}^{\gamma t_{n-1}};\bm{\theta}_{\phi'}(t_{n-1}^{-}),\bm{\Theta}(t_{n-1}^{-}))&&\\
&= -\int_{t_{n-1}}^{t_n}\mathsf{A}(\bm{\phi'}(t)\mathrm{e}^{\gamma t}; \bm{\dot\theta_{\phi'}}(t),\bm{\Theta}(t)) \,\dd t- \int_{t_{n-1}}^{t_n} \int_0^1 \sum_{i,\alpha,j,\beta=1}^d \bigg( \partial_t \tilde{A}_{i\alpha j\beta}^{\tau}\partial_{\beta}\bm{\theta}_{\phi',j}(t), \partial_{\alpha}\bm{\Theta}_i(t) \bigg)_{L^2}\,\dd \tau\, \dd t &&\\
&\quad+\mathsf{A}(\bm{\phi'}(t_{n}^{-})\mathrm{e}^{\gamma t_{n}};\bm{\theta}_{\phi'}(t_{n}^{-}),\bm{\Theta}(t_{n}^{-}))-\mathsf{A}(\bm{\phi'}(t_{n-1}^{-})\mathrm{e}^{\gamma t_{n-1}};\bm{\theta}_{\phi'}(t_{n-1}^{-}),\bm{\Theta}(t_{n-1}^{-})),
\end{flalign*}
where $\tilde{A}_{i\alpha j\beta}^{\tau}:=A_{i\alpha j \beta}(\nabla \mathbf{W}(t)\mathrm{e}^{\gamma t}+\tau(\nabla\bm{\phi'}(t)-\nabla\mathbf{W}(t)) \mathrm{e}^{\gamma t})$. Analogously, we have
\begin{flalign*}
&-\int_{t_{n-1}}^{t_n}\mathsf{A}(\bm{\phi}(t)\mathrm{e}^{\gamma t}; \bm{\theta_{\phi'}}(t),\bm{\dot{\Theta}}(t)) \,\dd t-\mathsf{A}(\bm{\phi}(t_{n-1}^{+})\mathrm{e}^{\gamma t_{n-1}};\bm{\theta}_{\phi'}(t_{n-1}^{+}), \bm{\Theta}(t_{n-1}^{+})) &&\nonumber\\
&\quad+\mathsf{A}(\bm{\phi}(t_{n-1}^{-})\mathrm{e}^{\gamma t_{n-1}};\bm{\theta}_{\phi'}(t_{n-1}^{-}),\bm{\Theta}(t_{n-1}^{-}))&&\\
&= \int_{t_{n-1}}^{t_n}\mathsf{A}(\bm{\phi}(t)\mathrm{e}^{\gamma t}; \bm{\dot\theta_{\phi'}}(t),\bm{\Theta}(t))\, \dd t+ \int_{t_{n-1}}^{t_n} \int_0^1 \sum_{i,\alpha,j,\beta=1}^d \bigg( \partial_t A_{i\alpha j \beta}^{\tau}\partial_{\beta}\bm{\theta}_{\phi',j}(t), \partial_{\alpha}\bm{\Theta}_i(t) \bigg)_{L^2}\,\dd \tau \,\dd t &&\\
&\quad-\mathsf{A}(\bm{\phi}(t_{n}^{-})\mathrm{e}^{\gamma t_{n}};\bm{\theta}_{\phi'}(t_{n}^{-}),\bm{\Theta}(t_{n}^{-}))+\mathsf{A}(\bm{\phi}(t_{n-1}^{-})\mathrm{e}^{\gamma t_{n-1}};\bm{\theta}_{\phi'}(t_{n-1}^{-}),\bm{\Theta}(t_{n-1}^{-})).
\end{flalign*}
This implies that 
\begin{flalign}\label{beforesum}
&\|\bm{\dot \Theta}(t_n^{-})\|_{L^2}^2+\|\bm{\dot \Theta}(t_{n-1}^{+})\|_{L^2}^2+\gamma^2\|\bm{\Theta}(t_{n}^{-})\|_{L^2}^2+\gamma^2\|\bm{\Theta}(t_{n-1}^{+})\|_{L^2}^2+4\gamma\int_{t_{n-1}}^{t_n} \|\bm{\dot \Theta}(t)\|_{L^2}^2 \,\dd t &&\nonumber\\
&\quad+\mathsf{A}(\bm{\phi}(t_n^{-})\mathrm{e}^{\gamma t_n}; \bm{\Theta}(t_{n}^{-}), \bm{\Theta}(t_n^{-}))-\mathsf{A}(\bm{\phi}(t_{n-1}^{-})\mathrm{e}^{\gamma t_{n-1}};\bm{\Theta}(t_{n-1}^{-}),\bm{\Theta}(t_{n-1}^{-}))&&\nonumber\\
&=2\left( \bm{\dot \Theta}(t_{n-1}^{-}),\bm{\dot \Theta} (t_{n-1}^{+})\right)_{L^2}+2\gamma^2 \left( \bm{\Theta}(t_{n-1}^{-}),\bm{\Theta} (t_{n-1}^{+})\right)_{L^2} &&\\
&\quad+\int_{t_{n-1}}^{t_n} \int_0^1 \sum_{i,\alpha,j,\beta=1}^d \bigg( \partial_t A_{i\alpha j \beta}^{\tau}\partial_{\beta}\bm{\Theta}_j(t), \partial_{\alpha}\bm{\Theta}_i(t) \bigg)_{L^2}\,\dd \tau \,\dd t &&\nonumber\\
&\quad+\int_{t_{n-1}}^{t_n}  \bigg(\mathsf{A}(\bm{\phi'}(t)\mathrm{e}^{\gamma t}; \bm{\theta}_{\phi'}(t), \bm{\dot \Theta}(t))- \mathsf{A}(\bm{\phi}(t)\mathrm{e}^{\gamma t}; \bm{\theta}_{\phi'}(t), \bm{\dot \Theta}(t))\bigg) \,\dd t &&\nonumber\\
&\quad+\int_{t_{n-1}}^{t_n}  \bigg(\mathsf{A}(\bm{\phi}(t)\mathrm{e}^{\gamma t}; \bm{\dot \theta}_{\phi'}(t), \bm{\Theta}(t))- \mathsf{A}(\bm{\phi'}(t)\mathrm{e}^{\gamma t}; \bm{\dot \theta}_{\phi'}(t), \bm{\Theta}(t))\bigg) \,\dd t &&\nonumber\\
&\quad+ 2\int_{t_{n-1}}^{t_n} \bigg(\mathsf{A}(\bm{\phi}(t)\mathrm{e}^{\gamma t}; \mathbf{W}(t)-\Pi_k\mathbf{W}(t),\bm{\dot \Theta}(t))-\mathsf{A}(\bm{\phi'}(t)\mathrm{e}^{\gamma t}; \mathbf{W}(t)-\Pi_k\mathbf{W}(t),\bm{\dot \Theta}(t)) \bigg )\,\dd t&&\nonumber\\
&\quad+\int_{t_{n-1}}^{t_n} \int_0^1 \sum_{i,\alpha,j,\beta=1}^d \bigg( (\partial_t A_{i\alpha j\beta}^{\tau} -\partial_t \tilde{A}_{i\alpha j\beta}^{\tau})\partial_{\beta}\bm{\theta}_{\phi',j}(t), \partial_{\alpha}\bm{\Theta}_i(t) \bigg)_{L^2}\,\dd \tau \,\dd t &&\nonumber\\
&\quad+\mathsf{A}(\bm{\phi'}(t_{n}^{-})\mathrm{e}^{\gamma t_{n}};\bm{\theta}_{\phi'}(t_{n}^{-}),\bm{\Theta}(t_{n}^{-}))-\mathsf{A}(\bm{\phi'}(t_{n-1}^{-})\mathrm{e}^{\gamma t_{n-1}};\bm{\theta}_{\phi'}(t_{n-1}^{-}),\bm{\Theta}(t_{n-1}^{-}))&&\nonumber\\
&\quad-\mathsf{A}(\bm{\phi}(t_{n}^{-})\mathrm{e}^{\gamma t_{n}};\bm{\theta}_{\phi'}(t_{n}^{-}),\bm{\Theta}(t_{n}^{-}))+\mathsf{A}(\bm{\phi}(t_{n-1}^{-})\mathrm{e}^{\gamma t_{n-1}};\bm{\theta}_{\phi'}(t_{n-1}^{-}),\bm{\Theta}(t_{n-1}^{-})).&&\nonumber
\end{flalign}
Again, we need to bound the terms on the right-hand side of the equation (\ref{beforesum}).
By the Cauchy--Schwarz inequality, we have
\begin{equation}\label{cauchy1}
2\left( \bm{\dot \Theta}(t_{n-1}^{-}),\bm{\dot \Theta} (t_{n-1}^{+})\right)_{L^2}\leq  \|  \bm{\dot \Theta}(t_{n-1}^{-})\|_{L^2}^2 +  \|  \bm{\dot \Theta}(t_{n-1}^{+})\|_{L^2}^2, 
\end{equation}
\begin{equation}\label{cauchy2}
2\gamma^2 \left( \bm{\Theta}(t_{n-1}^{-}),\bm{\Theta} (t_{n-1}^{+})\right)_{L^2}\leq \gamma^2 \|  \bm{\Theta}(t_{n-1}^{-})\|_{L^2}^2+ \gamma^2 \|  \bm{\Theta}(t_{n-1}^{+})\|_{L^2}^2.
\end{equation}
By Lemma \ref{bdlemma}, for $\bm{\phi}\in \mathcal{F}$, we have
\begin{equation}
\bigg\vert \int_{t_{n-1}}^{t_n}\int_0^1 \sum_{i,\alpha, j,\beta=1}^d\bigg( \partial_t A_{i\alpha j \beta}^{\tau}\partial_{\beta}\bm{\Theta}_j(t), \partial_{\alpha}\bm{\Theta}_i(t)\bigg)_{L^2} \, \dd \tau \,\dd t\bigg\vert \leq \:C_{\tau}  \int_{t_{n-1}}^{t_n}\|\nabla \bm{\Theta}(t)\|_{L^2}^2 \,\dd t.
\end{equation}
Recalling that the values of $\nabla \mathbf{W}(t)\mathrm{e}^{\gamma t}+\tau\nabla(\bm{\phi}(t)-\mathbf{W}(t))\mathrm{e}^{\gamma t}$ and $\nabla \mathbf{W}(t)\mathrm{e}^{\gamma t}+\tau\nabla(\bm{\phi'}(t)-\mathbf{W}(t))\mathrm{e}^{\gamma t}$ belong to the convex set $\mathcal{M}_{\delta}$, and that $A_{i\alpha j \beta}$ is Lipschitz continuous on $\mathcal{M}_{\delta}$, we have
\begin{flalign}\label{for_FTC}
&\bigg\vert \int_{t_{n-1}}^{t_n} \mathsf{A}(\bm{\phi'}(t)\mathrm{e}^{\gamma t}; \bm{\theta}_{\phi'}(t), \bm{\dot\Theta}(t))-\mathsf{A}(\bm{\phi}(t)\mathrm{e}^{\gamma t}; \bm{\theta}_{\phi'}(t), \bm{\dot \Theta}(t))\,\dd t\bigg\vert &&\nonumber\\
&\leq  \int_{t_{n-1}}^{t_n} \big\vert\mathsf{A}(\bm{\phi'}(t)\mathrm{e}^{\gamma t}; \bm{\theta}_{\phi'}(t), \bm{\dot\Theta}(t))-\mathsf{A}(\bm{\phi}(t)\mathrm{e}^{\gamma t}; \bm{\theta}_{\phi'}(t), \bm{\dot \Theta}(t))\big\vert\, \dd t &&\nonumber\\
&\leq L_{\delta} \int_{t_{n-1}}^{t_n}\|\nabla \bm{\phi}(t)-\nabla\bm{\phi'}(t)\|_{L^2}\|\nabla \bm{\theta}_{\phi'}(t)\|_{L^{\infty}}\|\nabla\bm{\dot \Theta}(t)\|_{L^2}\,\dd t &&\nonumber\\
&\leq L_{\delta}C_0C_1 \int_{t_{n-1}}^{t_n} \|\nabla \mathbf{R}(t)\|_{L^2}h^{-1-\frac{d}{2}}\|\nabla\bm{\theta}_{\phi'}(t)\|_{L^2}\|\bm{\dot \Theta}(t)\|_{L^2} \,\dd t &&\nonumber\\
&\leq 2\gamma\int_{t_{n-1}}^{t_n}\|\bm{\dot \Theta}(t)\|_{L^2}^2 \dd t+C(\gamma)\max_{t\in I_n} \|\nabla \mathbf{R}(t)\|_{L^2}^2 h^{-2-d} \int_{t_{n-1}}^{t_n} \|\nabla\bm{\theta}_{\phi '}(t)\|_{L^2}^2\,\dd t, &&
\end{flalign}
where we have used the inverse inequalities (ii,a) and (ii,b), Young's inequality. To approximate the $\nabla \bm{\theta}_{\phi'}$ term, we apply the fundamental theorem of calculus as in the proof of (a) to have 
\begin{equation}\label{FTC}
\int_{t_{n-1}}^{t_n}\|\nabla \bm{\theta}_{\phi'}(t)\|_{L^2}^2 \,\dd t \leq 2k_n \|\nabla \bm{\theta}_{\phi'}(t_n^{-})\|_{L^2}^2 + 2C_0^2k_n^2h^{-2}\int_{t_{n-1}}^{t_n}\| \bm{\dot{\theta}}_{\phi'}(t)\|_{L^2}^2 \,\dd t. 
\end{equation}
Substituting (\ref{FTC}) into inequality (\ref{for_FTC}) gives
\begin{flalign}
&\bigg\vert \int_{t_{n-1}}^{t_n} \mathsf{A}(\bm{\phi'}(t)\mathrm{e}^{\gamma t}; \bm{\theta}_{\phi'}(t), \bm{\dot\Theta}(t))-\mathsf{A}(\bm{\phi}(t)\mathrm{e}^{\gamma t}; \bm{\theta}_{\phi'}(t), \bm{\dot \Theta}(t))\,\dd t\bigg\vert &&\nonumber\\
&\leq 2\gamma\int_{t_{n-1}}^{t_n}\|\bm{\dot \Theta}(t)\|_{L^2}^2 \,\dd t+C\max_{t\in I_n} \|\nabla \mathbf{R}(t)\|_{L^2}^2 h^{-2-d} k_n  \left(\|\nabla \bm{\theta}_{\phi'}(t_n^{-})\|_{L^2}^2+\int_{t_{n-1}}^{t_n}\| \bm{\dot{\theta}}_{\phi'}(t)\|_{L^2}^2 \,\dd t \right),&&
\end{flalign}
where we have used the assumption that $\mu_i k_i\leq  h^2$ for each $i=1,\ldots, N.$ Here $C$ is a generic positive constant.
Analogously, we obtain
\begin{flalign}
&\bigg\vert  \int_{t_{n-1}}^{t_n} \mathsf{A}(\bm{\phi'}(t)\mathrm{e}^{\gamma t}; \bm{\dot \theta}_{\phi'}(t), \bm{\Theta}(t))-\mathsf{A}(\bm{\phi}(t)\mathrm{e}^{\gamma t}; \bm{\dot \theta}_{\phi'}(t), \bm{\Theta}(t))\,\dd t\bigg\vert &&\nonumber\\
&\leq   \int_{t_{n-1}}^{t_n} \big\vert\mathsf{A}(\bm{\phi'}(t)\mathrm{e}^{\gamma t}; \bm{\dot \theta}_{\phi'}(t), \bm{\Theta}(t))-\mathsf{A}(\bm{\phi}(t)\mathrm{e}^{\gamma t}; \bm{\dot\theta}_{\phi'}(t), \bm{\Theta}(t))\big\vert\: \dd t &&\nonumber\\
&\leq L_{\delta} \int_{t_{n-1}}^{t_n}\|\nabla \bm{\phi}(t)-\nabla\bm{\phi'}(t)\|_{L^2}\|\nabla \bm{\dot \theta}_{\phi'}(t)\|_{L^{\infty}}\|\nabla\bm{\Theta}(t)\|_{L^2}\: \dd t &&\nonumber\\
&\leq L_{\delta}C_1 \int_{t_{n-1}}^{t_n} \|\nabla \mathbf{R}(t)\|_{L^2}h^{-\frac{d}{2}}\|\nabla\bm{\dot \theta}_{\phi'}(t)\|_{L^2}\|\nabla\bm{\Theta}(t)\|_{L^2}\: \dd t &&\nonumber\\
&\leq \tilde{C} \int_{t_{n-1}}^{t_n}\|\nabla\bm{\Theta}(t)\|_{L^2}^2\: \dd t+C\max_{t\in I_n} \|\nabla\mathbf{R}(t)\|_{L^2}^2 h^{-d-2}  \left(\int_{t_{n-1}}^{t_n}\| \bm{\dot{\theta}}_{\phi'}(t)\|_{L^2}^2\: \dd t \right).
\end{flalign}
Applying the Lipschitz continuity of $A_{i\alpha j \beta}$ again, we have
\begin{flalign}\label{Wterm}
& 2\bigg\vert  \int_{t_{n-1}}^{t_n} \big( \mathsf{A}(\bm{\phi}(t)\mathrm{e}^{\gamma t}; \mathbf{W}(t)-\Pi_k\mathbf{W}(t),\bm{\dot\Theta}(t))- \mathsf{A}(\bm{\phi'}(t)\mathrm{e}^{\gamma t}; \mathbf{W}(t)-\Pi_k\mathbf{W}(t),\bm{\dot\Theta}(t))\big) \,\dd t\bigg\vert &&\nonumber\\
&\leq 2  \int_{t_{n-1}}^{t_n} \bigg\vert \big( \mathsf{A}(\bm{\phi}(t)\mathrm{e}^{\gamma t}; \mathbf{W}(t)-\Pi_k\mathbf{W}(t),\bm{\dot\Theta}(t))- \mathsf{A}(\bm{\phi'}(t)\mathrm{e}^{\gamma t}; \mathbf{W}(t)-\Pi_k\mathbf{W}(t),\bm{\dot\Theta}(t))\big) \bigg\vert\,\dd t &&\nonumber\\
&\leq 2L_{\delta} \int_{t_{n-1}}^{t_n} \|\nabla\mathbf{R}(t)\|_{L^2}\|\nabla \mathbf{W}(t)-\nabla(\Pi_k\mathbf{W})(t)\|_{L^{\infty}}\|\nabla \bm{\dot \Theta}(t)\|_{L^2}\,\dd t &&\nonumber\\
&\leq 2 L_{\delta}C_0C_1h^{-\frac{d}{2}-1} \int_{t_{n-1}}^{t_n} \|\nabla\mathbf{R}(t)\|_{L^2}\|\nabla \mathbf{W}(t)-\nabla(\Pi_k\mathbf{W})(t)\|_{L^{2}}\|\bm{\dot\Theta}(t)\|_{L^2}\,\dd t &&\nonumber\\
&\leq  C(\gamma)h^{-d-2}\max_{t\in I_n}\|\nabla\mathbf{R}(t)\|_{L^2}^2\int_{t_{n-1}}^{t_n}\|\nabla \mathbf{W}(t)-\nabla(\Pi_k\mathbf{W})(t)\|_{L^{2}}^2\,\dd t+\gamma\int_{t_{n-1}}^{t_n} \|\bm{\dot \Theta}(t)\|_{L^2}^2 \,\dd t &&\nonumber\\
&\leq \gamma \int_{t_{n-1}}^{t_n} \|\bm{\dot \Theta}(t)\|_{L^2}^2\: \dd t+C(\gamma,\mathbf{W})h^{-d-2} \frac{k_n^{2(q_n+1)+1}}{q_n^{2(s-1)}}\max_{t\in I_n}\|\nabla\mathbf{R}(t)\|_{L^2}^2,  \quad\text{ for } \mathbf{W}\in W^{s, \infty}([0,T]; H_0^1),
\end{flalign}
where $C(\gamma,\mathbf{W})$ is a positive constant depending on both $\gamma$ and the nonlinear projection $\mathbf{W}.$ Next, we need to bound the term involving $(\partial_t A_{i\alpha j\beta}^{\tau}-\partial_t \tilde{A}_{i \alpha j\beta}^{\tau}).$ Recall that 
$$\partial_t A_{i\alpha j \beta}^{\tau}:=\partial_t A_{i\alpha j\beta}(\nabla \mathbf{W}(t)\mathrm{e}^{\gamma t} +\tau(\nabla \bm{\phi}(t)-\nabla\mathbf{W}(t))\mathrm{e}^{\gamma t})$$
and $$\partial_t \tilde{A}_{i \alpha j\beta}^{\tau}:=\partial_t A_{i\alpha j\beta}(\nabla \mathbf{W}(t)\mathrm{e}^{\gamma t} +\tau(\nabla \bm{\phi'}(t)-\nabla\mathbf{W}(t))\mathrm{e}^{\gamma t}).$$
By Taylor's theorem with an integral remainder, we have 
\begin{flalign*}
&\partial_t A_{i\alpha j\beta}^{\tau} -\partial_t \tilde{A}_{i \alpha j\beta}^{\tau} &&\\
&=\int_0^1 \sum_{k,\gamma, l,\delta=1}^d  \frac{\partial^2 }{\partial \eta_{k \gamma} \partial\eta_{l \delta}}A_{i\alpha j \beta}(\nabla \mathbf{W}(t)\mathrm{e}^{\gamma t} +\tau(\nabla \bm{\phi'}(t)-\nabla\mathbf{W}(t))\mathrm{e}^{\gamma t} +\tilde{\tau}\tau (\nabla \bm{\phi}-\nabla \bm{\phi'})\mathrm{e}^{\gamma t})&&\\
&\quad\times \partial_t \partial_ {\delta}(\mathbf{W}_{l}(t)\mathrm{e}^{\gamma t} +\tau( \bm{\phi'}_{l}(t)-\mathbf{W}_{l}(t))\mathrm{e}^{\gamma t} +\tilde{\tau}\tau (\bm{\phi}_{l}- \bm{\phi'}_{l})\mathrm{e}^{\gamma t}) \tau \partial_{\gamma} (\bm{\phi}_{k}(t)-\bm{\phi'}_k(t))\,\dd \tilde{\tau}.
\end{flalign*}
Since $A_{i\alpha j\beta}$ is sufficiently smooth (in particular, twice continuously differentiable), we can estimate the above difference term by 
\begin{flalign*}
|\partial_t A_{i\alpha j\beta}^{\tau} -\partial_t \tilde{A}_{i\alpha j\beta}^{\tau}| \leq & \:\hat{C} \|\nabla\mathbf{R}(t)\|_{L^{\infty}}\big( \|\nabla\mathbf{W}(t)\|_{L^{\infty}}+\|\nabla\mathbf{\dot{W}}(t)\|_{L^{\infty}}\big) &&\\
&+\hat{C} \|\nabla\mathbf{R}(t)\|_{L^{\infty}}\big( \| \nabla \bm{\phi}(t)-\nabla\bm{\phi'}(t)\|_{L^{\infty}}+\|\nabla \bm{\dot{\phi}}(t)-\nabla\bm{\dot{\phi}'}(t)\|_{L^{\infty}}\big)&&\\
&+\hat{C} \|\nabla\mathbf{R}(t)\|_{L^{\infty}}\big(\|\nabla\bm{\phi'}-\nabla \mathbf{W}(t)\|_{L^{\infty}}+\|\nabla \bm{\dot{\phi'}}(t)-\nabla\mathbf{\dot{W}}(t)\|_{L^{\infty}}\big).&&
\end{flalign*}
Similarly to the proof of Lemma \ref{bdlemma}, we can show that 
\begin{equation}\label{bd1}
\|\nabla\bm{\phi'}-\nabla \mathbf{W}(t)\|_{L^{\infty}}+\|\nabla \bm{\dot{\phi'}}(t)-\nabla\mathbf{\dot{W}}(t)\|_{L^{\infty}}\leq C_{\tau}.
\end{equation}
Property (iii, c) of the nonlinear projection $\mathbf{W}$ implies that 
\begin{equation}\label{bd2}
\|\nabla \mathbf{W}\|_{L^{\infty}}+\|\nabla \mathbf{\dot{W}}(t)\|_{L^{\infty}}\leq c_0+c_1.
\end{equation}
In the view of the triangle inequality, we have
\begin{align}\label{bd3}
&\| \nabla \bm{\phi}(t)-\nabla\bm{\phi'}(t)\|_{L^{\infty}}+\|\nabla \bm{\dot{\phi}}(t)-\nabla\bm{\dot{\phi}'}(t)\|_{L^{\infty}} \nonumber\\
&\leq  \|\nabla \bm{\phi}(t)-\nabla\mathbf{W}(t)\|_{L^{\infty}}+\|\nabla \bm{\phi'}(t)-\nabla\mathbf{W}(t)\|_{L^{\infty}}+\|\nabla \bm{\dot{\phi}}(t)-\nabla\mathbf{\dot{W}}(t)\|_{L^{\infty}}\nonumber\\
&\quad +\|\nabla \bm{\dot{\phi}'}(t)-\nabla\mathbf{\dot{W}}(t)\|_{L^{\infty}}\nonumber\\
& \leq 2C_{\tau}.
\end{align}
Combining (\ref{bd1})--(\ref{bd3}) and applying the inverse inequality (ii,b), Young's inequality, we obtain 
\begin{flalign}\label{partialA}
&\bigg\vert  \int_{t_{n-1}}^{t_n} \int_0^1 \sum_{i,\alpha,j,\beta=1}^d \bigg( (\partial_t A_{i\alpha j\beta}^{\tau} -\partial_t \tilde{A}_{i\alpha j\beta}^{\tau})\partial_{\beta}\bm{\theta}_{\phi',j}(t), \partial_{\alpha}\bm{\Theta}_i(t) \bigg)_{L^2}\:\dd \tau \,\dd t \bigg\vert &&\nonumber\\
&\leq  \tilde{C}_{\mathrm{Lip}} \int_{t_{n-1}}^{t_n} \| \nabla\mathbf{R}(t)\|_{L^{\infty}}\|\nabla \bm{\theta}_{\phi'}(t)\|_{L^2}\|\nabla \bm{\Theta}(t)\|_{L^2}\,\dd t && \nonumber\\
&\leq  \tilde{C}_{\mathrm{Lip}}C_1 h^{-\frac{d}{2}} \int_{t_{n-1}}^{t_n} \| \nabla\mathbf{R}(t)\|_{L^2}\|\nabla\bm{\theta}_{\phi'}(t)\|_{L^2}\|\nabla \bm{\Theta}(t)\|_{L^2} \,\dd t && \nonumber\\
&\leq  \tilde{C}  \int_{t_{n-1}}^{t_n} \|\nabla \bm{\Theta}(t)\|_{L^2}^2 \,\dd t+C h^{-d} \int_{t_{n-1}}^{t_n} \| \nabla\mathbf{R}(t)\|_{L^2}^2\|\nabla\bm{ \theta}_{\phi '}(t)\|_{L^2}^2 \,\dd t \nonumber&&\\
&\leq  \tilde{C} \int_{t_{n-1}}^{t_n} \|\nabla \bm{\Theta}(t)\|_{L^2}^2 \,\dd t+ 2C\max_{t\in I_n} \|\nabla\mathbf{R}(t)\|_{L^2}^2 h^{-d} k_n  \left(\|\nabla \bm{\theta}_{\phi'}(t_n^{-})\|_{L^2}^2+\int_{t_{n-1}}^{t_n}\| \bm{\dot{\theta}}_{\phi'}(t)\|_{L^2}^2\, \dd t \right),&&
\end{flalign}
where $\tilde{C}_{\mathrm{Lip}}= \hat{C}(3C_{\tau}+c_0+c_1)$ and $C$ is a generic positive constant. Combining the estimates (\ref{cauchy1})--(\ref{Wterm}) and (\ref{partialA}), we obtain 
\begin{flalign}
&\|\bm{\dot \Theta}(t_n^{-})\|_{L^2}^2+\gamma^2\|\bm{\Theta}(t_n^{-})\|_{L^2}^2+\mathsf{A}(\bm{\phi}(t_n^{-})\mathrm{e}^{\gamma t_n}; \bm{\Theta}(t_n^{-}),\bm{\Theta}(t_n^{-}))+\gamma\int_{t_{n-1}}^{t_n} \|\bm{\dot \Theta}(t)\|_{L^2}^2\,\dd t &&\nonumber\\
&\leq \mathsf{A}(\bm{\phi}(t_{n-1}^{-})\mathrm{e}^{\gamma t_{n-1}}; \bm{\Theta}(t_{n-1}^{-}),\bm{\Theta}(t_{n-1}^{-}))+\|\bm{\dot \Theta}(t_{n-1}^{-})\|_{L^2}^2+\gamma^2\|\bm{\Theta}(t_{n-1}^{-})\|_{L^2}^2+C\int_{t_{n-1}}^{t_n} \|\nabla \bm{\Theta}(t)\|_{L^2}^2 \,\dd t &&\nonumber\\
&\quad+C(h^{-d-2}+h^{-d})k_n \left(\|\nabla \bm{\theta}_{\phi'}(t_n^{-})\|_{L^2}^2+\int_{t_{n-1}}^{t_n}\| \bm{\dot{\theta}}_{\phi'}(t)\|_{L^2}^2 \,\dd t \right)\max_{t\in I_n} \|\nabla\mathbf{R}(t)\|_{L^2}^2 &&\nonumber\\
&\quad +Ch^{-d-2}\left( \int_{t_{n-1}}^{t_n} \|\bm{\dot{\theta}_{\phi'}}(t)\|_{L^2}^2 \, \dd t\right)\max_{t\in I_n}\|\nabla \mathbf{R}(t)\|_{L^2}^2+C(\gamma,\mathbf{W})h^{-d-2}\frac{k_n^{2q_n+3}}{q_n^{2(s-1)}} \max_{t\in I_n}\|\nabla\mathbf{R}(t)\|_{L^2}^2&&\nonumber\\
&\quad+\mathsf{A}(\bm{\phi'}(t_{n}^{-})\mathrm{e}^{\gamma t_{n}};\bm{\theta}_{\phi'}(t_{n}^{-}),\bm{\Theta}(t_{n}^{-}))-\mathsf{A}(\bm{\phi'}(t_{n-1}^{-})\mathrm{e}^{\gamma t_{n-1}};\bm{\theta}_{\phi'}(t_{n-1}^{-}),\bm{\Theta}(t_{n-1}^{-})) &&\nonumber\\
&\quad-\mathsf{A}(\bm{\phi}(t_{n}^{-})\mathrm{e}^{\gamma t_{n}};\bm{\theta}_{\phi'}(t_{n}^{-}),\bm{\Theta}(t_{n}^{-}))+\mathsf{A}(\bm{\phi}(t_{n-1}^{-})\mathrm{e}^{\gamma t_{n-1}};\bm{\theta}_{\phi'}(t_{n-1}^{-}),\bm{\Theta}(t_{n-1}^{-})).
\end{flalign}
Summing up over $n=1,\ldots,j$, we have
\begin{flalign}\label{contractionsum}
&\|\bm{\dot \Theta}(t_j^{-})\|_{L^2}^2+\gamma^2\|\bm{\Theta}(t_j^{-})\|_{L^2}^2+\mathsf{A}(\bm{\phi}(t_j^{-})\mathrm{e}^{\gamma t_j}; \bm{\Theta}(t_j^{-}),\bm{\Theta}(t_j^{-}))+\gamma\sum_{n=1}^j\int_{t_{n-1}}^{t_n} \|\bm{\dot \Theta}(t)\|_{L^2}^2\,\dd t &&\nonumber\\
&\leq C\sum_{n=1}^j\int_{t_{n-1}}^{t_n} \|\nabla \bm{\Theta}(t)\|_{L^2}^2 \,\dd t +C \left(h^{-d-2}+h^{-d}\right)\left(\sum_{n=1}^j k_n\|\nabla\bm{\theta}_{\phi'}(t_{n-1}^{-})\|_{L^2}^2\right) \max_{t\in I_n,1\leq n\leq j}\|\nabla\mathbf{R}(t)\|_{L^2}^2&&\nonumber\\
&\quad+C \left(h^{-d-2}+h^{-d}\right)\left(\sum_{n=1}^j k_n\int_{t_{n-1}}^{t_n}\|\bm{\dot{\theta}}_{\phi'}\|_{L^2}^2 \dd t\right) \max_{t\in I_n, 1\leq n\leq j}\|\nabla\mathbf{R}(t)\|_{L^2}^2&&\nonumber\\
&\quad+Ch^{-d-2}\left( \sum_{n=1}^j\int_{t_{n-1}}^{t_n} \|\bm{\dot{\theta}_{\phi'}}(t)\|_{L^2}^2 \, \dd t\right)\max_{t\in I_n}\|\nabla \mathbf{R}(t)\|_{L^2}^2+C(\gamma,\mathbf{W})h^{-d-2}\sum_{n=1}^j\frac{k_n^{2q_n+3}}{q_n^{2(s-1)}} \max_{t\in I_n}\|\nabla\mathbf{R}(t)\|_{L^2}^2&&\nonumber\\
&\quad +\mathsf{A}(\bm{\phi'}(t_{j}^{-})\mathrm{e}^{\gamma t_{j}};\bm{\theta}_{\phi'}(t_{j}^{-}),\bm{\Theta}(t_{j}^{-}))-\mathsf{A}(\bm{\phi}(t_j^{-})\mathrm{e}^{\gamma t_j};\bm{\theta}_{\phi'}(t_j^{-}),\bm{\Theta}(t_j^{-}))&&\nonumber\\
&\leq C\sum_{n=1}^j\int_{t_{n-1}}^{t_n} \|\nabla \bm{\Theta}(t)\|_{L^2}^2 \,\dd t +C (\mathbf{v})\left(h^{-d-2}+h^{-d}\right)\left(\sum_{n=1}^j k_n h^{2r+2}+\frac{k_n^{2q_n+1}}{q_n^{2(s-1)}}\right) \max_{t\in I_n, 1\leq n\leq j}\|\nabla\mathbf{R}(t)\|_{L^2}^2&&\nonumber\\
&\quad+C(\mathbf{v}) \left(h^{-d-2}+h^{-d}\right)\max_{1\leq n\leq j}k_n\left(\sum_{n=1}^j k_n h^{2r+2}+\frac{k_n^{2q_n+1}}{q_n^{2(s-1)}}\right) \max_{t\in I_n,1\leq n\leq j}\|\nabla\mathbf{R}(t)\|_{L^2}^2&&\nonumber\\
&\quad +Ch^{-d-2}\left( \sum_{n=1}^j\int_{t_{n-1}}^{t_n} \|\bm{\dot{\theta}_{\phi'}}(t)\|_{L^2}^2 \, \dd t\right)\max_{t\in I_n}\|\nabla \mathbf{R}(t)\|_{L^2}^2+C(\gamma,\mathbf{W})h^{-d-2}\sum_{n=1}^j\frac{k_n^{2q_n+3}}{q_n^{2(s-1)}} \max_{t\in I_n, 1\leq n\leq j}\|\nabla\mathbf{R}(t)\|_{L^2}^2 &&\nonumber\\
&\quad+\mathsf{A}(\bm{\phi'}(t_{j}^{-})\mathrm{e}^{\gamma t_{j}};\bm{\theta}_{\phi'}(t_{j}^{-}),\bm{\Theta}(t_{j}^{-}))-\mathsf{A}(\bm{\phi}(t_{j}^{-})\mathrm{e}^{\gamma t_{j}};\bm{\theta}_{\phi'}(t_{j}^{-}),\bm{\Theta}(t_{j}^{-}))&&\nonumber\\
&\leq C\sum_{n=1}^j\int_{t_{n-1}}^{t_n} \|\nabla \bm{\Theta}(t)\|_{L^2}^2 \dd t +C(\mathbf{v}) h^{-d-2}\left(\sum_{n=1}^j k_n h^{2r+2}+\frac{k_n^{2q_n+1}}{q_n^{2(s-1)}}\right) \max_{t\in I_n, 1\leq n\leq j}\|\nabla\mathbf{R}(t)\|_{L^2}^2&&\nonumber\\
&\quad+\mathsf{A}(\bm{\phi'}(t_{j}^{-})\mathrm{e}^{\gamma t_{j}};\bm{\theta}_{\phi'}(t_{j}^{-}),\bm{\Theta}(t_{j}^{-}))-\mathsf{A}(\bm{\phi}(t_j^{-})\mathrm{e}^{\gamma t_j};\bm{\theta}_{\phi'}(t_j^{-}),\bm{\Theta}(t_j^{-})),
\end{flalign}
where $C$ is a generic positive constant and $C(\mathbf{v})$ is a positive constant depending on the exact solution $\mathbf{v}.$ These constants may change from line to line. Using the Lipschitz continuity of $A(\cdot; \bm{\theta}_{\phi'}(t_j^{-}),\bm{\Theta}(t_j^{-}))$ and the inverse inequality (ii,b), we obtain
\begin{align}\label{lipestimate}
&\vert A(\bm{\phi'}(t_j^{-})\mathrm{e}^{\gamma t_j}; \bm{\theta}_{\phi'}(t_j^{-}),\bm{\Theta}(t_j^{-}))-A(\bm{\phi}(t_j^{-})\mathrm{e}^{\gamma t_j}; \bm{\theta}_{\phi'}(t_j^{-}),\bm{\Theta}(t_j^{-}))\vert \nonumber\\
&\leq L_{\delta} \|\nabla\mathbf{R}(t_j^{-})\|_{L^2} \|\nabla \bm{\theta}_{\phi'}(t_j^{-})\|_{L^{\infty}} \|\nabla\bm{\Theta}(t_j^{-}) \|_{L^2}&&\nonumber\\
&\leq  C(\mathbf{v})h^{-\frac{d}{2}}\left(\sum_{n=1}^j k_n h^{2r+2}+\frac{k_n^{2q_n+1}}{q_n^{2(s-1)}} \right)^{\frac{1}{2}}\max_{t\in I_n,1\leq n\leq j}\|\nabla\mathbf{R}(t)\|_{L^2}\|\nabla \bm{\Theta}(t_j^{-})\|_{L^2}.
\end{align}
Combining the estimates (\ref{contractionsum}) and (\ref{lipestimate}) and applying the assumption (S2b) to $\mathsf{A}(\bm{\phi}(t_j^{-})\mathrm{e}^{\gamma t_j}; \bm{\Theta}(t_j^{-}),\bm{\Theta}(t_j^{-}))$ on the left-hand-side of the resulting inequality yield 
\begin{flalign}\label{contraction_1}
&\|\bm{\dot \Theta}(t_j^{-})\|_{L^2}^2+\gamma^2\|\bm{\Theta}(t_j^{-})\|_{L^2}^2+M_1\|\nabla \bm{\Theta}(t_j^{-})\|^2+\gamma\sum_{n=1}^j\int_{t_{n-1}}^{t_n} \|\bm{\dot \Theta}(t)\|_{L^2}^2\,\dd t &&\nonumber\\
&\leq C\sum_{n=1}^j\int_{t_{n-1}}^{t_n} \|\nabla \bm{\Theta}(t)\|_{L^2}^2 \dd t +C(\mathbf{v})h^{-d-2}\left(\sum_{n=1}^j k_n h^{2r+2}+\frac{k_n^{2q_n+1}}{q_n^{2(s-1)}}\right) \max_{t\in I_n,\:1\leq n\leq j}\|\nabla\mathbf{R}(t)\|_{L^2}^2&&\nonumber\\
& \quad+C(\mathbf{v})h^{-\frac{d}{2}}\left(\sum_{n=1}^j k_i h^{2r+2}+\frac{k_i^{2q_i+1}}{q_i^{2(s-1)}} \right)^{\frac{1}{2}}\max_{t\in I_n,\:1\leq n\leq j}\|\nabla\mathbf{R}(t)\|_{L^2}\|\nabla \bm{\Theta}(t_j^{-})\|_{L^2}.&& 
\end{flalign}
By applying Young's inequality on the right-hand side of (\ref{contraction_1}), we have 
\begin{flalign}\label{contraction_2}
&\|\bm{\dot \Theta}(t_j^{-})\|_{L^2}^2+\gamma^2\|\bm{\Theta}(t_j^{-})\|_{L^2}^2+M_1\|\nabla \bm{\Theta}(t_j^{-})\|^2+\gamma\sum_{n=1}^j\int_{t_{n-1}}^{t_n} \|\bm{\dot \Theta}(t)\|_{L^2}^2\,\dd t &&\nonumber\\
&\leq C\sum_{n=1}^j\int_{t_{n-1}}^{t_n} \|\nabla \bm{\Theta}(t)\|_{L^2}^2 \dd t +C(\mathbf{v}) h^{-d-2}\left(\sum_{n=1}^j k_n h^{2r+2}+\frac{k_n^{2q_n+1}}{q_n^{2(s-1)}}\right) \max_{t\in I_n,\:1\leq n\leq j}\|\nabla\mathbf{R}(t)\|_{L^2}^2 &&\nonumber\\
&\quad + C(M_1,\mathbf{v})h^{-d}\left(\sum_{n=1}^j k_i h^{2r+2}+\frac{k_i^{2q_i+1}}{q_i^{2(s-1)}} \right)\max_{t\in I_n,\:1\leq n\leq j}\|\nabla\mathbf{R}(t)\|_{L^2}^2+\frac{M_1}{2}\|\nabla \bm{\Theta}(t_j^{-})\|_{L^2}^2, 
\end{flalign}
where $C(M_1, \mathbf{v})$ is a constant depending on $M_1$ and the exact solution $\mathbf{v}$. This implies that 
\begin{flalign}\label{contraction_3}
&\|\bm{\dot \Theta}(t_j^{-})\|_{L^2}^2+\gamma^2\|\bm{\Theta}(t_j^{-})\|_{L^2}^2+\frac{M_1}{2}\|\nabla \bm{\Theta}(t_j^{-})\|^2+\gamma\sum_{n=1}^j\int_{t_{n-1}}^{t_n} \|\bm{\dot \Theta}(t)\|_{L^2}^2\,\dd t&& \nonumber\\
&\leq C\sum_{n=1}^j\int_{t_{n-1}}^{t_n} \|\nabla \bm{\Theta}(t)\|_{L^2}^2 \,\dd t +C(\mathbf{v}) h^{-d-2}\left(\sum_{n=1}^j k_n h^{2r+2}+\frac{k_n^{2q_n+1}}{q_n^{2(s-1)}}\right) \max_{t\in I_n,\:1\leq n\leq j}\|\nabla\mathbf{R}(t)\|_{L^2}^2 &&\nonumber\\
&\quad+C(M_1,\mathbf{v})h^{-d}\left(\sum_{n=1}^j k_i h^{2r+2}+\frac{k_i^{2q_i+1}}{q_i^{2(s-1)}} \right)\max_{t\in I_n,\:1\leq n\leq j}\|\nabla\mathbf{R}(t)\|_{L^2}^2 &&\nonumber\\ 
&\leq C\sum_{n=1}^j\int_{t_{n-1}}^{t_n} \|\nabla \bm{\Theta}(t)\|_{L^2}^2 \,\dd t +C(\mathbf{v}) h^{-d-2}\left(\sum_{n=1}^j k_n h^{2r+2}+\frac{k_n^{2q_n+1}}{q_n^{2(s-1)}}\right) \max_{t\in I_n,\:1\leq n\leq j}\|\nabla\mathbf{R}(t)\|_{L^2}^2. &&
\end{flalign}
By an analogous application of the discrete Gr\"{o}nwall lemma as in the proof of (a), we can deduce that, for $k_n$ sufficiently small for each $n=1,\ldots,j$,
\begin{align}\label{aftergronwall}
&\|\bm{\dot \Theta}(t_j^{-})\|_{L^2}^2+\|\bm{\Theta}(t_j^{-})\|_{L^2}^2+\|\nabla \bm{\Theta}(t_j^{-})\|^2+\sum_{n=1}^j\int_{t_{n-1}}^{t_n} \|\bm{\dot \Theta}(t)\|_{L^2}^2\,\dd t\nonumber \\
&\leq   \tilde{C}(\mathbf{v}) h^{-d-2}\left(\sum_{n=1}^j k_n h^{2r+2}+\frac{k_n^{2q_n+1}}{q_n^{2(s-1)}}\right) \max_{t\in I_n,\:1\leq n\leq j}\|\nabla\mathbf{R}(t)\|_{L^2}^2,
\end{align}
where $\tilde{C}(\mathbf{v})$ is a positive constant depending on $\mathbf{v}$ which may vary from line to line.
By the fundamental theorem of calculus and the triangle inequality, we have for each $t\in I_n$, with $n=1,\ldots N$,
\begin{align*}
\|\nabla \bm{\Theta}(t)\|_{L^2}&\leq \|\nabla \bm{\Theta}(t_n^{-})\|_{L^2}+\int_{t_{n-1}}^{t_n}\|\nabla\bm{\dot \Theta}(t)\|_{L^2}\, \dd t\\
&\leq  \|\nabla\bm{\Theta}(t_n^{-})\|_{L^2}+C_0h^{-1}\int_{t_{n-1}}^{t_n}\|\bm{\dot \Theta}(t)\|_{L^2}\,\dd t\\
&\leq \|\nabla\bm{\Theta}(t_n^{-})\|_{L^2}+C_0h^{-1}k_n^{-\frac{1}{2}}\left(\int_{t_{n-1}}^{t_n}\|\bm{\dot \Theta}(t)\|_{L^2}^2\,\dd t\right)^{\frac{1}{2}}.
\end{align*}
This implies that, for each $t\in I_n$ with $1\leq n\leq N$, 
\begin{align}\label{h1semin}
\|\nabla \bm{\Theta}(t)\|_{L^2}\leq &
\:(1+C_0 h^{-1} k_n^{-\frac{1}{2}})\tilde{C}(\mathbf{v})^{\frac{1}{2}} h^{-\frac{d}{2}-1}\left(\sum_{i=1}^N k_i h^{2r+2}+\frac{k_i^{2q_i+1}}{q_i^{2(s-1)}} \right)^{\frac{1}{2}}\max_{t\in I_n, \:1\leq n\leq N}\|\nabla\mathbf{R}(t)\|_{L^2}\nonumber\\
\leq &\:\tilde{C}_1(\mathbf{v}) h^{-\frac{d}{2}-1}\left(\sum_{i=1}^N k_i h^{2r+2}+\frac{k_i^{2q_i+1}}{q_i^{2(s-1)}} \right)^{\frac{1}{2}}\max_{t\in I_n, \:1\leq n\leq N}\|\nabla\mathbf{R}(t)\|_{L^2},
\end{align}
where the last inequality follows from the assumption that $\mu_i k_i\leq h^2$ for each $1\leq i\leq N.$ Analogously, we have 
\begin{align}\label{l2}
\|\bm{\Theta}(t)\|_{L^2} \leq & \: \|\bm{\Theta}(t_n)\|_{L^2} +k_n ^{\frac{1}{2}}\left( \int_{t_{n-1}}^{t_n} \|\bm{\dot{\Theta}}(t)\|_{L^2}^2 \,\dd t\right)^{\frac{1}{2}} \nonumber\\
{}\leq & \:\tilde{C}_2(\mathbf{v})h^{-\frac{d}{2}-1}\left(\sum_{i=1}^N k_i h^{2r+2}+\frac{k_i^{2q_i+1}}{q_i^{2(s-1)}} \right)^{\frac{1}{2}}\max_{t\in I_n, \:1\leq n\leq N}\|\nabla\mathbf{R}(t)\|_{L^2},
\end{align}
and 
\begin{align}\label{firstderivative}
\|\bm{\dot{\Theta}}(t)\|_{L^2} \leq & \: \|\bm{\dot{\Theta}}(t_n)\|_{L^2} +  \int_{t_{n-1}}^{t_n} \|\bm{\ddot{\Theta}}(t)\|_{L^2} \dd t\nonumber\\
{}\leq & \: \|\bm{\dot{\Theta}}(t_n)\|_{L^2} +C_2k_n^{-\frac{1}{2}}  \left(\int_{t_{n-1}}^{t_n} \|\bm{\dot{\Theta}}(t)\|_{L^2}^2\, \dd t \right)^{\frac{1}{2}}\nonumber\\
{}\leq & \:\tilde{C}_3(\mathbf{v})h^{-\frac{d}{2}-1}k_n^{-\frac{1}{2}}\left(\sum_{i=1}^N k_i h^{2r+2}+\frac{k_i^{2q_i+1}}{q_i^{2(s-1)}} \right)^{\frac{1}{2}}\max_{t\in I_n, \:1\leq n\leq N}\|\nabla\mathbf{R}(t)\|_{L^2},
\end{align}
Summing up (\ref{h1semin})--(\ref{firstderivative}) and taking maximum on the left-hand for $t\in I_n$, with $1\leq  n\leq N$, we have 
\begin{align}
&\max_{t\in I_n, 1\leq n \leq N}\left(\|\bm{\Theta}(t)\|_{H^1}+\|\bm{\dot{\Theta}}(t)\|_{L^2} \right)\nonumber\\
&\leq \widehat{C}(\mathbf{v}) h^{-\frac{d}{2}-1} \max_{1\leq n\leq N} k_n^{-\frac{1}{2}}\left(\sum_{i=1}^N k_i h^{2r+2}+\frac{k_i^{2q_i+1}}{q_i^{2(s-1)}} \right)^{\frac{1}{2}}\max_{t\in I_n,1\leq n\leq N}\|\nabla\mathbf{R}(t)\|_{L^2}.
\end{align}
By choosing the mesh size $h$ and time steps $\{k_i\}_{i=1}^N<1$ for each $i=1,2,\ldots, N$ small enough, and $r$ and $\{q_i\}_{i=1}^N$ for each $i=1,2,\ldots, N$ large enough such that 
\begin{equation}\label{contractionconstant}
\widehat{C}(\mathbf{v}) h^{-\frac{d}{2}-1} \max_{1\leq n\leq N} k_n^{-\frac{1}{2}}\left(\sum_{i=1}^N k_i h^{2r+2}+\frac{k_i^{2q_i+1}}{q_i^{2(s-1)}} \right)^{\frac{1}{2}}<1, 
\end{equation}
we obtain 
\begin{equation}
\max_{t\in I_n, 1 \leq n \leq N}\left( \|\bm{\Theta}(t)\|_{H^1}+\| \bm{\dot{\Theta}}(t)\|_{L^2}\right)< \max_{t\in I_n, 1 \leq n \leq N}\left( \|\mathbf{R}(t)\|_{H^1}+\|\dot{\mathbf{R}}(t)\|_{L^2}\right).
\end{equation}
Indeed, the inequality (\ref{contractionconstant}) follows from our assumptions that $r>\frac{d}{2}+1$, $k_i^{q_i-\frac{1}{2}}=o(h^{1+\frac{d}{2}})$ for each $i=1,2,\ldots,N$. 
Therefore, by Banach's fixed point theorem, $\mathbf{v}_{\mathrm{DG}}=\mathbf{v}_{\phi}$ is the unique solution to (\ref{discrete_vform}). By the triangle inequality, (\ref{prop1}), (\ref{prop2}) and property (iii,b) of the nonlinear projection $\mathbf{W}$, we have 
\begin{align*}
&\|\mathbf{v}_{\mathrm{DG}}(t_j^{-})-\mathbf{v}(t_j^{-})\|_{L^2}+ \|\mathbf{\dot v}_{\mathrm{DG}}(t_j^{-})-\mathbf{\dot v}(t_j^{-})\|_{L^2}\\
&\leq  \|\bm{\theta}(t_j^{-})\|_{L^2}+\|\bm{\dot \theta}(t_j^{-})\|_{L^2}+\|\mathbf{W}(t_j^{-})-\mathbf{v}(t_j^{-})\|_{L^2}+\|\mathbf{\dot W}(t_j^{-})-\mathbf{\dot v}(t_j^{-})\|_{L^2}\\
&\leq  C_{\ast}(\mathbf{v})\bigg(\sum_{n=1}^j k_n h^{2r+2}+\frac{k_n^{2q_n+1}}{q_n^{2(s-1)}} \bigg)^{\frac{1}{2}} + 2 \tilde{C}_r(\mathbf{v}) h^{r+1} \quad \text{ (by (iii,b))}\\
&\leq  C(\mathbf{v})\bigg(h^{2r+2}+\sum_{n=1}^j \frac{k_n^{2q_n+1}}{q_n^{2(s-1)}} \bigg)^{\frac{1}{2}}.
\end{align*}
\section{Numerical experiments}\label{numerical}
In this section, we show some numerical experiments on a simple version of the nonlinear elastodynamics equation to verify the error estimates proved in Section \ref{elasto_convergence}.

\subsection{Numerical results for a one-dimensional nonlinear elastodynamics problem}\label{test1}
We consider the one-dimensional nonlinear equation
\begin{equation}
\ddot{u}(x,t)+2\gamma\dot{u}(x,t)+\gamma^2 u(x,t)-\partial_{x}[S(\partial_x u(x,t))]= f(x,t) \quad \mbox{ in } (0,1)\times (0,T],
\end{equation}
\begin{equation}\label{wavebc}
u(0, t) = u(1,t)= 0\quad \mbox{for all } t\in (0,T],
\end{equation}
\begin{equation}\label{waveic}
u(x,0)= u_0(x), \quad \dot{u}(x,0)= u_1(x).
\end{equation}
We take $S(\partial_x u)=\frac{1}{3}\partial_x u^3$ and the time interval to be $I=(0,T]$ with $T=1$. Let $u_0$, $u_1$ and $f$ be chosen such that the exact solution is 
$$u(x,t)=\sin(\sqrt{2}\pi t)\sin(\pi x).$$
That is, $u_0(x)\equiv 0$, $u_1(x)= \sqrt{2}\pi \sin(\pi x)$, and 
$$f(x,t) =[(-2\pi^2+\gamma^2)\sin(\sqrt{2}\pi t)+2\sqrt{2}\gamma\pi \cos(\sqrt{2}\pi t)]\sin(\pi x)+\pi^4\sin^3(\sqrt{2}\pi t)\cos^2(\pi x)\sin(\pi x).$$
We first discretize the problem in the spatial direction using continuous piecewise polynomials of degree $p\geq 1$. Let $\mathcal{V}_h$ be the finite element function space with $h$ being the spatial discretization parameter. The numerical approximation of the nonlinear wave-type equation following a Picard-type linearization in the nonlinear term is the following: find $u_h\in \mathcal{V}_h$ such that 
$$\int_{\Omega} \ddot{u}_h\cdot v_h\: \mathrm{d} x + \int_{\Omega} 2\gamma\dot{u}_h\cdot v_h\: \mathrm{d} x+\int_{\Omega} \gamma^2 u_h\cdot v_h\: \mathrm{d} x+ \frac{1}{3} \int_{\Omega}(\partial_x u_h^{\ast})^2\partial_x u_h\cdot \partial_x v_h\: \mathrm{d} x = \int_{\Omega} f\cdot v_h\: \mathrm{d} x,$$
for all $v_h\in\mathcal{V}_h$.
Here we assume that $\partial_x u_h^{\ast}$ is known at each time step $I_{n}$ either as an initial guess by using $u_h$ over the previous time interval, or as a previous iterate in the Picard iteration.
Now the problem results in the following second-order differential system for the nodal displacement $\mathbf{U}(t)$:
$$\begin{cases} \tilde{M} \mathbf{\ddot{U}}(t)+2\gamma \tilde{M}{\mathbf{\dot{U}}}(t)+\gamma^2\tilde{M}\mathbf{U}(t)+\frac{1}{3}\tilde{K}(t)\mathbf{U}(t)= \mathbf{F}(t), \quad t\in (0,T], &\\
\mathbf{\dot{U}}(0)=\mathbf{U}_1,\quad \mathbf{U}(0)=\mathbf{U}_0,
\end{cases}$$
where $\mathbf{U}_0=[0, \ldots,0]^{\mathrm{T}}\in \mathbb{R}^{\hat{d}}$ and $\mathbf{U}_1$ is the $\hat{d}$-vector corresponding to $u_1$ at the grid points, $\mathbf{\ddot{U}}(t)$ (respectively $\mathbf{\dot{U}}(t)$) represents the vector of nodal acceleration (respectively velocity) and $\mathbf{F}(t)$ is the vector of externally applied loads. $\tilde{M}$ is the mass matrix which is defined as $$\tilde{M}_{ij}:=\int_0^1 \psi_i(x)\psi_j(x)\, \dd x.$$ The time-dependent stiffness matrix $\tilde{K}(t)$ is defined as
$$\tilde{K}_{ij}(t):=\int_0^1 (\partial_x u_h^{\ast}(t))^2 \partial_x\psi_i(x)\partial_x\psi_j(x) \, \dd x,$$
where $\{\psi_i\}_{i=1}^{\hat{d}}$ are the basis functions in the spatial direction.

Multiplying the above algebraic formulation by $\tilde{M}^{-\frac{1}{2}}$ and setting $\mathbf{Z}(t)=\tilde{M}^{\frac{1}{2}}\mathbf{U}(t)$, we obtain
\begin{equation}\label{algebra_eq_nonlinear}
\mathbf{\ddot{Z}}(t) + L\mathbf{\dot{Z}}(t)+ K_0\mathbf{Z}(t) + K(t)\mathbf{Z}(t)=\mathbf{G}(t), \quad t\in (0,T],
\end{equation}
\begin{equation}\label{algebra_ic_nonlinear}
\mathbf{\dot{Z}}(0)= \tilde{M}^{\frac{1}{2}}\mathbf{U}_1, \quad \mathbf{Z}(0)=\tilde{M}^{\frac{1}{2}} \mathbf{U}_0.
\end{equation}
Here 
\begin{align*}
L= 2\gamma\mathrm{Id},&\quad K_0 = \gamma^2\mathrm{Id},\\
K(t)= \frac{1}{3}\tilde{M}^{-\frac{1}{2}}\tilde{K}(t)\tilde{M}^{-\frac{1}{2}}, &\quad \mathbf{G}(t)= \tilde{M}^{-\frac{1}{2}}\mathbf{F}(t).
\end{align*}
Note that both $\tilde{K}(t)$ and $K(t)$ are time-dependent. We subdivided $[0,T)$ into $N$ subintervals $I_n$, for $n=1,\ldots,N,$ of uniform length $k$. We assume that the polynomial degree in time is constant at each time step. That is, $q_1=\cdots=q_N\geq 2.$ If we consider the time integration on a generic time interval $I_n$ for each $n=1,\ldots N$, our discontinuous-in-time formulation reads as:
find $\mathbf{Z}\in\mathcal{V}_{kh}^{q_n}$ such that 
\begin{align}\label{weak_form_t2}
&\left(\mathbf{\ddot{Z}}(t), \mathbf{\dot{v}}\right)_{L^2(I_n)}+\left(L\mathbf{\dot{Z}}(t),\mathbf{\dot{v}}\right)_{L^2(I_n)}+\left(K_0\mathbf{Z}(t),\mathbf{\dot{v}}\right)_{L^2(I_n)}+\left(K(t)\mathbf{Z}(t),\mathbf{\dot{v}}\right)_{L^2(I_n)}\nonumber\\
&\quad+\mathbf{\dot{Z}}(t_{n-1}^{+})\cdot \mathbf{\dot{v}}(t_{n-1}^{+})+K_0\mathbf{Z}(t_{n-1}^{+})\cdot\mathbf{v}(t_{n-1}^{+})+K(t_{n-1}^{+})\mathbf{Z}(t_{n-1}^{+})\cdot\mathbf{v}(t_{n-1}^{+})\nonumber\\
&=\left( \mathbf{G}(t), \mathbf{\dot{v}}\right)_{L^2(I_n)}+\mathbf{\dot{Z}}(t_{n-1}^{-})\cdot\mathbf{\dot{v}}(t_{n-1}^{+})+K_0\mathbf{Z}(t_{n-1}^{-})\cdot\mathbf{v}(t_{n-1}^{+})+K(t_{n-1}^{-})\mathbf{Z}(t_{n-1}^{-})\cdot\mathbf{v}(t_{n-1}^{+}),
\end{align}
for all $\mathbf{v}\in \mathcal{V}_{kh}^{q_n}$, where on the right-hand side the values $\mathbf{\dot{Z}}(t_{n-1}^{-})$ and $\mathbf{Z}(t_{n-1}^{-})$ computed for $I_{n-1}$ are used as initial conditions for the current time interval. For $I_1$, we set $\mathbf{\dot{Z}}(t_0^{-})=\mathbf{\dot{Z}}(0)$ and $\mathbf{Z}(t_0^{-})=\mathbf{Z}(0).$ Focusing on the generic time interval $I_n$, we introduce the basis functions in the time direction $\{\phi^{j}(t)\}_{j=1}^{q_n+1}$ for the polynomial space $\mathbb{P}^{q_n}(I_n)$ and define $D=\hat{d}(q_n+1),$ the dimension of the local finite element space $\mathcal{V}_{kh}^{q_n}.$ We also introduce the vectorial basis $\{\Phi_{m}^j(t)\}_{m=1,\ldots,\hat{d}}^{j=1,\ldots, q_n+1}$, where $\Phi_{m}^j(t)$ is the $\hat{d}$-dimensional vector whose $m$-th component is $\phi^j(t)$ and the other components are zero. We write 
\begin{equation}\label{basis_decom_in_t_nonlinear}
\mathbf{Z}(t)=\sum_{m=1}^{\hat{d}} \sum_{j=1}^{q_n+1} \alpha_{m}^j \Phi_m^j(t),
\end{equation}
where $\alpha_m^j\in\mathbb{R}$ for $m=1,\ldots, \hat{d}$, $j=1,\ldots,q_n+1$.
By choosing $\mathbf{v}(t)= \Phi_{m}^j(t)$ for each $m=1,\ldots,\hat{d}$, $j=1,\ldots, q_n+1$, we obtain the following algebraic system
\begin{equation}\label{final_algebraic_nonlinear}
\mathbf{A}\mathbf{z}=\mathbf{b},
\end{equation}
where $\mathbf{z}\in\mathbb{R}^{D}=\mathbb{R}^{(q_n+1)\hat{d}}$ is the solution vector (whose entries are the values of $\alpha_m^j$); $\mathbf{b}\in\mathbb{R}^D$ corresponds to the right-hand side, which is given componentwise as 
\begin{equation}\label{rhsb_nonlinear}
\mathbf{b}_m^j =  \left( \mathbf{G}(t), \dot{\Phi}_m^j\right)_{L^2(I_n)}+\mathbf{\dot{Z}}(t_{n-1}^{-})\cdot\dot{\Phi}_m^j(t_{n-1}^{+})+K_0\mathbf{Z}(t_{n-1}^{-})\cdot\Phi_m^j(t_{n-1}^{+})+K(t_{n-1}^{-})\mathbf{Z}(t_{n-1}^{-})\cdot\Phi_m^j(t_{n-1}^{+}), 
\end{equation}
for $m=1,\ldots,\hat{d}, j=1,\ldots,q_n+1.$
$\mathbf{A}$ is the local stiffness matrix defined with its structure being discussed below.  For $l,j=1,\ldots,q_n+1$, 
\begin{align*}
&M_{lj}^1=\left(\ddot{\phi}^j, \dot{\phi}^l\right)_{L^2(I_n)},\quad M_{lj}^2=\left(\dot{\phi}^j, \dot{\phi}^l\right)_{L^2(I_n)},\quad M_{lj}^3=\left(\phi^j, \dot{\phi}^l\right)_{L^2(I_n)},\\
& \tilde{M}_{lj}^3=\left(K(t)\phi^j, \dot{\phi}^l\right)_{L^2(I_n)},\quad M_{lj}^4=\dot{\phi}^j(t_{n-1}^{+})\cdot \dot{\phi}^l(t_{n-1}^{+}), \\
&M_{lj}^5=\phi^j(t_{n-1}^{+})\cdot \phi^l(t_{n-1}^{+}), \quad \tilde{M}_{lj}^5=K(t_{n-1}^{+})\phi^j(t_{n-1}^{+})\cdot \phi^l(t_{n-1}^{+}).
\end{align*}
Setting 
\begin{align*}
M= M^1+M^4, &\quad B_{ij}=  L_{ij} M^2+{K_0}_{ij}(M^3+M^5)+(\tilde{M}^3+\tilde{M}^5),
\end{align*}
with $M, B_{ij}\in\mathbb{R}^{(q_n+1)\times (q_n+1)}$ for any $i,j=1,\ldots,\hat{d}$, we can rewrite the matrix $\mathbf{A}$ as 
$$\mathbf{A}= \begin{bmatrix}
M & 0 & 0  &\cdots & 0  \\
0 & M & 0 &\cdots & 0 \\
\vdots& \ddots  & \ddots &\ddots& \vdots & \\
0  & 0 & 0 &\cdots & M
\end{bmatrix}+\begin{bmatrix}
B_{1,1} & B_{1,2} & \cdots & B_{1,\hat{d}}  \\
B_{2,1} & B_{2,2}& \cdots & B_{2,\hat{d}} \\
\vdots& \ddots & \ddots  & \vdots & \\
B_{\hat{d},1} & B_{\hat{d},2}  &\cdots & B_{\hat{d},\hat{d}}
\end{bmatrix}.$$
For each time interval $I_n=(t_{n-1}, t_n]$, we use the following shifted Legendre polynomials $\{\phi_i\}$ as the basis polynomial for the discontinuous Galerkin formulation:
\begin{align*}
& \phi^1(t)=1, \quad \phi^2(t)=\frac{2(t-t_{n-1}^{+})}{k_n}-1, \quad \phi^3(t)=\frac{6(t-t_{n-1}^{+})^2}{k_n^2}-\frac{6(t-t_{n-1}^{+})}{k_n}+1 \\
{}& \phi^4(t)=\frac{20(t-t_{n-1}^{+})^3}{k_n^3}-\frac{30(t-t_{n-1}^{+})^2}{k_n^2}+\frac{12(t-t_{n-1}^{+})}{k_n}-1\\
{} & \phi^5(t)=\frac{70(t-t_{n-1}^{+})^4}{k_n^4}-\frac{140(t-t_{n-1}^{+})^3}{k_n^3}+\frac{90(t-t_{n-1}^{+})^2}{k_n^2}-\frac{20(t-t_{n-1}^{+})}{k_n}+1.
\end{align*}
In order to compute the time-dependent matrices $\tilde{M}_3(t)$ and $\tilde{M}_5(t)$, which also depend on the gradient of the solution, we apply a \emph{Picard iteration} at each time interval. We set the maximal number of Picard iterations to be $30$ at each time step and the tolerance to be $1\mathrm{e}-10$. The details of the algorithm are summarized below.
\begin{algorithm}[ht]
\caption{Iterative Algorithm (Multiple Picard iterations at each time interval)} 
\begin{algorithmic}
  \STATE \emph{Initialization}: 
  $\partial_x u_h^{\ast}= \partial_x u_0$ and  $$[\tilde{K}^0]_{ij}=\int_0^1 (\partial_x u_0)^2 \partial_x\psi_i(x)\partial_x \psi_j(x) \: \dd x.$$
  
  \STATE \emph{Iteration}:
  On each interval $I_n=(t_{n-1}, t_n]$ for $n=1,2,\ldots, N$, we solve  
 $$\tilde{M}\mathbf{\ddot{U}}(t)+2\gamma \tilde{M}\mathbf{\dot{U}}(t)+\gamma^2\tilde{M}\mathbf{U}(t)+\frac{1}{3}\tilde{K}^n(t)\mathbf{U}(t)=\mathbf{F}(t)$$
iteratively (using Picard iterations) by applying the discontinuous-in-time integration. Here $$[\tilde{K}_0^n]_{ij}=\int_0^1 [\partial_x u_{\mathrm{DG}}^{n-1}(t)]^2 \partial_x\psi_i(x)\partial_x \psi_j(x) \: \dd x,$$ where $u_{\mathrm{DG}}^{n-1}(t)$ is the solution we obtained from the previous time interval $I_{n-1}$.
 $$[\tilde{K}_k^n]_{ij}=\int_0^1 [\partial_x u_{\mathrm{DG}}^{n-1,k-1}(t)]^2 \partial_x\psi_i(x)\partial_x \psi_j(x) \: \dd x,$$ for $k=1,2,\ldots,$
where $u_{\mathrm{DG}}^{n-1,k-1}(t)$ is computed from the previous Picard iteration by using the stiffness matrix $[\tilde{K}_{k-1}^n]_{ij}.$

  \STATE \emph{Update}: $$[\tilde{K}_0^{n+1}]_{ij}=\int_0^1 [\partial_x u_{\mathrm{DG}}^n(t)]^2 \partial_x\psi_i(x)\partial_x \psi_j(x) \: \dd x,$$
where $\partial_x u_{\mathrm{DG}}^n(t)$ is computed using $[\tilde{K}_{k_{\mathrm{end}}}^{n}]$. Here $k_{{\mathrm{end}}}$ is either the maximal (final) Picard iteration number or the iteration at which a certain tolerance is achieved.

\STATE Now move to the next time interval $I_{n+1}$.
\end{algorithmic}
\end{algorithm}

Take $q=2$ as an example; we note that on each time interval $I_n$ the solution $\mathbf{U}(t)^{n}$ can be defined using the Legendre basis via 
$$\mathbf{U}(t)^{n}=\tilde{M}^{-\frac{1}{2}}\mathbf{Z}^n(t)=\tilde{M}^{-\frac{1}{2}}\left( \bm{\alpha}_1^n\phi^1(t)+ \bm{\alpha}_2^n\phi^2(t)+\bm{\alpha}_3^n\phi^3(t)\right).$$
Here $\bm{\alpha}_i^n$ for $i=1,2$ and $3$ are the coefficient vectors computed by extracting appropriate entries from $\mathbf{z}$.
This implies that we can formulate $\tilde{K}^{n+1}(t)$ by considering 
\begin{align*}
\tilde{K}^{n+1}(t) &= \int_0^1 \vert\nabla \mathbf{U}^n(t)\vert^2\partial_x \psi_i \partial_x \psi_j \:\mathrm{d} x\\
&= \int_0^1 \vert \tilde{M}^{-\frac{1}{2}}\mathrm{grad}(\bm{\alpha}_1^n)\phi^1(t)+ \tilde{M}^{-\frac{1}{2}}\mathrm{grad}(\bm{\alpha}_2^n)\phi^2(t)+\tilde{M}^{-\frac{1}{2}}\mathrm{grad}(\bm{\alpha}_3^n)\phi^3(t)\vert ^2\partial_x \psi_i \partial_x \psi_j \:\mathrm{d} x\\
&=  \int_0^1 \vert \tilde{M}^{-\frac{1}{2}}\mathrm{grad}(\bm{\alpha}_1^n)\phi^1(t)\vert^2  \partial_x\psi_i \partial_x \psi_j \:\mathrm{d} x +\int_0^1 \vert\tilde{M}^{-\frac{1}{2}}\mathrm{grad}(\bm{\alpha}_2^n)\phi^2(t)\vert^2  \partial_x \psi_i \partial_x \psi_j \:\mathrm{d} x \\
&\quad+\int_0^1 \vert\tilde{M}^{-\frac{1}{2}}\mathrm{grad}(\bm{\alpha}_3^n)\phi^3(t)\vert ^2  \partial_x \psi_i \partial_x \psi_j \:\mathrm{d} x\\
&\quad+ \int_0^1 2\left(\tilde{M}^{-\frac{1}{2}}\mathrm{grad}(\bm{\alpha}_1^n)\phi^1(t)\cdot\tilde{M}^{-\frac{1}{2}}\mathrm{grad}(\bm{\alpha}_2^n)\phi^2(t)\right)  \partial_x\psi_i \partial_x \psi_j \:\mathrm{d} x\\
&\quad+ \int_0^1 2\left(\tilde{M}^{-\frac{1}{2}}\mathrm{grad}(\bm{\alpha}_1^n)\phi^1(t)\cdot \tilde{M}^{-\frac{1}{2}}\mathrm{grad}(\bm{\alpha}_3^n)\phi^3(t)\right)  \partial_x\psi_i \partial_x \psi_j \:\mathrm{d} x\\
&\quad+ \int_0^1 2\left(\tilde{M}^{-\frac{1}{2}}\mathrm{grad}(\bm{\alpha}_2^n)\phi^2(t)\cdot \tilde{M}^{-\frac{1}{2}}\mathrm{grad}(\bm{\alpha}_3^n)\phi^3(t)\right) \partial_x \psi_i \partial_x \psi_j \:\mathrm{d} x\\
{}&:= K_{11}^{n+1}(t)+K_{22}^{n+1}(t)+K_{33}^{n+1}(t)+K_{12}^{n+1}(t)+K_{13}^{n+1}(t)+K_{23}^{n+1}(t).
\end{align*}
\begin{table}[ht]
  \begin{center}
   \caption{$||u(T)-u_{\mathrm{DG}}(t_N^{-})||_{L^2}+||\dot{u}(T)-\dot{u}_{\mathrm{DG}}(t_N^{-})||_{L^2} $   and corresponding convergence rates with respect to polynomial degrees $q=2, 3, 4$.}\label{tab:table1}
    \begin{tabular}{c c c c c } 
     \hline
      $ q$  & $h$  &  $k=h^2$   & $L^2$-error & rate \\
      \hline 
  $2$  & $2.50000\mathrm{e}-1$    & $6.25000\mathrm{e}-2$   & $1.2123\mathrm{e}-2$ & ---  \\
       & $2.00000\mathrm{e}-1$    &  $4.00000\mathrm{e}-2$   &$4.9774\mathrm{e}-3$ &  $1.9948$\\
      & $1.25000\mathrm{e}-1$    &  $ 1.56250\mathrm{e}-2$   &$1.1643\mathrm{e}-3$ &  $1.5455$\\
&    $6.25000\mathrm{e}-2$    &  $3.90625\mathrm{e}-3$      & $1.2454\mathrm{e}-4$ &  $1.6124$\\
  $3$  & $2.50000\mathrm{e}-1$    & $6.25000\mathrm{e}-2$   & $4.1533\mathrm{e}-4$ & ---  \\
       & $2.00000\mathrm{e}-1$    &  $4.00000\mathrm{e}-2$   &$1.8590\mathrm{e}-4$ &  $1.8012$\\
      & $1.25000\mathrm{e}-1$    &  $ 1.56250\mathrm{e}-2$   &$2.5283\mathrm{e}-5$ &  $2.1224$\\
&    $6.25000\mathrm{e}-2$    &  $3.90625\mathrm{e}-3$      &$1.5609\mathrm{e}-6$ &  $2.0089$\\
  $4$  & $2.50000\mathrm{e}-1$    & $6.25000\mathrm{e}-2$   &$2.5498\mathrm{e}-5$ & ---  \\
       & $2.00000\mathrm{e}-1$    &  $4.00000\mathrm{e}-2$   &$8.3999\mathrm{e}-6$ &  $2.4881$\\
      & $1.25000\mathrm{e}-1$    &  $ 1.56250\mathrm{e}-2$   &$9.7790\mathrm{e}-7$ &  $2.2878$\\
&    $6.25000\mathrm{e}-2$    &  $3.90625\mathrm{e}-3$      &$3.1115\mathrm{e}-8$ &  $2.4870$\\
\hline        
    \end{tabular}
    \end{center}
\end{table}
\begin{figure}[ht]
\centering
    \includegraphics[scale=0.7]{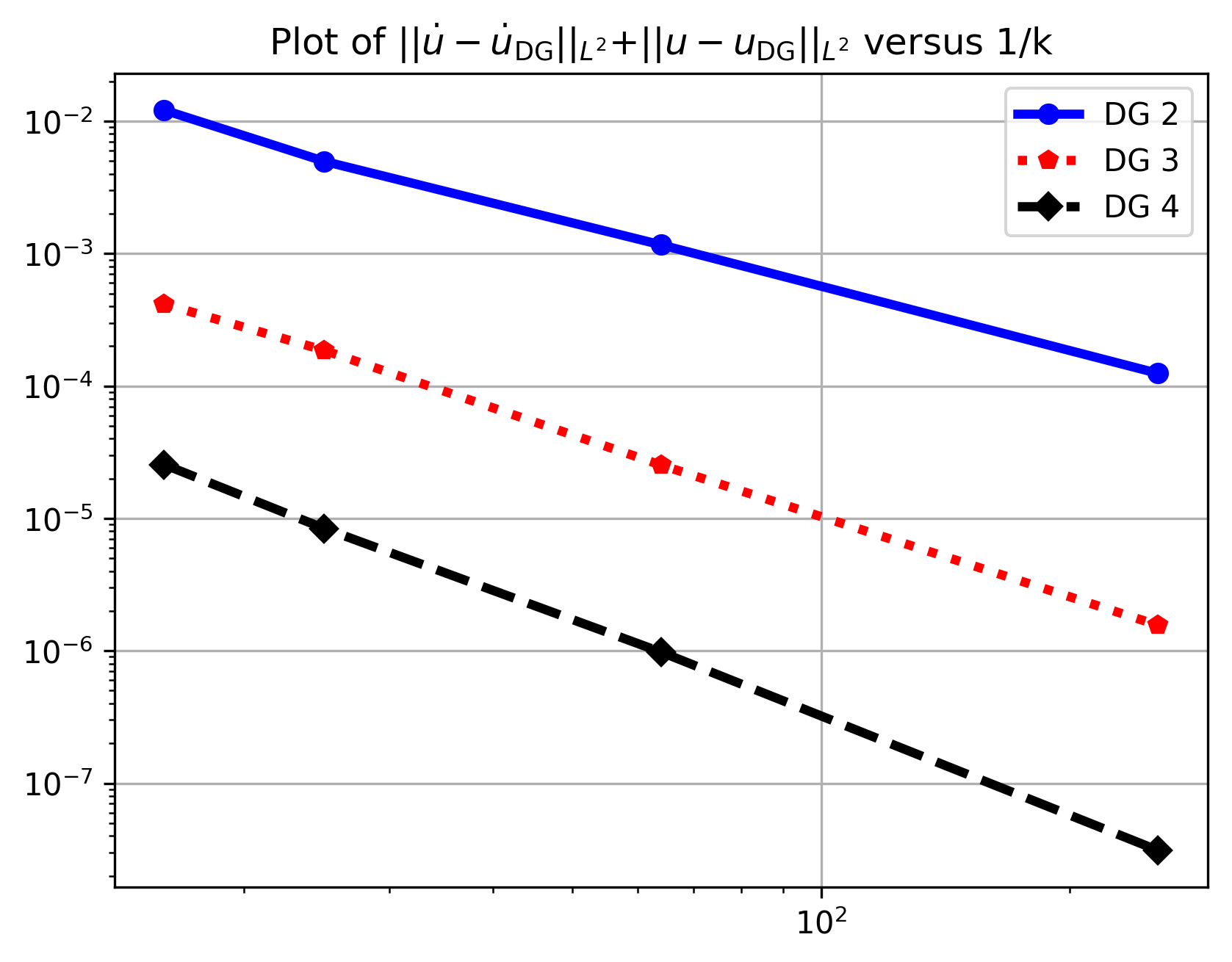}
\caption{Computed error  $\|u(T)-u_{\mathrm{DG}}(t_N^{-})\|_{L^2}+\|\dot{u}(T)-\dot{u}_{\mathrm{DG}}(t_N^{-})\|_{L^2} $ plotted against $1/k$ for polynomial degrees $q=2, 3, 4$.}\label{fig:plot6}
\end{figure}

We use $\mathrm{CG}$--$p$ elements where $p=q$ in space with $k=h^2$, $T=1$ and $\gamma=1$, and compute the errors $\|u(T)-u_{\mathrm{DG}}(t_N^{-})\|_{L^2}+\|\dot{u}(T)-\dot{u}_{\mathrm{DG}}(t_N^{-})\|_{L^2} $ versus $k$ for $k=h^2=6.25000\mathrm{e}-2, 4.00000\mathrm{e}-2, 1.56250\mathrm{e}-2$ and $3.90625\mathrm{e}-3$ with respect to polynomial degrees $2,3,4$ in Table \ref{tab:table1}. Note that here we use $h=2.50000\mathrm{e}-1$, $2.00000\mathrm{e}-2$, $1.25000\mathrm{e}-1$ and $6.25000\mathrm{e}-2$ instead of the conventional halving procedure; this is to avoid the accumulation of any unnecessary floating point errors resulting from a large number of time steps while still having sufficient data to compute the convergence rates. The computed errors are shown in Figure \ref{fig:plot6} in a log-log scale. As expected, the error decreases as we increase the polynomial degree $q$ or decrease the time step $k$. By Remark \ref{uniform_grid}, we expect convergence rates of order $1.5, 2.0$ and $2.5$ for $q=p=2,3$ and $4$ respectively, which are consistent with the numerical results shown in Table \ref{tab:table1}.

\section{Appendices}
\subsection{Proof of the auxiliary lemma}\label{lemm_pf}
\begin{lmm}
Under the assumptions stated in Theorem \ref{convergence_thm}, there exists a constant $C_{\tau}>0$ such that, for $t\in I_n$, $n=1,2\ldots N$,
\begin{equation}
\vert \mathsf{A}_{t}(\bm{\phi}(t)\mathrm{e}^{\gamma t}; \bm{\varphi}(t),\bm{\psi}(t))\vert\leq C_{\tau}\|\nabla \bm{\varphi}(t)\|_{L^2}\|\nabla\bm{\psi}(t)\|_{L^2}.
\end{equation}
\end{lmm}
\begin{proof}
Note that for $t\in I_n, n=1,2,\ldots,N$,
\begin{flalign*}
&\partial_t [A_{i\alpha j \beta}(\nabla\mathbf{W}(t) \mathrm{e}^{\gamma t}+\tau (\nabla\bm{\phi}(t)-\nabla\mathbf{W}(t)) \mathrm{e}^{\gamma t})]&&\\
&=\sum_{k,m=1}^d \frac{\partial A_{i\alpha j\beta}}{\partial \eta_{km}}\left(\nabla \mathbf{W}(t) e^{\gamma t}+\tau \nabla(\bm{\phi}(t)-\mathbf{W}(t))\mathrm{e}^{\gamma t}\right)  \partial_{m}\left(\partial_t (\mathbf{W}_k(t)\mathrm{e}^{\gamma t})+\tau \partial_t((\bm{\phi}_{k}(t)-\mathbf{W}_k(t))\mathrm{e}^{\gamma t})\right).
\end{flalign*}
Since the values of the function $\nabla\mathbf{W}(t)\mathrm{e}^{\gamma t}+\tau\nabla(\bm{\phi}(t)-\mathbf{W}(t))\mathrm{e}^{\gamma t}$ for $t\in [0,T]$, $\tau\in (0,1)$, belong to the compact convex subset $\mathcal{M}_{\delta}$ of $\mathbb{R}^{d\times d}$, and because $A_{i\alpha j \beta}$ is sufficiently smooth (and in particular continuously differentiable on $\mathcal{M}_{\delta}$), we have 
\begin{flalign*}
&\bigg\vert\sum_{i,\alpha,j,\beta=1}^d \frac{1}{2}\int_0^1 \left( \partial_t [A_{i\alpha j \beta}(\nabla \mathbf{W}(t) \mathrm{e}^{\gamma t}+\tau\nabla(\bm{\phi}(t)-\mathbf{W}(t))\mathrm{e}^{\gamma t})]\partial_j \bm{\varphi}_{\beta}(t),\partial_{i}\bm{\psi}_{\alpha}(t)\right)_{L^2} \dd \tau \bigg\vert &&\\
&\leq  c\left(\|\nabla \mathbf{\dot W}(t)\|_{L^{\infty}}+\|\nabla\mathbf{W}(t)\|_{L^{\infty}}\right)\|\nabla \bm{\varphi}(t)\|_{L^2}\|\nabla \bm{\psi}(t)\|_{L^2} &&\\
& \quad+c\left(\|\nabla(\bm{\dot \phi}(t) - \mathbf{\dot W}(t))\|_{L^{\infty}}+\|\nabla(\bm{\phi}(t) -\mathbf{W}(t))\|_{L^{\infty}}\right)\|\nabla \bm{\varphi}(t)\|_{L^2}\|\nabla \bm{\psi}(t)\|_{L^2}&&\\
&\leq \tilde{C}h^{-\frac{d}{2}}\left(\|\nabla (\bm{\dot \phi}(t) -\mathbf{\dot W}(t))\|_{L^2}+\|\nabla(\bm{\phi}(t) -\mathbf{W}(t))\|_{L^2}  \right)\|\nabla \bm{\varphi}(t)\|_{L^2}\|\nabla \bm{\psi}(t)\|_{L^2} &&\\
&\quad+ \tilde{C}\|\nabla \bm{\varphi}(t)\|_{L^2}\|\nabla \bm{\psi}(t)\|_{L^2},&&
\end{flalign*}
where we have applied the inverse inequality (ii,b) and property (iii,c) of the nonlinear projection $\mathbf{W}$. We shall bound $\|\nabla(\bm{\dot \phi}(t) -\mathbf{\dot W}(t))\|_{L^2}$ and $\|\nabla(\bm{\phi}(t) -\mathbf{W}(t))\|_{L^2}$ for $t\in I_n, n=1,2\ldots, N$. Applying the triangle inequality, we have
\begin{align*}
\|\nabla (\bm{\phi}(t) -\mathbf{W}(t))\|_{L^2}\leq \|\nabla (\bm{\phi}(t) -\Pi_k\mathbf{W}(t))\|_{L^2}+\|\nabla(\Pi_k\mathbf{W}(t)-\mathbf{W}(t))\|_{L^2}.
\end{align*}
Note that for $t\in I_n, n=1,2\ldots N$,
\begin{flalign*}
\|\nabla(\bm{\phi}(t)-\Pi_k\mathbf{W}(t))\|_{L^2} \leq & \:\|\nabla(\bm{\phi}(t_n^{-})-\Pi_k\mathbf{W}(t_n^{-}))\|_{L^2}+\int_{t}^{t_n}\|\partial_s( \nabla\bm{\phi}(s)-\nabla\Pi_k\mathbf{W}(s))\|_{L^2}\:\dd s &&\\
\leq& \:\|\nabla(\bm{\phi}(t_n^{-})-\Pi_k\mathbf{W}(t_n^{-}))\|_{L^2}+\int_{t_{n-1}}^{t_n}\|\partial_t( \nabla\bm{\phi}(t)-\nabla\Pi_k\mathbf{W}(t))\|_{L^2}\,\dd t &&\\
\leq &  \:\|\nabla(\bm{\phi}(t_n^{-})-\Pi_k\mathbf{W}(t_n^{-}))\|_{L^2}+C_0h^{-1}\sqrt{k_n}\left(\int_{t_{n-1}}^{t_n}\|\partial_t( \bm{\phi}(t)-\Pi_k\mathbf{W}(t))\|_{L^2}^2 \,\dd t\right)^{\frac{1}{2}} &&\\
\leq&  \:C(\mathbf{v})\left(\sum_{i=1}^n k_i h^{2r+2}+ \frac{k_i^{2q_i+1}}{q_i^{2(s-1)}} \right)^{\frac{1}{2}},&&
\end{flalign*}
where we have used the inverse inequality (ii,a), H\"{o}lder's inequality, the fact that $\bm{\phi}\in \mathcal{F}$ and the assumption that $\mu_i k_i
\leq h^2$ for each $i=1,\ldots,N$. Here $C(\mathbf{v})$ denotes a constant depending on $\mathbf{v}$, which may vary throughout this proof. On the other hand,
\begin{flalign*}
\|\nabla(\Pi_k\mathbf{W}(t)-\mathbf{W}(t))\|_{L^2} &\leq  \|\nabla(\Pi_k\mathbf{W}(t_n^{-})-\mathbf{W}(t_n^{-}))\|_{L^2}+\int_{t}^{t_n} \|\partial_s(\nabla\Pi_k\mathbf{W}(s)-\nabla\mathbf{W}(s))\|_{L^2} \dd s\\
&\leq \int_{t_{n-1}}^{t_n} \|\partial_t(\nabla\Pi_k\mathbf{W}(t)-\nabla\mathbf{W}(t))\|_{L^2}\dd t \quad\mbox{ ( since } \Pi_k\mathbf{W}(t_n^{-})=\mathbf{W}(t_n^{-})\: ) &&\\
&\leq \sqrt{k_n}\bigg( \int_{t_{n-1}}^{t_n}\|\partial_t(\nabla\Pi_k\mathbf{W}(t)-\nabla\mathbf{W}(t))\|_{L^2}^2 \dd t \bigg)^{\frac{1}{2}} &&\\
&\leq   C \frac{k_n^{q_n+\frac{1}{2}}}{q_n^{s-1}}\|\mathbf{W}\|_{H^{s}(I_n; H_0^1)},\quad  \text{ for } \mathbf{W}\in H^{s}([0,T]; H_0^1), &&
\end{flalign*}
where we have used inequality (\ref{prop6}) with the $L^2$ norm in space replaced by the $H^1$ semi-norm. Thus,
\begin{equation}\label{L2bound}
\max_{t\in I_n, 1\leq n\leq N}\|\nabla(\bm{\phi}(t)-\mathbf{W}(t))\|_{L^2}\leq C(\mathbf{v})\left( \sum_{i=1}^n k_i h^{2r+2}+ \frac{k_i^{2q_i+1}}{q_i^{2(s-1)}} \right)^{\frac{1}{2}}+ C \frac{k_n^{q_n+\frac{1}{2}}}{q_n^{s-1}}\|\mathbf{W}\|_{H^{s}(I_n; H_0^1)}.
\end{equation}
Applying the triangle inequality to the time derivative term, we have
\begin{align*}
\|\nabla(\bm{\dot\phi}(t) -\mathbf{\dot W}(t))\|_{L^2}\leq \|\partial_t(\nabla\bm{\phi}(t) -\nabla\Pi_k\mathbf{W}(t))\|_{L^2}+\|\partial_t(\nabla\Pi_k\mathbf{W}(t)-\nabla\mathbf{W}(t))\|_{L^2}.
\end{align*}
Note that for $t\in I_n, n=1,2,\ldots N$,
\begin{align*}
\|\partial_{t}(\nabla\bm{\phi}(t)-\nabla\Pi_k\mathbf{W}(t))\|_{L^2}
&\leq  C_0h^{-1}\|\partial_{t}(\bm{\phi}(t)-\Pi_k\mathbf{W}(t))\|_{L^2}\\
&\leq C_0h^{-1}\|\partial_t (\bm{\phi}(t_n^{-})-\Pi_k\mathbf{W}(t_n^{-}))\|_{L^2}+C_0h^{-1}\int_{t}^{t_n}\|\partial_{ss}( \bm{\phi}(s)-\Pi_k\mathbf{W}(s))\|_{L^2} \dd s \\
&\leq  C_0h^{-1}\|\partial_t(\bm{\phi}(t_n^{-})-\Pi_k\mathbf{W}(t_n^{-}))\|_{L^2}+C_0h^{-1}\int_{t_{n-1}}^{t_n}\|\partial_{tt}( \bm{\phi}(t)-\Pi_k\mathbf{W}(t))\|_{L^2} \dd t.
\end{align*}
Since $\bm{\phi}\in\mathcal{F}$, we have $$\|\partial_t(\bm{\phi}(t_n^{-})-\Pi_k\mathbf{W}(t_n^{-}))\|_{L^2}\leq C_{\ast}(\mathbf{v})\left( \sum_{i=1}^n k_i h^{2r+2}+\frac{k_i^{2q_i+1}}{q_i^{2(s-1)}} \right)^{\frac{1}{2}}.$$ Using the inverse inequality in time (\ref{inveqtime}), we obtain
\begin{flalign*}
\int_{t_{n-1}}^{t_n}\|\partial_{tt}( \bm{\phi}(t)-\Pi_k\mathbf{W}(t))\|_{L^2} \dd t\leq &\: \sqrt{k_n}\left(\int_{t_{n-1}}^{t_n}\|\partial_{tt}( \bm{\phi}(t)-\Pi_k\mathbf{W}(t))\|_{L^2}^2 \dd t\right)^{\frac{1}{2}}&&\\
\leq & \: C_2\frac{1}{\sqrt{k_n}}\left(\int_{t_{n-1}}^{t_n}\|\partial_{t}( \bm{\phi}(t)-\Pi_k\mathbf{W}(t))\|_{L^2}^2 \dd t\right)^{\frac{1}{2}}&&\\
\leq &\: C(\mathbf{v})\frac{1}{\sqrt{k_n}}\left(  \sum_{i=1}^n k_i h^{2r+2}+ \frac{k_i^{2q_i+1}}{q_i^{2(s-1)}}\right)^{\frac{1}{2}}.&&
\end{flalign*}
On the other hand, for $t\in I_n, n=1,2,\ldots, N$, we have
\begin{align*}
&\|\partial_t (\nabla\Pi_k \mathbf{W}(t)-\nabla\mathbf{W}(t))\|_{L^2}\\
&\leq \|\partial_t (\nabla\Pi_k \mathbf{W}(t_n^{-})-\nabla\mathbf{W}(t_n^{-}))\|_{L^2} +\int_{t}^{t_n}\|\partial_{ss}(\nabla\Pi_k\mathbf{W}(s)-\nabla\mathbf{W}(s))\|_{L^2}\:\dd s &&\\
&\leq \int_{t_{n-1}}^{t_n}\|\partial_{tt}(\nabla\Pi_k\mathbf{W}(t)-\nabla\mathbf{W}(t))\|_{L^2}\,\dd t \quad (\mbox{ since } \Pi_k\mathbf{W}(t_n^{-})=\mathbf{W}(t_n^{-}) ) \\
&\leq \sqrt{k_n}\left(\int_{t_{n-1}}^{t_n}\|\partial_{tt}(\nabla\Pi_k\mathbf{W}(t)-\nabla\mathbf{W}(t))\|_{L^2}^2 \,\dd t\right)^{\frac{1}{2}}\\
&\leq C \frac{k_n^{q_n-\frac{1}{2}}}{q_n^{s-3}}\|\mathbf{W}\|_{H^{s}(I_n; H_0^1)},\quad  \text{ for } \mathbf{W}\in H^{s}([0,T]; H_0^1).
\end{align*}
Thus 
\begin{align}\label{derivativebound}
\max_{t\in I_n, 1\leq n\leq N}\|\nabla\bm{\dot\phi}(t)-\nabla\mathbf{\dot W}(t)\|_{L^2} &\leq \:C(\mathbf{v})(h^{-1}+ h^{-1}k_n^{-\frac{1}{2}})\left( \sum_{i=1}^n k_i h^{2r+2}+ \frac{k_i^{2q_i+1}}{q_i^{2(s-1)}}\right)^{\frac{1}{2}} \nonumber\\
&\quad+C\frac{k_n^{q_n-\frac{1}{2}}}{q_n^{s-3}}\|\mathbf{W}\|_{H^{s}([0,T]; H_0^1)}.
\end{align}
Combining (\ref{L2bound}) and (\ref{derivativebound}), we have 
\begin{align}\label{constbound}
&\tilde{C}h^{-\frac{d}{2}}\max_{t\in I_n, 1\leq n\leq N} \left(\|\nabla\bm{\phi}(t)-\nabla\mathbf{ W}(t)\|_{L^2}+\|\nabla\bm{\dot\phi}(t)-\nabla\mathbf{\dot W}(t)\|_{L^2}\right) \nonumber\\
&\leq  \tilde{C}h^{-\frac{d}{2}} C(\mathbf{v})(1+h^{-1}+h^{-1}k_n^{-\frac{1}{2}})\left(  \sum_{i=1}^n k_i h^{2r+2}+ \frac{k_i^{2q_i+1}}{q_i^{2(s-1)}} \right)^{\frac{1}{2}} \nonumber\\
&\quad+\tilde{C}h^{-\frac{d}{2}}C\bigg(\frac{k_n^{q_n+\frac{1}{2}}}{q_n^{s-1}}+\frac{k_n^{q_n-\frac{1}{2}}}{q_n^{s-3}}\bigg)\|\mathbf{W}\|_{H^{s}([0,T]; H_0^1)}.
\end{align}
Since $r> \frac{d}{2}+1$, $k_i^{q_i-\frac{1}{2}}=o(h^{1+\frac{d}{2}})$ and $\mu k_i\leq h^2\leq \nu_i k_i$ for each $i=1,2,\ldots N,$ we can choose $h_0>0$ such that for $h\leq h_0$, the right-hand side of (\ref{constbound}) is bounded by $1$. Thus, (\ref{lemmaie}) follows by taking $C_{\tau}=\tilde{C}+1$. The constant $C_{\tau}$ defined in this way does not depend on $C_{\ast}(\mathbf{v})$.
\end{proof}
\subsection{Approximation properties of the elliptic projection}\label{W}
Here we derive the properties (iii,a)--(iii,c) of the nonlinear projection $\mathbf{W}$. We write $a\lesssim b$ if there exists a universal constant $C>0$ independent of the spatial discretization parameter $h$ such that $a\leq Cb.$
\subsubsection{\mathinhead{L^2}{L2} bound on \mathinhead{\nabla(\mathbf{v}-\mathbf{W})}{nabla(v-W)} and \mathinhead{L^{\infty}}{Linfty} bound on \mathinhead{\nabla\mathbf{W}}{nablaW}}
Recall that for each $t\in [0,T]$, $a(\mathbf{W}(t),\bm{\varphi})=a(\mathbf{v}(t),\bm{\varphi})$  for all $\bm{\varphi}\in \mathcal{V}_h.$ Let $\mathcal{P}_h\colon L^2\to \mathcal{V}_h$ denote the standard $L^2$-projection operator in the spatial direction. Then we have 
\begin{equation}
a(\mathbf{W}(t),\bm{\varphi})-a(\mathcal{P}_h\mathbf{v}(t),\bm{\varphi})= a(\mathbf{v}(t),\bm{\varphi})-a(\mathcal{P}_h\mathbf{v}(t),\bm{\varphi}) \quad \mbox{ for all } \bm{\varphi}\in \mathcal{V}_h.
\end{equation}
That is, 
\begin{align*}
{}&\int_0^1 \tilde{a}(\mathcal{P}_h\mathbf{v}(t)\mathrm{e}^{\gamma t}+\tau(\mathbf{W}(t)-\mathcal{P}_h\mathbf{v}(t))\mathrm{e}^{\gamma t}; \mathbf{W}(t)-\mathcal{P}_h\mathbf{v}(t),\bm{\varphi}) \:\dd \tau\\
&= \int_0^1 \tilde{a}(\mathcal{P}_h\mathbf{v}(t)\mathrm{e}^{\gamma t}+\tau(\mathbf{v}(t)-\mathcal{P}_h\mathbf{v}(t))\mathrm{e}^{\gamma t}; \mathbf{v}(t)-\mathcal{P}_h\mathbf{v}(t),\bm{\varphi}) \:\dd \tau\
\end{align*}
Define the following subset of $H_0^1$,
$$\mathcal{F}=\{\bm{\phi}\in \mathcal{V}_h \colon \|\nabla (\bm{\phi}-\mathcal{P}_h\mathbf{v})\|_{L^2}\leq C_{\ast} h^r \|\mathbf{v}\|_{H^{r+1}} \mbox{ for }\frac{d}{2}< r \leq \min(p,m-1)\}$$
where $C_{\ast}$ is a constant independent of $h$. The set $\mathcal{F}$ is non-empty since for each fixed $t\in [0,T]$, $\mathcal{P}_h\mathbf{v}(t)\in \mathcal{F}$. Furthermore, $\mathcal{F}$ is a closed and convex subset of $H_0^1$. We define the fixed point mapping $\mathcal{N}$ on $\mathcal{F}$ as follows. Given $\bm{\phi}\in \mathcal{F}$, we denote by $\mathbf{W}_{\phi}\in\mathcal{V}_h$, the solution to the following linear variational problem: find $\mathbf{W}_{\phi}\in\mathcal{V}_h$ such that 
\begin{align*}
{}&\int_0^1 \tilde{a}(\mathcal{P}_h\mathbf{v}(t)\mathrm{e}^{\gamma t}+\tau(\bm{\phi}(t)-\mathcal{P}_h\mathbf{v}(t))\mathrm{e}^{\gamma t}; \mathbf{W}_{\phi}(t)-\mathcal{P}_h\mathbf{v}(t),\bm{\varphi}) \:\dd \tau\\
&= \int_0^1 \tilde{a}(\mathcal{P}_h\mathbf{v}(t)\mathrm{e}^{\gamma t}+\tau(\mathbf{v}(t)-\mathcal{P}_h\mathbf{v}(t))\mathrm{e}^{\gamma t}; \mathbf{v}(t)-\mathcal{P}_h\mathbf{v}(t),\bm{\varphi}) \:\dd \tau\ \quad \mbox{ for all } \bm{\varphi}\in\mathcal{V}_h.
\end{align*}
Since $\mathcal{V}_h$ is a finite dimensional linear space, the existence and uniqueness of $\mathbf{W}_{\phi}(t)\in\mathcal{V}_h$ for each $t\in [0,T]$ follows if we can show that $\int_0^1 \tilde{a}(\mathcal{P}_h\mathbf{v}(t)\mathrm{e}^{\gamma t}+\tau(\bm{\phi}(t)-\mathcal{P}_h\mathbf{v}(t)); \cdot, \cdot)$ is coercive on $\mathcal{V}_h\times \mathcal{V}_h$ in the $|\cdot |_{H^1}$ semi-norm. This is indeed true in view of the assumption (S2b).
For each $t\in [0,T]$, if we take $\mathbf{W}(t)=\mathbf{W}_{\phi}(t)$, we have 
\begin{equation}
\|\nabla(\mathbf{W}(t)-\mathcal{P}_h\mathbf{v}(t))\|_{L^2}\leq C_{\ast} h^r\| \mathbf{v}(t)\|_{H^{r+1}}, \quad\frac{d}{2}<r\leq \min(p,m-1).
\end{equation}
By the approximation properties of $\mathcal{P}_h$ in the $|\cdot|_{H^1}$ semi-norm, we have 
\begin{equation}\label{approx_prop}
\|\nabla( \mathbf{v}(t)-\mathcal{P}_h\mathbf{v}(t))\|_{L^2}\lesssim h^r\|\mathbf{v}(t)\|_{H^{r+1}}, \quad\frac{d}{2}<r\leq \min(p,m-1).
\end{equation}
It follows from the triangle inequality that 
\begin{equation}\label{grad_error_bound}
\| \nabla (\mathbf{W}(t)-\mathbf{v}(t))\|_{L^2}\lesssim h^r\|\mathbf{v}(t)\|_{H^{r+1}}, \quad\frac{d}{2}<r\leq \min(p,m-1).
\end{equation}
By the approximation properties of $\mathcal{P}_h$ in the $|\cdot|_{W^{1,\infty}}$ semi-norm, we have 
\begin{equation}\label{approx_prop_inf}
\|\nabla( \mathbf{v}(t)-\mathcal{P}_h\mathbf{v}(t))\|_{L^{\infty}}\lesssim h^{r-\frac{d}{2}}\|\mathbf{v}(t)\|_{H^{r+1}}, \quad\frac{d}{2}<r\leq \min(p,m-1).
\end{equation}
Combining (\ref{grad_error_bound}) and (\ref{approx_prop_inf}), we obtain
\begin{align*}
\|\nabla \mathbf{W}(t)\|_{L^{\infty}}&\leq \|\nabla \mathbf{v}(t)\|_{L^{\infty}}+\|\nabla (\mathbf{W}(t)-\mathbf{v}(t))\|_{L^{\infty}} \\
&\leq \|\nabla \mathbf{v}(t)\|_{L^{\infty}}+\|\nabla (\mathbf{W}(t)-\mathcal{P}_h\mathbf{v}(t))\|_{L^{\infty}}+\|\nabla (\mathcal{P}_h\mathbf{v}(t)-\mathbf{v}(t))\|_{L^{\infty}} \\
&\leq \|\nabla\mathbf{v}(t)\|_{L^{\infty}}+C_1 h^{-\frac{d}{2}}\| \nabla(\mathbf{W}(t)-\mathcal{P}_h\mathbf{v}(t))\|_{L^2}+ C(\mathbf{v})h^{r-\frac{d}{2}} \\
{}&\leq c_0,
\end{align*}
for some constant $c_0$. The last inequality follows from the boundedness of $\nabla\mathbf{v}$ and the fact that $r>\frac{d}{2}$, while the second last line follows from (ii,b) and (\ref{approx_prop_inf}).
\subsubsection{\mathinhead{L^2}{L2} bound on \mathinhead{\nabla(\mathbf{\dot{W}}-\mathbf{\dot{v}})}{nabladotWdotv} and \mathinhead{L^{\infty}}{Linfty} bound on \mathinhead{\nabla\mathbf{\dot{W}}}{nabladotW}}
For the estimate of the $L^2$ bound on $\nabla(\mathbf{\dot{W}}-\mathbf{\dot{v}})$, we follow the proof from Section 6 in \cite{OS2}. We need to show that $t\mapsto \mathbf{W}(t)$ is differentiable with respect to $t.$ For $\mathbf{U}\in\mathcal{V}_h$ and $t\in [0,T]$, we notice that the mapping $\bm{\varphi}\mapsto a(\mathbf{U},\bm{\varphi})-a(\mathbf{v}(t),\bm{\varphi})$ is a bounded linear functional on $\mathcal{V}_h$; hence by \emph{Riesz representation theorem}, there exists a unique $\mathcal{A}(t,\mathbf{U})\in \mathcal{V}_h$ such that 
$$
(\mathcal{A}(t,\mathbf{U}),\bm{\varphi})=a(\mathbf{U}, \bm{\varphi})-a(\mathbf{v}(t),\bm{\varphi}).
$$
It follows from the linearization process that the derivative of the nonlinear mapping $(t,\mathbf{U})\mapsto \mathcal{A}(t, \mathbf{U})$ with respect to $\mathbf{U}$, evaluated at $\mathbf{U}=\mathbf{W}(t)$, exists and is invertible for any $t\in [0,T]$. We also have $\mathcal{A}(t,\mathbf{W}(t))=0$. Since $\mathbf{v}(t)$ is differentiable with respect to $t$, it follows that $\mathcal{A}(t,\mathbf{U})$ is differentiable in a neighbourhood of $(t_0, \mathbf{W}(t_0))$ for any $t_0\in (0,T)$. We then deduce from the \emph{implicit function theorem} that $t\mapsto \mathbf{W}(t)$ is differentiable in $(0,T)$.
Next, we derive the error bound of $\|\nabla (\mathbf{\dot W}(t)-\mathbf{\dot v}(t)) \|_{L^2}.$ By definition of $\mathbf{W}(t)$, we have 
\begin{align*}
&\int_0^1 \tilde{a}(\mathcal{P}_h\mathbf{v}(t)\mathrm{e}^{\gamma t}+\tau (\mathbf{W}(t)-\mathcal{P}_h\mathbf{v}(t))\mathrm{e}^{\gamma t};\mathbf{W}(t)-\mathcal{P}_h\mathbf{v}(t), \bm{\varphi}) \:\dd \tau \\
&= \int_0^1 \tilde{a}(\mathcal{P}_h\mathbf{v}(t)\mathrm{e}^{\gamma t}+\tau (\mathbf{v}(t)-\mathcal{P}_h\mathbf{v}(t))\mathrm{e}^{\gamma t};\mathbf{ v}(t)-\mathcal{P}_h\mathbf{ v}(t), \bm{\varphi}) \:\dd \tau \quad \mbox{ for all }\bm{\varphi}\in\mathcal{V}_h. 
\end{align*}
After differentiation with respect to $t$, we have 
\begin{align*}
&\int_0^1 \tilde{a}(\mathcal{P}_h\mathbf{v}(t)\mathrm{e}^{\gamma t}+\tau (\mathbf{W}(t)-\mathcal{P}_h\mathbf{v}(t))\mathrm{e}^{\gamma t};\mathbf{\dot W}(t)-\mathcal{P}_h\mathbf{\dot v}(t), \bm{\varphi}) \:\dd \tau \\
{}&\quad+\int_0^1\int_{\Omega} \sum_{i,\alpha,j,\beta, k,m=1}^d \bigg(\bigg\{\frac{\partial A_{i\alpha j \beta}}{\partial \eta _{km}}(\nabla \mathcal{P}_h\mathbf{v}(t)\mathrm{e}^{\gamma t}+\tau(
\nabla\mathbf{W}(t)-\nabla\mathcal{P}_h\mathbf{v}(t))\mathrm{e}^{\gamma t})\\
{}&\quad\times\partial_m\left( \partial_t[\mathcal{P}_h\mathbf{ v}_k(t)\mathrm{e}^{\gamma t}]+\tau \partial_t[(\mathbf{ W}_k(t)-\mathcal{P}_h\mathbf{ v}_k(t))\mathrm{e}^{\gamma t}]\right)\bigg\} \partial_j (\mathbf{W}-\mathcal{P}_h\mathbf{v})_{\beta},\partial_j\bm{\varphi}_{\alpha} \bigg)_{L^2}\:\dd x \:\dd \tau\\
&= \int_0^1 \tilde{a}(\mathcal{P}_h\mathbf{v}(t)\mathrm{e}^{\gamma t}+\tau (\mathbf{ v}(t)-\mathcal{P}_h\mathbf{ v}(t))\mathrm{e}^{\gamma t};\mathbf{\dot v}(t)-\mathcal{P}_h\mathbf{\dot v}(t), \bm{\varphi}) \:\dd \tau\\
{}&\quad+\int_0^1\int_{\Omega}\sum_{i,\alpha,j,\beta, k,m=1}^d\bigg(\bigg\{ \frac{\partial A_{i\alpha j \beta}}{\partial \eta _{k,m}}(\nabla \mathcal{P}_h\mathbf{v}(t)\mathrm{e}^{\gamma t}+\tau(\nabla\mathbf{v}(t)-\nabla\mathcal{P}_h\mathbf{v}(t))\mathrm{e}^{\gamma t})\\
{}&\quad\times\partial_m\left( \partial_t [\mathcal{P}_h\mathbf{v}_k(t)\mathrm{e}^{\gamma t}]+\tau \partial_t [(\mathbf{v}_k(t)-\mathcal{P}_h\mathbf{v}_k(t))\mathrm{e}^{\gamma t}]\right)\bigg\} \partial_j (\mathbf{ v}-\mathcal{P}_h\mathbf{ v})_{\beta},\partial_j\bm{\varphi}_{\alpha} \bigg)_{L^2}\:\dd x \:\dd \tau,
\end{align*}
for all $\bm{\varphi}\in\mathcal{V}_h.$
Rearranging gives
\begin{flalign*}
&\int_0^1 \tilde{a}(\mathcal{P}_h\mathbf{v}(t)\mathrm{e}^{\gamma t}+\tau (\mathbf{W}(t)-\mathcal{P}_h\mathbf{v}(t))\mathrm{e}^{\gamma t};\mathbf{\dot W}(t)-\mathcal{P}_h\mathbf{\dot v}(t), \bm{\varphi})\: \dd \tau &&\\
&= \int_0^1 \tilde{a}(\mathcal{P}_h\mathbf{v}(t)\mathrm{e}^{\gamma t}+\tau (\mathbf{v}(t)-\mathcal{P}_h\mathbf{v}(t))\mathrm{e}^{\gamma t};\mathbf{\dot v}(t)-\mathcal{P}_h\mathbf{\dot v}(t), \bm{\varphi}) \:\dd \tau&&\\
&\quad+\int_0^1\int_{\Omega}  \sum_{i,\alpha,j,\beta, k,m=1}^d \bigg(\bigg\{ \frac{\partial A_{i\alpha j \beta}}{\partial \eta _{k,m}}(\nabla \mathcal{P}_h\mathbf{v}(t)\mathrm{e}^{\gamma t}+\tau(
\nabla\mathbf{v}(t)-\nabla\mathcal{P}_h\mathbf{v}(t))\mathrm{e}^{\gamma t})&&\\
&\quad\times\partial_m\left( \partial_t [\mathcal{P}_h\mathbf{v}_k(t)\mathrm{e}^{\gamma t}]+\tau\partial_t [(\mathbf{v}_k(t)-\mathcal{P}_h\mathbf{v}_k(t))\mathrm{e}^{\gamma t}]\right)\bigg\} \partial_j (\mathbf{v}-\mathcal{P}_h\mathbf{ v})_{\beta},\partial_j\bm{\varphi}_{\alpha} \bigg)_{L^2} \:\dd x\: \dd \tau &&\\
&\quad-\int_0^1\int_{\Omega} \sum_{i,\alpha,j,\beta, k,m=1}^d \bigg(\bigg\{\frac{\partial A_{i\alpha j \beta}}{\partial \eta _{k,m}}(\nabla \mathcal{P}_h\mathbf{v}(t)\mathrm{e}^{\gamma t}+\tau(
\nabla\mathbf{W}(t)-\nabla\mathcal{P}_h\mathbf{v}(t))\mathrm{e}^{\gamma t})&&\\
&\quad\times\partial_m\left( \partial_t [\mathcal{P}_h\mathbf{ v}_k(t)\mathrm{e}^{\gamma t}]+\tau\partial_t [(\mathbf{W}_k(t)-\mathcal{P}_h\mathbf{ v}_k(t))\mathrm{e}^{\gamma t}]\right)\bigg\}\partial_j (\mathbf{ W}-\mathcal{P}_h\mathbf{ v})_{\beta},\partial_j\bm{\varphi}_{\alpha} \bigg)_{L^2} \:\dd x \: \dd \tau &&\\
&:=  \:T_1+T_2+T_3.
\end{flalign*}
Taking $\bm{\varphi}(t)=\mathbf{\dot W}(t)-\mathcal{P}_h\mathbf{\dot v}(t)$, we have
\begin{flalign*}
T_1\lesssim & \:h^r \|\mathbf{\dot v}(t)\|_{H^{r+1}}\|\nabla(\mathbf{\dot W}(t)-\mathcal{P}_h\mathbf{\dot v}(t))\|_{L^2}&&\\
T_2\lesssim &\: \bigg( \|\nabla\mathcal{P}_h\mathbf{\dot v}(t)\|_{L^{\infty}}+\|\nabla  (\mathbf{\dot v}(t)-\mathcal{P}_h\mathbf{\dot v}(t))\|_{L^{\infty}}+\|\nabla\mathcal{P}_h\mathbf{v}(t)\|_{L^{\infty}}+\|\nabla ( \mathbf{v}(t)-\mathcal{P}_h\mathbf{v}(t))\|_{L^{\infty}}\bigg)&&\\
{}& \times \|\nabla \mathbf{v}(t)-\nabla\mathcal{P}_h\mathbf{v}(t)\|_{L^2}\|\nabla ( \mathbf{\dot W}(t)-\mathcal{P}_h\mathbf{\dot v}(t))\|_{L^2}&&\\
{}\lesssim& \: h^r \|\mathbf{v}(t)\|_{H^{r+1}} \|\nabla ( \mathbf{\dot W}(t)-\mathcal{P}_h\mathbf{\dot v}(t))\|_{L^2}.&&\\
T_3 \lesssim & \: \bigg( \|\nabla\mathcal{P}_h\mathbf{\dot v}(t)\|_{L^{\infty}}+\|\nabla ( \mathbf{\dot W}(t)-\mathcal{P}_h\mathbf{\dot v}(t))\|_{L^{\infty}}+\|\nabla\mathcal{P}_h\mathbf{v}(t)\|_{L^{\infty}}+\|\nabla ( \mathbf{W}(t)-\mathcal{P}_h\mathbf{v}(t))\|_{L^{\infty}}\bigg) &&\\
{}& \times \|\nabla (\mathbf{W}(t)-\mathcal{P}_h\mathbf{v}(t))\|_{L^2}\|\nabla ( \mathbf{\dot W}(t)-\mathcal{P}_h\mathbf{\dot v}(t))\|_{L^2}\\
{}\lesssim &\: h^r \|\mathbf{v}(t)\|_{H^{r+1}}\bigg( \|\nabla\mathcal{P}_h\mathbf{\dot v}(t)\|_{L^{\infty}}+\|\nabla ( \mathbf{\dot W}(t)-\mathcal{P}_h\mathbf{\dot v}(t))\|_{L^{\infty}}+\|\nabla\mathcal{P}_h\mathbf{v}(t)\|_{L^{\infty}}&&\\
&+\|\nabla ( \mathbf{W}(t)-\mathcal{P}_h\mathbf{v}(t))\|_{L^{\infty}}\bigg)\times \| \nabla (\mathbf{\dot W}(t)-\mathcal{P}_h\mathbf{\dot  v}(t))\|_{L^2}&&\\
{}\lesssim & \: h^r \|\mathbf{v}(t)\|_{H^{r+1}}\bigg( \|\nabla\mathbf{\dot v}(t)\|_{L^{\infty}}+\|\nabla ( \mathbf{\dot v}(t)-\mathcal{P}_h\mathbf{\dot v}(t))\|_{L^{\infty}}+\|\nabla\mathbf{v}(t)\|_{L^{\infty}}+\|\nabla ( \mathbf{v}(t)-\mathcal{P}_h\mathbf{v}(t))\|_{L^{\infty}}\bigg)&&\\
{}& \times \|\nabla (\mathbf{\dot W}(t)-\mathcal{P}_h\mathbf{\dot v}(t))\|_{L^2}+ h^{r-\frac{d}{2}}\|\mathbf{v}(t)\|_{H^{r+1}}\|\nabla (\mathbf{ W}(t)-\mathcal{P}_h\mathbf{ v}(t) )\|_{L^2}\|\nabla (\mathbf{\dot W}(t)-\mathcal{P}_h\mathbf{\dot v}(t) )\|_{L^2}&&\\
{}& +h^{r-\frac{d}{2}}\|\mathbf{v}(t)\|_{H^{r+1}}\|\nabla (\mathbf{\dot W}(t)-\mathcal{P}_h\mathbf{\dot v}(t) )\|_{L^2}^2 &&\\
{}\lesssim &\: h^r \|\mathbf{v}(t)\|_{H^{r+1}}\|\nabla (\mathbf{\dot W}(t)-\mathcal{P}_h\mathbf{\dot v}(t))\|_{L^2}&&\\
&+ h^{r-\frac{d}{2}}\|\mathbf{v}(t)\|_{H^{r+1}}\|\nabla (\mathbf{ W}(t)-\mathcal{P}_h\mathbf{ v}(t) )\|_{L^2}\|\nabla (\mathbf{\dot W}(t)-\mathcal{P}_h\mathbf{\dot v}(t))\|_{L^2}&&\\
{}& +h^{r-\frac{d}{2}}\|\mathbf{v}(t)\|_{H^{r+1}}\|\nabla (\mathbf{\dot W}(t)-\mathcal{P}_h\mathbf{\dot v}(t))\|_{L^2}^2.
\end{flalign*}
Combining the estimates for $T_1$, $T_2$ and $T_3$, we have 
\begin{align}\label{T123}
&\int_0^1 \tilde{a}(\mathcal{P}_h\mathbf{v}(t)\mathrm{e}^{\gamma t}+\tau (\mathbf{W}(t)-\mathcal{P}_h\mathbf{v}(t))\mathrm{e}^{\gamma t};\mathbf{\dot W}(t)-\mathcal{P}_h\mathbf{\dot v}(t), \bm{\varphi}) \dd \tau \nonumber\\
&\lesssim   h^r \left(\|\mathbf{v}(t) \|_{H^{r+1}}+\|\mathbf{\dot v}\|_{H^{r+1}}\right)\|\nabla (\mathbf{\dot W}-\mathcal{P}_h\mathbf{\dot v}(t))\|_{L^2}\\
&\quad + h^{r-\frac{d}{2}}\|\mathbf{v}\|_{H^{r+1}}\|\nabla (\mathbf{ W}(t)-\mathcal{P}_h\mathbf{ v}(t) )\|_{L^2}\|\nabla (\mathbf{\dot W}(t)-\mathcal{P}_h\mathbf{\dot v}(t) )\|_{L^2}\nonumber\\
&\quad +h^{r-\frac{d}{2}}\|\mathbf{v}\|_{H^{r+1}}\|\nabla (\mathbf{\dot W}(t)-\mathcal{P}_h\mathbf{\dot v}(t) )\|_{L^2}^2\nonumber.
\end{align}
Applying the strong ellipticity condition (S2b) on the left-hand side of (\ref{T123}), we have 
\begin{flalign*}
M_1\|\nabla (\mathbf{\dot W}(t)- \mathcal{P}_h\mathbf{\dot v}(t))\|_{L^2}^2 \lesssim & \: h^r \left(\|\mathbf{v}(t) \|_{H^{r+1}}+\|\mathbf{\dot v}\|_{H^{r+1}}\right)\|\nabla (\mathbf{\dot W}-\mathcal{P}_h\mathbf{\dot v}(t))\|_{L^2}&&\\
{}& + h^{r-\frac{d}{2}}\|\mathbf{v}\|_{H^{r+1}}\|\nabla (\mathbf{ W}(t)-\mathcal{P}_h\mathbf{ v}(t) )\|_{L^2}\|\nabla (\mathbf{\dot W}(t)-\mathcal{P}_h\mathbf{\dot v}(t) )\|_{L^2}&&\\
{}& +h^{r-\frac{d}{2}}\|\mathbf{v}\|_{H^{r+1}}\|\nabla (\mathbf{\dot W}(t)-\mathcal{P}_h\mathbf{\dot v}(t) )\|_{L^2}^2.
\end{flalign*}
Dividing by $\|\nabla (\mathbf{\dot W}(t)-\mathcal{P}_h\mathbf{\dot v}(t) )\|_{L^2}$ on both sides yields \begin{flalign*}
M_1\|\nabla (\mathbf{\dot W}(t)- \mathcal{P}_h\mathbf{\dot v}(t))\|_{L^2} \lesssim & \: h^r \left(\|\mathbf{v}(t) \|_{H^{r+1}}+\|\mathbf{\dot v}\|_{H^{r+1}}\right)+ h^{r-\frac{d}{2}}\|\mathbf{v}\|_{H^{r+1}}\|\nabla (\mathbf{ W}(t)-\mathcal{P}_h\mathbf{ v}(t) )\|_{L^2}&&\\
{}& +h^{r-\frac{d}{2}}\|\mathbf{v}\|_{H^{r+1}}\|\nabla (\mathbf{\dot W}(t)-\mathcal{P}_h\mathbf{\dot v}(t) )\|_{L^2}.
\end{flalign*}
Since $r>\frac{d}{2}$, we can choose $h$ sufficiently small such that the last term on the right-hand side can be absorbed into the term on the left-hand side. This yields
\begin{equation}
\|\nabla(\mathbf{\dot W}(t)-\mathcal{P}_h\mathbf{\dot v}(t))\|_{L^2} \lesssim \:h^r \left( \|\mathbf{v}(t)\|_{H^{r+1}}+\|\mathbf{\dot v}(t)\|_{H^{r+1}}\right).
\end{equation}
Again, by the approximation property of $\mathcal{P}_h$, we have, for each $t\in [0,T]$, 
\begin{equation}
\|\nabla (\mathbf{\dot v}(t)-\mathcal{P}_h\mathbf{\dot v}(t))\|_{L^2} \lesssim \:h^r\|\mathbf{\dot v}(t)\|_{H^{r+1}}, \quad  \frac{d}{2}<r\leq \min(p,m-1).
\end{equation}
It follows from the triangle inequality that, for each $t\in [0,T]$,
\begin{equation}
\|\nabla\mathbf{\dot W}(t)-\nabla \mathbf{\dot v}(t)\|_{L^2}\lesssim \: h^r \left(\|\mathbf{v}(t)\|_{H^{r+1}}+\|\mathbf{\dot v}(t)\|_{H^{r+1}}\right), \quad  \frac{d}{2}<r\leq \min(p,m-1).
\end{equation}
By a similar argument as in the previous section, we can show that there exists a constant $c_1>0$ such that  
\begin{equation}
\|\nabla \mathbf{\dot W}(t)\|_{L^{\infty}}\leq c_1.
\end{equation}
\subsubsection{\mathinhead{L^2}{L2} bounds on \mathinhead{(\mathbf{v}-\mathbf{W})}{vminusW}, \mathinhead{(\mathbf{\dot{v}}-\mathbf{\dot{W}})}{dotvminusdotW} and \mathinhead{(\mathbf{\ddot {v}}-\mathbf{\ddot {W}})}{ddotvminusddotW}}
It was proved by Dobrowolski and Rannacher in \cite{DR} that for each $t\in [0,T]$,
\begin{equation}\label{DR}
\|\mathbf{v}(t)-\mathbf{W}(t)\|_{L^2}\leq C_r(\mathbf{v}) h^{r+1}, \quad  \frac{d}{2}<r\leq \min(p,m-1).
\end{equation}
We shall focus on proving the $L^2$ error bound of the time derivative using a duality argument in this section.
Consider the following boundary value problem: for a given $\mathbf{g}\in L^2$, solve $\bm{\psi}\in H_0^1$ such that 
\begin{equation}\label{adj_eq}
\tilde{a}(\mathbf{v}; \bm{\psi}, \bm{\phi})=(\mathbf{g}, \bm{\phi})_{L^2} \quad \mbox{ for all } \bm{\phi}\in H_0^1,
\end{equation}
where $\mathbf{v}$ is the solution of (\ref{setup2})--(\ref{ic2}) and
\begin{equation}\label{defoa}
\tilde{a}(\mathbf{v};\bm{\psi}, \bm{\phi})=\sum_{i,\alpha, j, \beta=1}^d (A_{i\alpha j\beta}(\nabla \mathbf{v})\partial_{\beta}\bm{\psi}, \partial_{\alpha} \bm{\phi}_i)_{L^2}.
\end{equation}
Since $A_{i\alpha j \beta}(\nabla \mathbf{v})\in W^{1,\infty}$ provided that $A_{i\alpha j \beta}$ is sufficiently smooth and $\nabla v\in C^{2,\alpha}$ (cf. Remark \ref{reg_v} ), the adjoint problem (\ref{adj_eq}) has a unique solution which satisfies the following elliptic regularity conditions, cf.  Theorem 1.1 and Theorem 2.6 of Chapter 8 in \cite{CW2},
\begin{equation}\label{elliptic_reg2}
\|\bm{\psi}\|_{H^{2}}\leq \hat{c}\left( \|\bm{\psi}\|_{L^2}+\|\mathbf{g}\|_{L^2}\right)
\end{equation}
for some positive constant $\hat{c}$. Taking $\bm{\phi}=\bm\psi\in H_0^1$ in (\ref{defoa}) and applying the coercive condition (S2b), we have 
\begin{equation}\label{coer}
M_1\|\nabla\bm{\psi}\|_{L^2}^2 \leq \|\mathbf{g}\|_{L^2}\|\bm{\psi}\|_{L^2}.
\end{equation}
Applying Poincar\'{e}'s inequality in (\ref{coer}), we deduce that 
$$\|\bm{\psi}\|_{L^2}\leq M_1^{-1} C_{\mathrm{poin}}\|\mathbf{g}\|_{L^2}.$$ Thus
\begin{equation}\label{elliptic_reg}
\|\bm{\psi}\|_{H^{2}}\leq c\|\mathbf{g}\|_{L^2},
\end{equation}
for some positive constant $c$.
The corresponding discrete problem is formulated as: find $\bm{\psi}_h\in\mathcal{V}_h$ such that 
\begin{equation}\label{adjoint_discrete}
\tilde{a}(\mathbf{v}; \bm{\psi}_h, \bm{\phi})=(\mathbf{g}, \bm{\phi})_{L^2}, \quad \mbox{ for all } \bm{\phi}\in \mathcal{V}_h.
\end{equation}
It is known that we have, cf., e.g., \cite{DR},
\begin{equation}\label{appro_prop}
\|\bm{\psi}-\bm{\psi}_h\|_{L^2}+h\|\bm{\psi}-\bm{\psi}_h\|_{H^1} \leq C h^{r+1} \|\bm{\psi}\|_{H^{r+1}},
\end{equation}
for some constant $C.$
Let $\mathbf{g}=\mathbf{\dot{v}}-\mathbf{\dot{W}}$, then (\ref{adj_eq}) becomes 
\begin{equation}\label{new_adj}
\tilde{a}(\mathbf{v},\bm{\psi}, \bm{\phi})=(\mathbf{\dot{v}}-\mathbf{\dot{W}}, \bm{\phi})_{L^2}.
\end{equation}
Plugging $\bm{\phi}=\mathbf{\dot{v}}-\mathbf{\dot{W}}
$ into (\ref{new_adj}), we obtain 
\begin{equation}
\|\mathbf{\dot{v}}-\mathbf{\dot{W}}\|_{L^2}^2=\tilde{a}(\mathbf{v}; \bm{\psi}, \mathbf{\dot{v}}-\mathbf{\dot{W}}).
\end{equation}
Using (S2a) and the definition of the elliptic projection, we have,
\begin{align}\label{estimate_l2}
\tilde{a}(\mathbf{v};\bm{\psi}, \mathbf{\dot{v}}-\mathbf{\dot{W}})=&\:\tilde{a}(\mathbf{v};\bm{\psi}-\bm{\psi}_h, \mathbf{\dot{v}}-\mathbf{\dot{W}})+ \tilde{a}(\mathbf{v};\bm{\psi}_h, \mathbf{\dot{v}}-\mathbf{\dot{W}})\nonumber\\
{}=& \:\tilde{a}(\mathbf{v};\bm{\psi}-\bm{\psi}_h, \mathbf{\dot{v}}-\mathbf{\dot{W}})+\tilde{a}(\mathbf{v}; \mathbf{\dot{v}}-\mathbf{\dot{W}},\bm{\psi}_h)\nonumber\\
{}=&\: \tilde{a}(\mathbf{v};\bm{\psi}-\bm{\psi}_h, \mathbf{\dot{v}}-\mathbf{\dot{W}})+\tilde{a}(\mathbf{W}; \mathbf{\dot{W}},\bm{\psi}_h)-\tilde{a}(\mathbf{v}; \mathbf{\dot{W}},\bm{\psi}_h).
\end{align}
By (iii,a), (\ref{elliptic_reg}) and (\ref{appro_prop}), we have 
\begin{align*}
|\tilde{a}(\mathbf{v};\bm{\psi}-\bm{\psi}_h, \mathbf{\dot{v}}-\mathbf{\dot{W}})|\leq&\: K_{\delta}\|\nabla(\mathbf{\dot{v}}-\mathbf{\dot{W}})\|_{L^2}\|\nabla(\bm{\psi}-\bm{\psi_h})\|_{L^2} \\
{}\leq &\:K_{\delta} C_r(\mathbf{v})h^r \cdot C h\|\bm{\psi}\|_{H^2} \quad \mbox{(by (iii,a) and (\ref{appro_prop}))}\\
{}\leq &\:K_{\delta} cC_r(\mathbf{v})Ch^{r+1}\|\mathbf{\dot{v}}-\mathbf{\dot{W}}\|_{L^2} \quad \mbox{(by (\ref{elliptic_reg}))}.
\end{align*}
For the remaining terms in (\ref{estimate_l2}), we observe that 
\begin{align*}
\vert\tilde{a}(\mathbf{W}; \mathbf{\dot{W}},\bm{\psi}_h)-\tilde{a}(\mathbf{v}; \mathbf{\dot{W}},\bm{\psi}_h)\vert\leq &\:\vert \tilde{a}(\mathbf{W}; \mathbf{\dot{W}},\bm{\psi}_h-\bm{\psi})-\tilde{a}(\mathbf{v}; \mathbf{\dot{W}},\bm{\psi}_h-\bm{\psi})\vert\\
{}&+ \vert\tilde{a}(\mathbf{W}; \mathbf{\dot{W}}-\mathbf{\dot{v}},\bm{\psi})-\tilde{a}(\mathbf{v}; \mathbf{\dot{W}}-\mathbf{\dot{v}},\bm{\psi})\vert\\
{}&+\vert\tilde{a}(\mathbf{W}; \mathbf{\dot{v}},\bm{\psi})-\tilde{a}(\mathbf{v}; \mathbf{\dot{v}},\bm{\psi})\vert\\
{} := &\: T_4+T_5+T_6.
\end{align*}
By Lipschitz continuity of $A_{i\alpha j \beta}$, we have 
\begin{align*}
T_4\leq &\:L_{\delta}\|\nabla (\mathbf{W}-\mathbf{v})\|_{L^2}\|\nabla \mathbf{\dot{W}}\|_{L^{\infty}}\|\nabla(\bm{\psi}_h-\bm{\psi})\|_{L^2}\\
{}\leq  &\: C_r(\mathbf{v})h^r c_1  \|\nabla(\bm{\psi}_h-\bm{\psi})\|_{L^2}\quad \mbox{ ( by (iii,a) and (iii,c))}\\
{}\leq  &\: C_r(\mathbf{v})c_1h^r\cdot Ch  \|\bm{\psi}\|_{H^2}\quad \mbox{ (by (\ref{appro_prop}))}\\
{}\leq & \: cC C_r(\mathbf{v})c_1h^{r+1}\|\mathbf{\dot{W}}-\mathbf{\dot{v}}\|_{L^2} \quad\mbox{ (by (\ref{elliptic_reg}))}.
\end{align*}
Similarly, we have
\begin{align*}
T_5\leq & \:L_{\delta} \|\nabla(\mathbf{W}-\mathbf{v})\|_{L^{\infty}}\|\nabla(\mathbf{\dot{W}}-\mathbf{\dot{v}})\|_{L^2} \|\nabla\bm{\psi}\|_{L^2}.\\
\end{align*}
Following the analysis in \cite{RS} and Chapter 8 of \cite{BS}, it can be shown that 
\begin{equation}\label{RS}
\|\mathbf{v}-\mathbf{W}\|_{W^{1,\infty}}\leq c(\mathbf{v})h^r,
\end{equation}
where $c(\mathbf{v})$ is a positive constant depending on the exact solution $\mathbf{v}.$
Therefore, we can bound $T_5$ by
\begin{equation}
T_5 \leq  L_{\delta} c(\mathbf{v}) h^r C_{r}(\mathbf{v}) h^r \|\bm{\psi}\|_{H^2}\leq   L_{\delta} c(\mathbf{v})  C_{r}(\mathbf{v}) c h^{r+1}\|\mathbf{\dot{W}}-\mathbf{\dot{v}}\|_{L^2},
\end{equation}
for any $r\geq 1$ provided that $h$ is sufficiently small. We bound $T_6$ by 
\begin{align*}
T_6 = & \:|\tilde{a}(\mathbf{W}; \dot{\mathbf{v}}, \bm{\psi})-\tilde{a}(\mathbf{v}; \dot{\mathbf{v}}, \bm{\psi})|\\
{} \leq & \:\bigg\vert \sum_{i,\alpha, j,\beta, k,\gamma=1}^d \bigg( \partial_{\gamma}(\mathbf{W}_k-\mathbf{v}_k)\partial_{\beta}\dot{\mathbf{v}}_j\frac{\partial A_{i\alpha j \beta}}{\partial \eta_{k\gamma}}(\nabla\mathbf{v}), \partial_{\alpha} \bm{\psi}_i\bigg)\bigg\vert\\
&+\bigg\vert\sum_{i,\alpha, j,\beta, k,\gamma, l, \delta=1}^d \bigg( \partial_{\delta}(\mathbf{W}_l-\mathbf{v}_l)\partial_{\gamma}(\mathbf{W}_k-\mathbf{v}_k)\partial_{\beta}\dot{\mathbf{v}}_j\int_0^1 \frac{\partial^2 A_{i\alpha j\beta}}{\partial \eta_{l\delta}\partial\eta_{k\gamma}}(\nabla\mathbf{v}+\tau\nabla(\mathbf{W}-\mathbf{v}))\dd \tau, \partial_{\alpha} \bm{\psi}_i \bigg)\bigg\vert\\
{}:= &\: b(\mathbf{W},\mathbf{v}; \mathbf{\dot{v}},\bm{\psi})+ d(\mathbf{W},\mathbf{v}; \mathbf{\dot{v}},\bm{\psi}).
\end{align*}
To ensure that  $\nabla \mathbf{W}\in \mathcal{Z}_{\delta}$, we take $h$ sufficiently small. By the convexity of $\mathcal{Z}_{\delta}$ and $\mathcal{M}_{\delta}$, we know that $\nabla\mathbf{v}(x)+\tau\nabla(\mathbf{W}(x)-\mathbf{v}(x))\in \mathcal{M}_{\delta}$ for $x\in \bar{\Omega}.$
Since $A_{i\alpha j \beta}$ is sufficiently smooth (in particular twice continuously differentiable on $\mathcal{M}_{\delta}$), we have
\begin{align*}
 d(\mathbf{W},\mathbf{v}; \dot{\mathbf{v}},\bm{\psi})\leq & \:C_A\|\nabla(\mathbf{W}-\mathbf{v})\|_{L^{\infty}}^2 \|\nabla \dot{\mathbf{v}}\|_{L^2}\|\nabla \bm{\psi}\|_{L^2}
\\
{}\leq &\: C_A c(\mathbf{v})^2 h^{2r}\|\nabla\mathbf{\dot{ v}}\|_{L^2}\|\nabla\bm{\psi}\|_{L^2} \quad \mbox{ (by (\ref{RS}) )}\\
{}\leq & \:C_A c(\mathbf{v})^2 h^{2r}\|\nabla\mathbf{\dot{v}}\|_{L^2}\|\bm{\psi}\|_{H^2}\\
{}\leq & \:C_A c(\mathbf{v})^2c h^{2r}\|\nabla\mathbf{\dot{v}}\|_{L^2}\|\mathbf{\dot{W}}-\mathbf{\dot{v}}\|_{L^2}\quad\mbox{ (by (\ref{elliptic_reg}) )}\\
{}\leq &\: C h^{r+1}\|\dot{\mathbf{W}}-\dot{\mathbf{v}}\|_{L^2}, \quad \mbox{ for } r\geq 1.
\end{align*}
For the estimation of $b(\mathbf{W},\mathbf{v}; \dot{\mathbf{v}},\bm{\psi})$, we apply integration by parts and the fact that $\mathcal{V}_h\subset H_0^1$ to obtain
\begin{align*}
b(\mathbf{W},\mathbf{v}; \dot{\mathbf{v}},\bm{\psi})= & \:\bigg\vert \sum_{i,\alpha, j,\beta, k,\gamma=1}^d \bigg( \partial_{\gamma}\big[\partial_{\beta}\dot{\mathbf{v}}_j\frac{\partial A_{i\alpha j \beta}}{\partial \eta_{k\gamma}}(\nabla\mathbf{v}), \partial_{\alpha} \bm{\psi}_i\big](\mathbf{W}_k-\mathbf{v}_k)\bigg)\bigg\vert\\
{}\leq & \:C\|\mathbf{W}-\mathbf{v}\|_{L^2} \|\mathbf{\dot{v}}\|_{W^{2,\infty}}\|\mathbf{v}\|_{W^{2,\infty}}\|\bm{\psi}\|_{H^2}\\
{}\leq & \:\tilde{C} h^{r+1}\|\mathbf{\dot{W}}-\mathbf{\dot{v}}\|_{L^2} \quad\mbox{ (by (\ref{elliptic_reg}) and (\ref{DR}) )}. 
\end{align*}
Combining the above estimates for $T_4, T_5$ and $T_6$, we have 
\begin{equation}
\|\mathbf{\dot{W}}-\mathbf{\dot{v}}\|_{L^2}\leq \tilde{C}_r(\mathbf{v}) h^{r+1},
\end{equation}
for some positive constant $\tilde{C}_r(\mathbf{v}).$

By a similar argument, we can easily show that $\mathbf{\dot W}(t)$ is differentiable with respect to $t$ and a similar $L^2$ error estimate for $ \mathbf{\ddot{W}}-\mathbf{\ddot{v}}$. The proof of this estimate can be found in \cite{Ma1} and \cite{Ma2}. We omit the details here.

\section*{Acknowledgement}
\begin{acknowledgement}
The author would like to thank Prof. Endre S\"{u}li for his helpful discussions and suggestions.
\end{acknowledgement}

\bibliographystyle{siam}
\bibliography{refs}
\end{document}